\documentclass[11pt]{article} % use larger type; default would be 10pt
\usepackage[T1]{fontenc}
\usepackage[english]{babel}

\usepackage[utf8]{inputenc} % set input encoding (not needed with XeLaTeX)

\usepackage{geometry} % to change the page dimensions
\usepackage{graphicx} % support the \includegraphics command and options

\usepackage{booktabs} % for much better looking tables
\usepackage{array} % for better arrays (eg matrices) in maths
\usepackage{paralist} % very flexible & customisable lists (eg. enumerate/itemize, etc.)
\usepackage{verbatim} % adds environment for commenting out blocks of text & for better verbatim
\usepackage{subfig} % make it possible to include more than one captioned figure/table in a single float

\usepackage{csquotes}
\usepackage{graphicx} % support the \includegraphics command and options
\usepackage{amssymb}
\usepackage{amsmath, amsfonts, bbm, dsfont}
\usepackage[all]{xy}
\usepackage{MnSymbol}

\usepackage{fancyhdr} % This should be set AFTER setting up the page geometry
\usepackage{sectsty}
\allsectionsfont{\sffamily\mdseries\upshape} % (See the fntguide.pdf for font help)
\usepackage[nottoc,notlof,notlot]{tocbibind} % Put the bibliography in the ToC
\usepackage[titles,subfigure]{tocloft} % Alter the style of the Table of Contents

\usepackage[
backend=biber,
style=alphabetic
]{biblatex}

\usepackage{amsthm}
\usepackage{tikz-cd}

\tikzcdset{row sep/normal= 0.8 cm,
column sep/normal= 0.8 cm}

\newtheorem{theo}{Theorem}[section]
\newtheorem{prop}[theo]{Proposition}
\newtheorem{lemma}[theo]{Lemma}
\newtheorem{coro}[theo]{Corollary}

\newtheorem{fact}[theo]{Fact}
\newtheorem{claim}[theo]{Claim}
\newtheorem{claim*}{Claim}

\theoremstyle{definition}
\newtheorem{defi}[theo]{Definition}
\newtheorem{rem}[theo]{Remark}
\newtheorem{notation}[theo]{Notation}

\newtheorem{contre-ex}[theo]{Contre-example}

\DeclareMathOperator{\flex}{\textbf{\textasciicircum}}

\usepackage{hyperref}
\hypersetup{
    colorlinks=true,
    linkcolor=black,
    filecolor=black,      
    urlcolor=black,
    citecolor=black,
}

\title{The group configuration theorem for generically stable types}
\author{Paul $\mathrm{Wang}^1$}
\date{}
\begin{document}
\setcounter{tocdepth}{2}
\maketitle

  %\tableofcontents

\footnotetext[1]{partially funded by ANR GeoMod  (AAPG2019, ANR-DFG)}
 \footnotetext[2]{\copyright \,  2023. This manuscript version is made available under the CC-BY-NC-ND 4.0 license http://creativecommons.org/licenses/by-nc-nd/4.0/}

\section*{Introduction}

In his thesis \cite{hrushovski-thesis}, Ehud Hrushovski proved a group configuration theorem, building a type-definable group from combinatorial data, in a stable setting. The aim of this paper is to generalize the theorem, using only hypotheses on the type of the configuration, without assuming tameness of the theory.

First, we shall introduce generically stable types, and state some of their known properties. Then, we will define some notions of genericity in definable groups and definable homogeneous spaces, and show a couple of results regarding groups with generically stable generics. Having done that, we shall state and prove a group configuration theorem (Theorem \ref{theo_config_groupe_faible}) for generically stable types. The proofs will be similar to the stable case, although a bit trickier. In passing, we also write down a uniqueness result (Proposition \ref{prop_recovering_group_from_config}), recovering a group with generically stable generics from its configuration, up to some notion of equivalence. That result is not new, and has actually been improved substantially. See for instance \cite[Theorem 2.15]{Stabilizers-NTP2}.

From now on, we fix a complete theory $T$,  in a language $\mathcal{L}$, and work inside $T^{eq}$ to ensure elimination of imaginaries. We let $acl$, resp. $dcl$, denote the algebraic closure, resp. definable closure, in $T^{eq}$.  Moreover, we let $\mathbb{U}$ denote a very saturated and strongly homogeneous model of $T$. A subset $A$ of a model $M$ is called \emph{small}, with respect to $M$, if $M$ is $|A|^+$-saturated and $|A|^+$-strongly homogeneous. Note that we might consider models $M \subset \mathbb{U}$ and sets $A \subset M$ such that $A$ is small with respect to $M$, and $M$ itself is small with respect to $\mathbb{U}$. By default, the sets of parameters we consider are small with respect to $\mathbb{U}$. If $a,b$ are small tuples, we may write $ab$ or $ a \flex b$ for the concatenation.

\section*{Acknowledgements}

This work is a continuation of my master's thesis, under the supervision of Silvain Rideau-Kikuchi. I would like to thank him for his guidance. I would also like to thank the anonymous referee for many helpful comments and suggestions.

\begin{section}{Generically stable types}

\setcounter{subsection}{0}

\begin{subsection}{Forking, invariant and definable types}

\begin{defi}\label{defi_deviation_et_division}

Let $A \subseteq \mathbb{U}$ be a set of parameters and $\phi(x,y)$ be a formula over $A$. Let $b$ be a tuple.

\begin{enumerate}
    \item The formula $\phi(x,b)$ divides over $A$ if there is an $A$-indiscernible sequence $(b_i)_{i < \omega}$, with $b_0=b$, such that the partial type $\lbrace \phi(x,b_i) \, | \, i < \omega \rbrace$ is inconsistent.
    \item The formula $\phi(x,b)$ forks over $A$ if $\phi(x,b)$ implies a finite disjunction of formulas, possibly with additional parameters, all of which divide over $A$.
\end{enumerate}

For any natural number $k$, a partial type $\pi$ is $k$-inconsistent if, for any choice of pairwise non-equivalent formulas $\phi_1, ..., \phi_k$ in $\pi$, the conjunction  $\bigwedge\limits_i\phi_i$ is not satisfiable. 

Then, by compactness and indiscernibility, one may replace "inconsistent" with "$k$-inconsistent for some $k$" in the above definition of dividing.

\end{defi}

\begin{defi}

\begin{enumerate}
    \item A partial type $\pi(x)$ divides (resp. forks) over a set $A$ if there exists a formula $\phi(x)$ (possibly with parameters outside of $A$) such that  $\pi(x) \models \phi(x)$ and $\phi(x)$ divides (resp. forks) over $A$.
    
    \item Let $a,b$ be tuples, and $C$ be a set. The tuple $a$ is independent from $b$ over $C$, which we denote $a\downfree_C b$, if $tp(a/ Cb)$ does not fork over $C$.
    
    \item Let $p \in S(A)$. The type $p$ is extensible (resp. stationary) if, for any $B \supseteq A$, there exists a (resp. a unique) $q \in S(B)$ such that $q|_A=p$ and $q$ does not fork over $A$. Given a stationary type $p \in S(A)$ and $B \supseteq A$, we let $p|_B \in S(B)$ denote its unique nonforking extension.
    
    %\item Let $p \in S(A)$, $q \in S(B)$ be stationary types. They are parallel if the unique nonforking extension of $p$ over $AB$ is equal to  the unique nonforking extension of $q$ over $AB$. In other words, $p$ and $q$ are parallel if $p|_{AB} = q|_{AB}$.
\end{enumerate}

\end{defi}

%\begin{rem}A partial type closed under conjunction divides (resp. forks) over a set $A$ if and only if it contains a formula which divides (resp. forks) over $A$.\end{rem}

\begin{defi}
Let $(a_i)_{i \in I}$ be a family of elements. Let $A$ be a set of parameters. We say that $(a_i)$ is an independent family over $A$ if, for all $i \in I$, we have $a_i \downfree_{A} (a_j)_{j \in I, \, j \neq i} .$

\end{defi}

%\begin{defi}A partial type $\pi$ is  finitely satisfiable in a set $A$ if, for any natural number $n$, for any choice of formulas $\phi_1,..., \phi_n$ in $\pi$, there exists $a \in A$ such that $\models \bigwedge\limits_{i=1}^n\phi_i(a)$.\end{defi}

%\begin{prop} A partial type finitely satisfiable in $A$ does not fork over $A$.\end{prop}

\begin{notation}
To simplify notations, if $A$ is a small set of parameters, we write $|A|$ instead of $|A|+|\mathcal{L}|+\aleph_0$.
\end{notation}

\begin{defi}\label{defi_definable_types}

Let $M$ be a model of $T$, let $p \in S(M)$, and $A \subseteq M$. 

\begin{enumerate}
    \item We say that $p$ is $A$-definable if, for all formulas without parameters $\phi(x,y)$, there exists a formula $d_p x \phi(x,y)$ with parameters in $A$ such that, for all $b \in M$, we have $M \models d_p x \, \phi(x, b)$ if and only if $\phi(x,b) \in p(x)$. We say that $d_p$ is ``the'' defining scheme for $p$. Indeed, since $M$ is a model, the defining scheme is unique up to equivalence.
    
    \item We say that $p$ is definable if it is $M$-definable.
    
    %\item If $d$ is a defining scheme, we let $d(B)$ denote the partial type over $B$ defined by "$\phi(x,b) \in d(B)$ if and only if $\models d x \, \phi(x, b)$", for any set $B$.
    
    \item If $p$ is definable, the canonical basis of $p$ is the smallest dcl-closed set $A \subseteq M$ such that $p$ is $A$-definable. By elimination of imaginaries, this set is well-defined.

    \item In the case where $M$ is  $|A|^+$-saturated, we say that $p$ is $A$-invariant if, for any formula without parameters $\phi(x,y)$, for all $b_1, b_2 \in M$, if $b_1 \equiv_A b_2$, then $\phi(x, b_1) \in p(x)$ if and only if $\phi(x, b_2) \in p(x)$. In other words, the formula $\phi(x,b)$ being in the type $p$ depends only on the type of $b$ over $A$.

    \item If $B \supseteq A$ is a set of parameters (not necessarily a model), and $q \in S(B)$, we say that $q$ is $A$-invariant (resp. $A$-definable, resp. definable) if it admits an $A$-invariant (resp. $A$-definable, resp. $B$-definable) extension $q_1$ to an $|A|^+$-saturated  model $N \supseteq B$.

\end{enumerate}

\end{defi}

%\begin{rem}In practice, we consider $A$-invariant types only when $M$ is an $|A|^+$-saturated and $|A|^+$-strongly homogeneous model. Then, being $A$-invariant is equivalent to being fixed by all the automorphisms of the group $Aut(M / A)$.\end{rem}

\begin{fact}\label{fact_unique_invariant_extension}\emph{(see \cite[Proposition 1.9]{Pil-Stab} and \cite{Sim-Book}(Section 2.2, discussion before Lemma 2.18))}
Let $A$ be a small subset of a model $M$. Let $p \in S(M)$ be $A$-invariant. Then, for any model $N \supseteq M$, the type $p$ has a unique extension $p|_N \in S(N)$ which is $A$-invariant

Moreover, if $p$ is $A$-definable, then $p|_N$ is $A$-definable, using the same defining scheme as $p$, and the whole conclusion holds even if $A$ is not small with respect to $M$.
\end{fact}

\begin{rem}\label{rem_invariant_types_as_processes}
Thanks to this fact, if $p \in S(M)$ is $A$-invariant, \emph{where $M$ is $|A|^+$-saturated}, and if $B \supseteq A$, we can write $p|_B \in S(B)$ for the restriction to $B$ of $p|_N$, where $N$ is a model containing $M B$. By uniqueness, the type $p|_B$ is well-defined, for it does not depend on the choice of the model $N$.

Because of this fact, we consider it useful to view invariant types as families of types, or as processes which construct complete types in a coherent way, and to identify two invariant types if they admit a common invariant extension to a sufficiently saturated and sufficiently strongly homogeneous model. 

\end{rem}

\begin{defi}
Let $A$ be a small subset of a model $M$. Let $p, q \in S(M)$, where $p$ is $A$-invariant. Let us assume that $p$ is in the variable $x$, and $q$ in the variable $y$. We define the tensor product $p \otimes q \in S(M)$ as follows:

If $\phi(x,y, z)$ is a formula without parameters, and if $c \in M$, then $\phi(x,y,c) \in p\otimes q$ if and only if, for some (equivalently, for all) element   $b \in M$ realizing $q|_{Ac}$, we have $\phi(x, b,c) \in p|_{Abc}$. 
\end{defi}

This tensor product is well-defined, and is a complete type over $M$. Indeed, we can check that the realizations of $p \otimes q$ are exactly the tuples of the form $ab$, where $b$ realizes $q$, and $a$ realizes $p|_{Mb}$. Note that if $q$ is also $A$-invariant, then $p \otimes q$ is $A$-invariant. If $p$ and $q$ are $A$-definable, then so is $p\otimes q$.

\begin{defi}
Let $A$ be a small subset of a model $M$. Let $p \in S(M)$ be an  $A$-invariant type. Let $\alpha$ be an ordinal. A sequence $(a_i)_{i < \alpha}$ is a \emph{Morley sequence} of $p$ over $A$ if, for all $i < \alpha$, the element $a_i$ realizes the type $p|_{A \cup (a_j)_{j < i}}$.
\end{defi}

%\begin{defi}
%Let $d_p$ be a defining scheme of a weakly definable type $p \in S(A)$. We say that $d_p$ is a good defining scheme if:

%\begin{enumerate}
    %\item For all $n$, for all finite family $\phi_1(x,y),..., \phi_n(x,y)$ of partitioned formulas without parameters, we have $\models \forall y [\bigwedge\limits_i d_p x \, \phi_i(x, y) \rightarrow \exists x \, \bigwedge\limits_i \phi_i(x, y)]$.
    
    %\item For all partitioned formula without parameters $\phi(x,y)$, we have
    
    %\noindent $\models \forall y [d_p x \, (\neg \phi (x,y)) \leftrightarrow \neg (d_p x \, \phi(x,y))]$.
%\end{enumerate}
%\end{defi}

%\begin{prop}\label{prop_bons_schémas}
%Let $d$ be a defining scheme, with parameters in $A$. The following are equivalent:

%\begin{enumerate}\item The scheme $d$ is a good defining scheme.    \item For all set of parameters $B$, the partial type $d(B)$ is a complete type over $B$.   \item There exists a model $M \supseteq A$ such that $d(M)$ is a complete type over $M$.\end{enumerate}\end{prop}

%\begin{defi}A weakly definable type is called \emph{definable} if it admits a good defining scheme.\end{defi}

For the following facts on forking, see for instance \cite{casanovas_2011}, \cite[Propositions 4.8 and 4.9]{Hossain-ExtBasHen}, and \cite[Section 7.1]{TenZie}.

\newpage

\begin{fact}\label{fact_forking}

Let $A \subseteq B \subseteq C$ be small sets of parameters, and $a, b, c$ be small tuples. 

\begin{enumerate}
\item Let $p \in S(C)$. Assume that $p$ does not fork over $A$. Then $p$ does not fork over $B$ and $p|_B$ does not fork over $A$.

\item Let $p \in S(B)$ be a type which does not fork over $A$. Then, there exists $q \in S(C)$ extending $p$ such that $q$ does not fork over $A$.

\item Assume that $a \downfree_A B$ and $b \downfree_{Aa} B$. Then $ab \downfree_A B$.

\item Assume $acl(A) \subseteq B$. Then:  \begin{enumerate}
    \item $tp(a/B)$ forks over $A$ if and only if $tp(a/B)$ forks over $acl(A)$.
    \item $a \downfree_A b$ if and only if $acl(Aa) \downfree_A acl(Ab)$.
\end{enumerate}

\item  Assume $a \in acl(Ab)\cap acl(Ac)$ and $b \downfree_A c$. Then $a \in acl(A)$. 

%TO DO: check that this is indeed true !

\item Let $M \supset A$ be an $|A|^+$-saturated and $|A|^+$-strongly homogeneous model. Let $p \in S(M)$ be an $A$-invariant type. Then $p$ does not fork over $A$.

\end{enumerate}
\end{fact}

%\begin{coro}\label{coro_un_type_definissable_ne_devie_pas}Let $M \supseteq A$ be a model, and $p \in S(M)$ an $A$-definable type. Then $p$ does not fork over $A$.\end{coro}

%\begin{prop}\label{deviation_et_dcl}Let $a$ be an element, let $A \subseteq B$ be sets of parameters. Let $c \in acl(aA)$. Assume that $tp(a / B)$ does not fork over $A$. Then $tp(c / B)$ does not fork over $A$.\end{prop}

%\begin{prop}\label{definissable_sur_acl_implique_non-deviant}Let $A$ be a set of parameters, $M$ a model containing $A$. Let $p \in S(M)$ be a type definable over $acl(A)$. Then $p$ does not fork over $A$.\end{prop}

\begin{lemma}\label{definissabilite_et_acl}
Let $a$ be an element, $C \subseteq D=acl(D)$, such that $tp(a/D)$ is definable over $C$. Then, for all $b \in acl(Ca)$, the type $tp(b/D)$ is definable over $acl(C)$. Similarly, if $b \in dcl(Ca)$, then $tp(b/D)$ is definable over $C$.
\end{lemma}

\begin{proof}
Let us prove the first point, the second one being easier. If $b \in acl(Ca)$, let $\phi(x,y)$ be a formula over $C$ such that $\models \phi(a, b)$ and $\models \forall x \, \exists^{\leq k} \, y \, \phi(x,y)$. By Definition \ref{defi_definable_types} (5), let $M \supset D$ be a sufficiently saturated model such that $p = tp(a/M)$ is $A$-definable. Let $\psi(y, z)$ be a formula without parameters. Let $q$ denote $tp(b/M)$.
Let us consider the following $C$-definable binary relation: $d_1 E d_2$ if and only if $\models d_p x [\forall y \,  \phi(x,y) \rightarrow [\psi(y,d_1)\leftrightarrow \psi(y,d_2)]]$.
Then, for $d_1, d_2 \in M$, if $d_1 E d_2$, then $\models \psi(b, d_1)\leftrightarrow \psi(b, d_2)$. Indeed, $a$ realizes $p|_{C d_1 d_2}$. 

\begin{claim}
The relation $E$ is a $C$-definable finite equivalence relation. 
\end{claim}
\begin{proof}
Reflexivity and symmetry are clear. Let us prove transitivity. Let $d_1 E d_2$ and $d_2 E d_3$. Let $\alpha$ realize $p|_{C d_1 d_2 d_3}$. Let $\beta $ be such that $\models \phi(\alpha, \beta)$. Then, we have $\models [\psi(\beta, d_1) \leftrightarrow \psi(\beta, d_2)] \wedge [\psi(\beta, d_2) \leftrightarrow \psi(\beta, d_3)]$. So $\models [\psi(\beta, d_1) \leftrightarrow \psi(\beta, d_3)]$. Since $\alpha \models p|_{C d_1 d_3}$, we have indeed $d_1 E d_3$.

Finally, we shall prove that $E$ has only finitely many classes. As $E$ is $C$-definable, and $C \subseteq M$, it is enough to prove that $E(M)$ has only finitely many classes. Since $a$ realizes $p|_M$, we know that, for $d_1, d_2 \in M$, we have $d_1 E d_2$ if and only if $ \models \forall y \,  (\phi(a,y) \rightarrow [\psi(y, d_1)\leftrightarrow \psi(y, d_2)])$. But $\models \exists^{\leq k} \, y \, \phi(a,y)$, therefore $E(M)$ has at most $2^k$ classes. So $E$ has at most $2^k$ classes, so it is a finite equivalence relation.
\end{proof}

Then, by elimination of imaginaries, let $c_1/E, \cdots, c_r/E$ be the codes of the classes modulo $E$, where the $c_i$ are in $M$. These codes are in $acl(C)$. To construct the definition $d_q y \psi(y,z)$, let $I \subseteq \lbrace 1, \cdots, r \rbrace$ be the set of indices $i$ such that $\models \psi(b, c_i)$. Using the property of the relation $E$, one can check that the formula $\bigvee_{i \in I} z E c_i$ is an appropriate definition of $q$ for the formula $\psi(y,z)$. Moreover, since the $c_i / E$ are in $acl(C)$, this formula is equivalent to a formula defined over $acl(C)$. 
\end{proof}

\end{subsection}

\begin{subsection}{General properties of generically stable types}

The definition of generically stable types below is from \cite{pillay-tanovic} (Definition 2.1). Most of the properties in this subsection come from \cite{garcia-onshuus-usvyatsov} (Appendix A), from \cite{pillay-tanovic} (Proposition 2.1) and \cite{adler-casanovas-pillay} (Fact 1.9, Lemma 2.1 and Theorem 2.2).

\begin{defi}\label{defi_gen_stable}

\begin{enumerate}
    \item Let $(a_i)_{i \in I}$ be a sequence of elements of the same sort. Let $B$ be a set. The mean of the types of the $a_i$ over $B$ is a partial (possibly complete) type, containing the formulas $\phi(x, b)$ over $B$ such that, for cofinitely many indices $i \in I$, we have $\models \phi(a_i,b)$.
    
    \item Let $A$ be a set, and $p \in S(M)$, where $M \supset A$ is a sufficiently saturated  model. The type $p$ is \emph{generically stable} over $A$ if $p$ is $A$-invariant and, for all ordinals $\alpha \geq \omega$, for all Morley sequences $(a_i)_{i < \alpha}$ of $p$ over $A$, the mean of the types of the $a_i$ over $M$ is a complete type over $M$.
\end{enumerate}

\end{defi}

\begin{rem}\label{rem_ops_alpha_denombrable}
The property for infinite ordinals in the definition above is equivalent to that for countably infinite ordinals. Indeed, a mean over an infinite index set is always a consistent partial type. If it is not complete, there exists a formula witnessing incompleteness. Then, countably many indices are enough to witness incompleteness of the mean for this formula. Thus, if the model $M$ is sufficiently saturated, it is enough to check the property for Morley sequences made of elements of $M$.
\end{rem}

\begin{prop}\label{proprietes_type_gen_stable}

Let $p \in S(M)$ be a complete type, generically stable over a small set $A \subset M$. Then:

\begin{enumerate}
    \item For any infinite Morley sequence $(a_i)_i$ of $p$ over $A$, the mean of the types of the $a_i$ is the type $p$ itself.
    
    \item For any $\phi(x,y)$ over $A$, there exists a natural number $n_{\phi}$ such that, for any infinite Morley sequence $(a_i)_i$ of $p$ over $A$, for any $b$, we have $\phi(x,b) \in p$ if and only if the set of indices $i$ such that $\models \neg \phi(a_i, b)$ contains at most $n_{\phi}$ elements.
    
    \item The type $p$ is definable over $A$.

    \item Any Morley sequence of $p$ over $A$ is an indiscernible set over $A$.
    
    \item If $B \subset M$ is a small set such that $p$ is $B$-invariant, then $p$ is generically  stable over $B$.
    
    \item The type $p|_A$ has a unique nonforking extension to $M$, which is $p$.

\end{enumerate}
\end{prop}

\begin{proof}
For points 2, 3, 4 and 6, see Proposition 2.1 in \cite{pillay-tanovic}. For point 1, see the proof of \cite[Proposition 2.1.i]{pillay-tanovic}. Let us now prove point 5. Let $B \subset M$ be a small set such that $p$ is $B$-invariant. Then, by point 3, $p$ is $A$-definable, so $B$-definable as well. By elimination if imaginaries, the type $p$ is thus $dcl(A)\cap dcl(B)$-definable. Let $C = dcl(A)\cap dcl(B)$.

\begin{claim} The type $p$ is generically stable over $C$. \end{claim} \begin{proof} By Remark \ref{rem_ops_alpha_denombrable}, let $\alpha$ be a countable infinite ordinal, let $(a_i)_{i < \alpha}$ be a Morley sequence of $p$ over $C$, made of elements of $M$. By contradiction, assume that the mean of the types of the $a_i$ over $M$ is not a complete type. Let $(a^{'}_i)_{i < \alpha}$ be a Morley sequence of $p$ over $A$, made of elements of $M$ as well. Then, the infinite tuples $(a_i)_{i < \alpha}$ and $(a^{'}_i)_{i < \alpha}$ have the same type over $C$. So, by strong homogeneity of $M$, there exists $\sigma \in Aut(M/C)$ such that $\sigma(a_i) = a^{'}_i$ for all $i < \alpha$. Since the mean of the types of the $a_i$ over $M$ is not a complete type, we deduce that the mean of the types of the $a^{'}_i$ over $M$ is not a complete type either, which contradicts the assumption of generic stability of $p$.\end{proof}

Then, unfolding the definition of generic stability, we deduce that $p$ is generically stable over $dcl(B) \supseteq C$, so over $B$ as well. \end{proof}

\begin{defi}\label{defi_gen_stable_abus}

\begin{enumerate}
\item Because of these properties, we may call a type $p \in S(A)$ generically stable if, for some, equivalently for all, sufficiently saturated model $M \supset A$, the type $p$ has a (necessarily unique) nonforking extension $q \in S(M)$ which is generically stable over $A$. 

\item If $B \supseteq A$, we may also say that a type $q \in S(B)$ is generically stable over $A$ if $q|_A$ is generically stable in the above sense, and $q$ does not fork over $A$.
\end{enumerate}
\end{defi}

The following Fact is a consequence of stationarity.

\begin{fact}\label{fact_gen_stable_abus}
If $q \in S(M)$ does not fork over $A$, where $M \supset A$ is a (not necessarily sufficiently saturated) model, then $q|_A$ is generically stable, in the above sense, if and only if $q$ is generically stable over $A$.
\end{fact}

\begin{rem}
The definition of generically stable types above is stronger than that of \cite{garcia-onshuus-usvyatsov}(Definition 1.8). More precisely, a type $p \in S(A)$ is generically stable, in the sense of \cite{garcia-onshuus-usvyatsov}, if and only if all its extensions to $acl(A)$ are generically stable, in the sense of Definition \ref{defi_gen_stable_abus}. This is why our definition implies stationarity, whereas that of \cite{garcia-onshuus-usvyatsov} does not. However, it is the only difference.
\end{rem}

%\begin{lemma}\label{passe-passe_base_canonique_cas_general}Let $C$ be a set, $\alpha$, $\alpha^{'}$, $\beta$ be tuples such that $tp(\beta / C)$ is generically stable in the sense of Definition \ref{defi_gen_stable_abus}, and $\alpha^{'} \equiv_C \alpha$. Assume that $\beta \downfree_C \alpha $ and $\beta \downfree_C \alpha^{'}$ . Then $\alpha^{'} \beta \equiv_C \alpha \beta $.\end{lemma}

%\begin{proof}The type $tp(\beta/C)$ is stationary, and its unique  nonforking extension is $C$-definable, thus $C$-invariant. The conclusion follows.\end{proof}

%\begin{prop}\label{prop_carac_gen_stable}Let $p \in S(M)$ be $A$-invariant, where $A \subset M$ is small. The following are equivalent:\begin{enumerate}    \item The type $p$ is generically stable over $A$.    \item For any Morley sequence $(a_i)_{i < \omega + \omega}$ of $p$ over $A$, the mean of the types $tp(a_i/M)$ is equal to $p$.    \item For any partitioned formula $\phi(x,y)$ over $A$, the partial type $\Pi_{\phi}(\overline{x}, \overline{y}):= p^{\otimes \omega}|_A(\overline{x}) \cup \lbrace \phi(x_i, y_j) \, | \, i <j < \omega \rbrace \cup \lbrace \neg\phi(x_i, y_j) \, | \, j \leq i < \omega \rbrace$ is inconsistent.  In other words, no Morley sequence of $p$ over $A$ witnesses the order property.\end{enumerate} \end{prop}

\end{subsection}

\begin{prop}\label{p|_B_est_gen_stable}

Let $p \in S(A)$ be a generically stable type. Let $B$ be a set of parameters containing $A$. Let $q \in S(B)$ be the unique nonforking extension of $p$. Then, $q$ is still generically stable and, for any $C \supseteq B$,  we have $q|_C = p|_C$.

\end{prop}

\begin{proof}Let $M \supset C$ be a sufficiently saturated model, and let $p^{'} \in S(M)$ be the nonforking extension of the type $p \in S(A)$. We know that $q$ does not fork over $A$. So, by Fact \ref{fact_forking}, $q$ has an extension $q' \in S(M)$ which does not fork over $A$. Then, $q^{'}|_A =q|_A = p$. So, by stationarity of $p$, we have $q^{'} = p^{'}$. So $q^{'}$ is generically stable over $A$, so a fortiori over $B$. Then, by Proposition \ref{proprietes_type_gen_stable}, $q^{'}|_B = q$ is stationary. Therefore, $q^{'}$ is indeed the nonforking extension of $q$, and $q^{'}$ is generically stable over $B$. Since we have also proved the equality   $q^{'} = p^{'}$, we are done.\end{proof}

%\begin{rem}\label{rem_q^n_est_definissable}
%If $q \in S(M)$ is generically stable over $A$, then, for all natural number $n$, the type $q^{\otimes n}$ is definable over $A$. In particular, the restriction $q^{\otimes n}|_A$ does not fork over $A$, i.e. is extensible.
%\end{rem}

\begin{prop}{Transitivity}\label{transitivite_gen_stable}

Let $p \in S(M)$ be a type generically stable over $A$. Let $B,C$ be sets of parameters such that $A \subseteq B \subseteq C \subset M$. Let $a \in M$ be a realization of $p|_A$, such that $a \downfree_A B$ and $a \downfree_B C$. Then $a \downfree_A C$. 
\end{prop}

\begin{proof}
Using stationarity and Fact \ref{p|_B_est_gen_stable} above, one can check that $tp(a/C) = p|_C$.
\end{proof}

\begin{prop}{Symmetry}\label{prop_symetrie} \emph{(\cite{garcia-onshuus-usvyatsov}, Theorem A.2 and Lemma A.5)}

Let $p \in S(A)$ be generically stable. Let $q \in S(A)$ be a type which does not fork over $A$. Let $a,b$ be such that $a \models p|_A$ and $b \models q$. Then $a \downfree_A b $ if and only if $ b \downfree_A a$. 

\end{prop}

The following lemma will be used repeatedly throughout the proof of the group configuration theorem.

\begin{lemma}\label{lemma_switch}{Swap}

Let $A$ be a set of parameters. Let $b,c,d$ be elements such that the types $tp(b/A), \,tp(c/A)$, and $tp(d / A)$ are generically stable.

\begin{enumerate}
    \item If $c \downfree_A d$ and $b \downfree_A cd$, then $bc \downfree_A d$.
    
    \item If $b \downfree_A c$ and $bc \downfree_A d$, then $b \downfree_A cd$.
\end{enumerate}

\end{lemma}

\begin{proof}

The first point is a consequence of Fact \ref{fact_forking} (3), and always holds. Let us prove the second point. We know that $bc \downfree_A d$, in particular $tp(bc / A)$ does not fork over $A$. Since $tp(d / A)$ is generically stable, we can apply symmetry, to deduce $d \downfree_A bc$. Also by symmetry, we have $c \downfree_A b$. So, by the first point applied to $(d, c, b)$, we have $dc \downfree_A b$. In particular, $dc \downfree_A A$. So, by symmetry again, we get $b \downfree_A cd$, as required.
%Let us prove the first point. We use symmetry, transitivity and stationarity: We know that $c \downfree_A d$ so, by symmetry for generically stable types, we have $d \downfree_A c$. Thus, by stationarity, $tp(d / Ac)$ is the unique nonforking extension of $tp(d/A)$, and is still generically stable. Besides, we have $b \downfree_A cd$, so $b \downfree_{Ac} d$. Since both $tp(b / Ac)$ and $tp(d / Ac)$ are generically stable, we may apply symmetry again, to get $d \downfree_{Ac} b$. Then, by transitivity, $d \downfree_A bc$. Now, by hypothesis and stationarity, $tp(bc / A)$ is the tensor product of $tp(b / A)$ and $tp(c / A)$, so is extensible. As $tp(d/A)$ is generically stable, we may apply symmetry one final time, to conclude that $bc \downfree_A d$, as desired.
\end{proof}

%\begin{lemma}\label{lemma_indep_descente_gen_stable}Let $C \subseteq D$. Let $a, b$ be such that $tp(a/C)$ is generically stable, in the sense of Definition \ref{defi_gen_stable_abus}, and $a \downfree_C D$. Assume that $a\downfree_D b$. Then $a\downfree_C b$.\end{lemma}

%\begin{proof}We have $a \downfree_D b$ and $a\downfree_C D$. Moreover, $tp(a/C)$ is generically stable. So, by transitivity, $a \downfree_C Db$, so $a\downfree_C b$, as desired.\end{proof}

\begin{prop}\label{prop_gen_stabilite_et_acl}

Let $p= tp(a/M)$  be a type generically stable over $A$, where $A \subseteq M$. 
\begin{enumerate}
    \item If $b \in dcl(Aa)$, then the type $tp(b / M)$ is generically stable over $A$.
    
    \item If $b \in acl(Aa)$, then the type $tp(b / M)$ is generically stable over $acl(A)$.
\end{enumerate}
\end{prop}

\begin{proof} Note that the hypothesis implies $a \downfree_A M$. Then, up to picking a bigger model $N \supseteq M$ such that $a \downfree_A N$, we may assume that $M$ is sufficiently saturated.
Let us prove the first point. Let $b \in dcl(Aa)$. Let $f$ an $A$-definable map  such that $f(a) =b$. We know that $tp(a/M)$ is definable over $A$. So $tp(b/M)$ is definable over $A$, so $A$-invariant. Let $q=tp(b/M)$. Thus, it remains to show that the property of the Morley sequences of Definition \ref{defi_gen_stable} holds. Let $\alpha$ be a countable ordinal.  Let $(a_i)_{i < \alpha}$ be a Morley sequence of $p$ over $A$, made of elements of $M$. 
Then, we will show that $(f(a_i))_i$ is a Morley sequence of $q$ over $A$. 

For all $i < \alpha$, the element $a_i$ realizes the type $p|_{A \cup(a_j)_{j < i}}$. In other words, $a_i \equiv_{A \cup(a_j)_{j < i}} a$. Since $f$ is $A$-definable, we deduce that $f(a_i) \equiv_{A \cup(a_j)_{j < i}} f(a)$. A fortiori, as $f(a)=b$, we have $f(a_i) \equiv_{A \cup(f(a_j))_{j < i}} b$. So $f(a_i)$ realizes the type $q|_{A \cup(f(a_j))_{j < i}}$. So $(f(a_i))_i$ is a Morley sequence of $q$ over $A$, made of elements of $M$. %Up to applying an automorphism of $M$ fixing $A$ pointwise, we may assume that $f(a_i)=b_i$ for all $i$.  
Then, the property of the mean of the types of the $f(a_i)$ is a consequence of that for the types of the $a_i$. So we have found a Morley sequence made of elements of $M$ with this property; by homogeneity it holds for all Morley sequences of $q$.

Let us then prove the second point. Let $b \in acl(Aa)$. By Lemma \ref{definissabilite_et_acl}, we know that $q=tp(b/M)$ is definable over $acl(A)$, so is a fortiori $acl(A)$-invariant.

Let $\phi(y,x,m)$ be a formula with parameters in $M$ such that $\phi(y,a,m)$ isolates the type of $b$ over $Ma$. Let $M_0 \prec M$ be a small model containing $Am$. We will show that $q$ is generically stable over $M_0$. Let $r = tp(ab / M_0)$. Then, by construction, we have $p(x) \cup r(x,y) \models q(y)$. So, by Theorem 3.5 (3) in \cite{mennui}, the type $q$ is generically stable over $M_0$. Since it is $acl(A)$-invariant, we conclude by point 5 of Proposition \ref{proprietes_type_gen_stable} that it is generically stable over $acl(A)$, as required.%
% Let $(c_i)_i$ be a Morley sequence of $q$ over $acl(A) m$. By contradiction, assume that there exists a formula $\psi(y, d)$ with parameters in $acl(A) m d$, where $d \in M$, such that the set of indices $i$ such that $\models \psi(c_i, d)$ is infinite and co-infinite. We write $c_i = a_i b_i$. Then $(a_i)_i$ is a Morley sequence of $p$ over $acl(A) m$. A fortiori, it is a Morley sequence over $A$. We also know that, for all $i$, we have $M \models \phi(c_i, a_i, m)$. Thus, there exist infinitely many $i$ such that $M \models \exists y \, \phi(y, a_i, m) \wedge \psi(y, d)$. Similarly, there exist infinitely many $i$ such that $M \models \exists y \, \phi(y, a_i, m) \wedge \neg \psi(y, d)$. But $(a_i)_i$ is a Morley sequence of $p$, which is generically stable. So $p(x) \models \exists y \, \phi(y, x, m) \wedge \psi(y, d)$ and $p(x) \models \exists y \, \phi(y, x, m) \wedge \neg \psi(y, d)$. So $\models \exists y \, \exists z \, \phi(y, a, m)\wedge \phi(z,a,m) \wedge \psi(y, d) \wedge \neg \psi(z,d)$. However, the formula $\phi(y,a,m)$ isolates the type of $c$ over $Ma$. This is a contradiction. So $q=tp(c/M)$ is generically stable over $acl(A) m$. Now, recall that $q$ is definable, so invariant, over $acl(A)$. So $q$ is generically stable over $acl(A)$, as required. 
\end{proof}

\begin{subsection}{Strong germs}

In this subsection, we state useful results on germs of definable maps at generically stable types.

\begin{fact}\label{fact_definissabilite_des_germs}
Let $p$ be an $A$-definable type. Let $X$ be an $A$-definable set, and $(f_b)_{b \in X}$ be an $A$-definable family of definable maps, such that $f_b$ is defined on $p|_{Ab}$, for all $b \in X$. Then, the equivalence relation on $X$ defined by $b_1 \sim b_2$ if and only if $p|_{A b_1 b_2} \models f_{b_1}(x)=f_{b_2}(x) $ is $A$-definable, since the type $p|_{A b_1 b_2}$ is definable by the good defining scheme of $p$. 
\end{fact}

\begin{defi}
In the above context, if $b$ is an element of $X$, we shall let $[f_{b}]_p$, or $[f_{b}]$ if the context is clear, denote the code of the class of the element $b$ for the equivalence relation $\sim$ defined above. We call this code the \emph{germ} of the function $f_b$ at the type $p$.
\end{defi}

In general, the germ of a definable map at a given definable type encodes less information than the code of said definable map. In some sense, it only captures the ``local'' (for the Stone topology) behavior of the map.

\begin{notation}
If a type $p$ is definable, we let $Cb(p)$ denote its canonical basis, i.e. the definable closure of the codes of the formulas in its defining scheme. Similarly, if $tp(a/B)$ is definable, we write $Cb(a/B)$ for the canonical basis of $tp(a/B)$.
\end{notation}

%\begin{lemma}\label{critère_germs_equal}Let $p \in S(A)$ be a definable type, $X$ a definable set, $f_y$ a definable map (with a distinguished parameter $y$) such that, for $y$ in $X$, the domain of $f_y$ contains the realizations of $p|_{Ay}$. Let $y_1, y_2$ be elements of $X$. The following are equivalent:\begin{enumerate}    \item There exists $a$ realizing $p|_{A y_1 y_2}$ such that  $f_{y_1}(a)=f_{y_2}(a)$.\item For all $a$ realizing $p|_{A y_1 y_2}$, we have $f_{y_1}(a)=f_{y_2}(a)$.\end{enumerate}\end{lemma}

\begin{defi}\label{defi_f_defined_at_p}

Let $p \in S(A)$ be a definable type, and $f$ be a definable map, possibly using parameters outside of $A$. We say that $f$ is defined at $p$, or well-defined at $p$, if, for some $B \supseteq A$ such that $f$ is $B$-definable, the function $f$ is defined at $p|_{B}$.

\end{defi}

\begin{prop}\label{prop_description_base_canonique_af_c(a)}
Let $p \in S(B)$ be an $A$-definable type, where $A \subseteq B = acl(B)$. Let $a$ be a realization of $p$, let $c \in B$ and $f_c$ be an $Ac$-definable map such that $f_c$ is defined at $p$. Then, the canonical basis of $tp(a f_c(a) / B)$ is interdefinable over $A$ with the definable closure of the set $ [f_c] Cb(a/B)$. In other words, we have the following equality: $dcl(A Cb(af_c(a) / B)) = dcl(A [f_c] Cb(a/B) )$.
\end{prop}

\begin{proof} To simplify notations, let $C = Cb(a f_c(a) / B)$. Also, since we are dealing with definable types, up to working with the unique definable extensions, we may assume that $B$ is equal to a sufficiently saturated model $M$. Let us first prove the following: $C \subseteq dcl(A [f_c]Cb(a/M))$. Let $b=f_c(a)$.
It suffices to show that $tp(ab / M)$ is definable over $A [f_c]$. By Lemma \ref{definissabilite_et_acl}, we know that $q=tp(ab / M)$ is definable over $Cb(a/M) c$. Let us show that $q$ is invariant over $[f_c]Cb(a/M) $, which will be enough to conclude. Let $\sigma \in Aut(M /  [f_c]Cb(a/M))$. Let us show that $\sigma(q)=q$. By hypothesis on $\sigma$, we then have $[f_{\sigma(c)}]=\sigma([f_{c}])=[f_c]$. Thus, $p(x) \models f_c(x)=f_{\sigma(c)}(x)$. 

Besides, we have $p(x) \cup \lbrace f_c(x) = y \rbrace \models q(x,y)$. Therefore $\sigma(p)(x) \cup \lbrace f_{\sigma(c)}(x) = y \rbrace \models \sigma(q)(x,y)$. But $p$ is $A$-invariant, and we proved that $p(x) \models f_c(x)=f_{\sigma(c)}(x)$. So $p(x) \cup \lbrace f_{c}(x) = y \rbrace \models \sigma(q)(x,y)$. Since we also know that $p(x) \cup \lbrace f_{c}(x) = y \rbrace \models q(x,y)$, we deduce that $\sigma(q)(x,y)=q(x,y)$. Thus, the type $tp(ab / M)$ is definable over $Cb(a/M) [f_c]$, as desired.

Let us now prove the converse inclusion. It is clear that $Cb(a / M) \subseteq C= Cb(a f_c(a) / M)$. It remains to show that $[f_c] \in dcl(AC)$.  By definition of $C$, the type $q$ defined by $$q(x,y):=tp(a f_c(a) / M)(x,y),$$ is $C$-definable. We use compactness. Let $d$ be such that $d \equiv_{AC} c$. It suffices to show that $[f_d]=[f_c]$. We know, by choice of the type $q$, that $q(x,y) \models f_c(x)=y$. In other words, $\models d_q xy \, (f_c(x) = y)$. Note that the formula $\phi(z) = d_q xy \, (f_z(x) = y)$ is over $AC$, since $q$ is $C$-definable and $f$ is $A$-definable. By choice of $d$, we have $d \equiv_{AC} c$. Thus $\models d_q xy \, (f_d(x) = y)$. So, if $q^{'}$ is the $C$-definable extension of $q$ to $Md$, we have $$q^{'}(x,y) \models f_d(x) = y \, \wedge \,f_c(x) = y.$$ Finally, $q^{'}(x,y) \models f_d(x) = f_c(x)$, so $p(x)|_{Md} \models f_d(x) = f_c(x)$,  i.e. $[f_c] = [f_d]$, as desired.
\end{proof}

\begin{prop}\label{prop_b_in_acl(Aa)_cas_gen_stable}\emph{(\cite{adler-casanovas-pillay}, Lemma 2.1)}
Let $q(x,y) \in S(M)$ be a type generically stable over $A\subseteq M$. Let $a,b$ be a realization of $q(x,y)$.

\begin{enumerate}
    \item If $b \in acl(Ma)$, then $b \in acl(Aa)$.
    \item If $b \in dcl(Ma)$, then $b \in dcl(Aa)$.
\end{enumerate}
\end{prop}

\begin{coro}{Strong germs}\label{germs_forts_cas_gen_stable} \emph{(\cite{adler-casanovas-pillay}, Theorem 2.2)}

Let $p \in S(M)$ be a type generically stable over a set $A\subseteq M$, $f_c$ an $A c$-definable map  (with a distinguished parameter $c \in M$) such that $f_c$ is defined at $p$. Then, for $a \models p|_{Ac}$, we have $f_c(a) \in dcl(A a [f_c])$. In fact, there exists an $A [f_c]$-definable map $F_{[f_c]}$ such that $p(x)|_{Ac} \models f_c(x)=F_{[f_c]}(x)$.

\end{coro}

\begin{rem}
\begin{enumerate}
\item In fact, using Proposition \ref{prop_description_base_canonique_af_c(a)} above, one can show that Corollary  \ref{germs_forts_cas_gen_stable} is essentially the same as Proposition \ref{prop_b_in_acl(Aa)_cas_gen_stable}. Thus, both statements could be seen as ``the strong germs property''.

\item The strong germs property, will be crucial in the proof of Theorem \ref{theo_config_groupe_faible}. In fact, most of Section 3 will be devoted to the study of the action of germs of definable maps on certain generically stable types. In some sense, considering germs of definable maps enables us to build a type-definable group, instead of an ind-type-definable one. However, defining the group operation, and the action on the space, relies heavily on the strong germs property.

 \end{enumerate}
\end{rem}

%\begin{proof}By hypothesis, $tp(a/M)$ is generically stable over $A$. Let $b= f_c(a)$. We know, by Proposition \ref{prop_gen_stabilite_et_acl}, that $tp(ab /M)$ is generically stable over $A c$. Besides, by the proposition above, we know that $tp(ab / M)$ is definable over $A [f_c]$. Then, $tp(ab / M)$ is generically stable over $A [f_c]$. By Proposition \ref{prop_b_in_acl(Aa)_cas_gen_stable}, this implies that $b \in dcl(A a [f_c])$, as desired. \end{proof}

%
%
%\begin{lemma}\label{lemme_germ_definable_sur_suite_Morley_cas_gen_stable}
%
%
%Let $p \in S(A)$ be a generically stable type, $h_b$ an $Ab$-definable map (with a distinguished parameter $b$) such that $p(x) \models x \in dom(h_b)$. Let $(e_i)_{i \in I}$ be a Morley sequence of $p$ over $A$, where $I$ is infinite. Then, the germ of $h_b$ at $p$ is definable over $(e_i, h_b(e_i))_{i \in I} \cup A$.
%\end{lemma}
%
%
%
%\begin{proof}
%By compactness, it suffices to prove that, if $b^{'}\equiv_F b$, then $[h_{b^{'}}]=[h_b]$, where $F = (e_i, h_b(e_i))_{i \in I} \cup A$. By generic stability of $p$, since we know that $\models h_b(e_i)=h_{b^{'}}(e_i)$ for all $i \in I$, and that $I$ is infinite, we deduce that $p|_{Abb^{'}}(x) \models h_b(x) = h_{b^{'}}(x)$.
%\end{proof}

\end{subsection}

\begin{subsection}{Commutativity}

\begin{fact}{Commutativity} \label{fact_gen_stable_commute_with_lui_meme}\emph{(see \cite[Remark 5.18]{Conant2021KeislerMI})}

Let $p, q \in S(M)$ be $A$-invariant types, where $A$ is a small set contained in $M$. Assume that $p$ is generically stable over $A$. Then, $p(x) \otimes q(y) = q(y) \otimes p(x)$, this equality being  between $A$-invariant types.

\end{fact}

\begin{defi}\label{defi_image_type}
If $p = tp(a / A)$ is a complete type, and $h$ is an $A$-definable map defined at $p$, we let $h_{*} p$, or $h(p)$, denote the type $h_{*} p = tp(h(a) / A)$. It is called the image of $p$ under $h$. Note that this does not depend on the choice of the realization $a$. 
\end{defi}

\begin{rem}\label{rem_h(p)}

In the definition above, if $p$ is $A$-invariant (resp. $A$-definable, resp. generically stable over $A$), then so is $h_*p$, and we have $h_*(p|_B)=(h_*p)|_B$ for all $B \supseteq A$.

\end{rem}

\begin{defi}\label{defi_acting_generically}
Let $p$ be a definable type, and $f$ a definable family of definable maps. We say that an element $a$ acts generically on $p$ via $f$, if  the definable map $f_a$ is well-defined at $p$, in the sense of Definition \ref{defi_f_defined_at_p}. If the definable family of definable maps $f$ is implicit, we just say that $a$ acts generically on $p$.

We say that a type-definable set $X$ acts generically on $p$ if, for some implicitly given $f$, all elements $a \in X$ act generically on $p$ via $f$.
\end{defi}

%\begin{lemma}\label{relevement_generique_cas_gen_stable}Let $B \subseteq C$, let $q \in S(C)$ and $p=q|_{B}$. Let $h$ be a $B$-definable map  such that, for all $a \models p$, the map $h$ is defined at $a$. Let $k \models h_{*} q$. Then, there exists a realization $\alpha$ of $q$ such that $h(\alpha)=k$.\end{lemma}

%\begin{proof}Let $a^{'}$ realize $q$. Then $h(a^{'})$ realizes $h_{*} q$, i.e. $h(a^{'}) \equiv_C h(a)$. So, there exists $\alpha$ such that $a^{'} h(a^{'}) \equiv_C \alpha h(a)$. Then, $\alpha$ realizes $q$, and $h(\alpha)=h(a)$. \end{proof}

\begin{defi}
Let $p(x) \in S(M)$ be an $A$-invariant type, where $A \subset M$ is a small set, and let $\mathcal{F}$ be a set of invariant types. We say that $p$  \emph{commutes with $\mathcal{F}$} if, for all invariant types  $q(y)$ in $\mathcal{F}$, we have $p(x) \otimes q(y) = q(y) \otimes p(x)$, this being an equality of invariant types. See Remark \ref{rem_invariant_types_as_processes} for an explanation of this idea.

\end{defi}

\begin{coro}\label{coro_limage_dun_type_symetrique_est_symetrique}

Let $p$ be a $B$-invariant type, and $h$ a $B$-definable map  such that $h$ is defined at $p$, and $\mathcal{F}$ a family of  invariant types. If $p$ commutes with $\mathcal{F}$, then $h_{*} p$  commutes with $\mathcal{F}$.

\end{coro}

\begin{proof}

Let $q(y)$ be an element of $\mathcal{F}$. By hypothesis on $\mathcal{F}$, there exists a small set $C \supseteq B$  such that $q(y)$ is $C$-invariant. Let us show that $h_{*} p(x)\otimes q(y) =  q(y)\otimes h_{*} p(x)$. Let $D \supseteq C$, and let $(k,b)$ realize $h_{*} p(x)\otimes q(y) |_{D}$. 

Then, $k$ realizes $h_{*} p|_{Db}$. So, by Remark \ref{rem_h(p)} applied to $Db$, there exists $a$ realizing $p|_{Db}$ such that $h(a)=k$. So $(a,b)$ realizes $p \otimes q|{D}$. Since $p$ commutes with $\mathcal{F}$, $(b,a)$ realizes $q\otimes p|_{D}$, so $b \models q|_{Da}$, a fortiori $b \models q|_{Dh(a)}=q|_{Dk}$. Therefore, $(b,k)$ realizes $q \otimes h_{*} p |_{D}$. We have proved that $h_{*} p$ commutes with $\mathcal{F}$.
\end{proof}

\begin{rem}
These notions give us some form of symmetry for tensor products of generically stable types; see Lemma \ref{composee_germs_indep_cas_gen_stable} for an example of how this commutativity is used. However, we do not know if the class of generically stable types is closed under tensor products, outside well-behaved theories, e.g. NIP.
\end{rem}

\end{subsection}

\end{section}

\begin{section}{The group configuration theorem}

\begin{subsection}{Genericity and group configurations}

Here, we define a notion of genericity for definable types concentrating on a type-definable group $G$, or $G$-space $X$. We then define group configurations, and explain how to build such using generic types. Few of the results are new, except maybe Propositions \ref{prop_gen_stable_generics_in_a_space} and \ref{prop_how_to_build_quadrangles} in the case of $G$-spaces, which is well-known for stable theories. For more results, and a more general framework allowing definable \emph{partial} types to be generic, see Section 3 in \cite{HruRK-MetaGp}.

\begin{defi}
\begin{enumerate}
\item A type-definable group $\Gamma$ is given by a type-definable set, along with a relatively definable map $m: \Gamma \times \Gamma \rightarrow \Gamma$ which defines a group operation.

\item Let $G$ be a type-definable group, and $X$ a type-definable space on which $G$ acts definably. Assume everything is defined over some set $A$. Let $\cdot$ denote both the group operation of $G$, and the action of $G$ on $X$.%, and $p \in S(B)$ be a definable type such that $p(x) \models X \in x$, with its defining scheme $d_p$. For all partitioned formulas without parameters $\phi(x,y)$, we define $Stab_{\phi}(p)$ via the following formula,  in the variable $g$ : $\forall y \,  d_p x \, [\phi(g\cdot x , y) \leftrightarrow \phi(x, y)]$.

%We then define the (left) stabilizer of $p$ in $G$, denoted by  $Stab(p)$, as the intersection of the $Stab_{\phi}(p)$. 
Let $B \supseteq A$, and $p, q \in S(B)$ be definable types concentrating on $X$. We define $Stab_{\phi}(p,q)$ by the formula $\forall y \,  [d_q x \, \phi(g\cdot x , y) \leftrightarrow d_p x \, \phi( x, y)]$. We then define  $Stab(p,q)$ as the intersection of all the $Stab_{\phi}(p,q)$ with $G$. If $p=q$, we write $Stab(p)$ instead of $Stab(p,p)$. In the case where $p,q \in G$, we also define the right stabilizer $Stab^r(p,q)$ by considering the right action by translations, and similarly for $Stab^r(p)$.

\end{enumerate}
 %More generally, we can define $Stab(p), Stab^r(p), Stab(p,q), Stab^r(p,q)$ for definable types $p, q \in X$, if $X$ is a space on which $G$ acts definably. 
\end{defi}

%\begin{fact}\label{fact_stab(p)_est_ss-grpe}If $p \in S(B)$ is a definable type such that $p(x) \models $ ``$x \in X$'', where $B$ contains the parameters defining $G$ and $X$, then the type-definable set $Stab_{}(p)\cap G $ is a subgroup of $G$.\end{fact}

\begin{rem}

Let $p,q \in S(B)$ be as in the definition above. Let $M$ be a model over which everything is defined. Then $Stab(p,q)(M)$ is precisely the set of elements $g \in G(M)$ such that $g\cdot p|_M = q|_M$. Also, $Stab(p)$ is a type-definable subgroup of $G$.

\end{rem}

%In this context, we let $Stab(p)$ or $Stab_G(p)$ denote the type-definable subgroup $Stab_{}(p)\cap G $. Similarly, we write $Stab(p,q)$ instead of $Stab(p,q) \cap G$.

\begin{defi}
Let $G$ be a type-definable group acting definably on a type-definable space $X$. Let $M$ be a  sufficiently saturated model over which everything is defined.

\begin{enumerate}
    \item Let $H \leq G$ be a type-definable subgroup.  We say that $H$ is of bounded index in $G$ if the cardinality of $G/H$ is bounded, i.e., does not grow beyond a fixed cardinal, regardless of the size of the model.
    \item Let $p \in S(M)$ be a definable type. We say that $p$ is a definable generic of the $G$-space $X$ if $p(x) \models $ \textquotedblleft $x \in X$\textquotedblright, and $Stab(p)$ is of bounded index in $G$. Letting $G$ act definably and regularly on itself by left translations, we can also speak of definable generic types in $G$.
    \item We say that the space $X$ is connected if it has a definable generic type over $M$ whose stabilizer is $G$ itself. It is generically stable if it has a generically stable generic. The group $G$ is connected (resp. generically stable) if it is connected (resp. generically stable) for the left regular action by translations.

\end{enumerate}

\end{defi}

\begin{rem}
Other, weaker notions of genericity have been developed. For instance, there is a notion of f-genericity, which relies on forking rather than definable types. See \cite[Definition 3.3]{Stabilizers-NTP2}. However, in this paper, we will only be interested in definable generics. Thus, we shall call them ``generics''.
\end{rem}

\begin{lemma}\label{lemma_pb_is_generic}
Let $G$ be a type-definable group. Let $X, Y$ be type-definable $G$-spaces, and $f : X \rightarrow Y$ be a definable $G$-equivariant map. Let $p \in S_X(M)$ be a definable type, where $M$ is a model containing all the parameters involved. Then, $Stab(f_*(p)) \geq Stab(p)$. In particular, if $p$ is generic in $X$, then $f_*(p)$ is generic in $Y$.

\end{lemma}

\begin{proof}
We may assume that $M$ is sufficiently saturated. Let $c \in Stab(p)(M)$. Then, we compute $c \cdot f_*(p) = f_*(c \cdot p)= f_*(p)$, so $c \in Stab(f_*(p))$, as required.
\end{proof}

%\begin{rem}Il existe une autre notion of genericity, defined à partir of recouvrements finis by des translatés d'sets definable. Dans the cas of a stable theory, ces deux notions coïncident. Voir l'annexe for plus of détails. \end{rem}

\begin{prop}\label{prop_criterion_connectedness}
Let $G$ be a type-definable group with a definable generic type. The following are equivalent:

\begin{enumerate}
\item The group $G$ is connected.

\item The group $G$ has no type-definable proper subgroup of bounded index.
\end{enumerate}

It these hold, then, for any definable generic type $p$, we have $Stab(p) = G$.
\end{prop}

\begin{proof}
The implication $2. \implies 1.$ is straightforward: by definition, the stabilizer of any generic type is of bounded index. Let us prove $1. \implies 2.$ Let $p$ be a generic type for $G$, whose stabilizer is $G$ itself. Let $H \leq G$ be a type-definable subgroup of bounded index. Then, if $M$ is a sufficiently saturated model containg all the parameters involved, it represents every coset of $H$. Then, $p|_M$ concentrates on a coset of $H$. So, the stabilizer $Stab(p)$ is contained in a conjugate of $H$. As $Stab(p) = G$, we deduce that $H=G$, as desired.
\end{proof}

\begin{lemma}\label{lemma_gen_stable_generics_in_a_group} \emph{(\cite{pillay-tanovic} Lemma 2.1)}

Let $G$ be a type-definable group, defined over a set of parameters $A$.
Let $p \in S(A)$ be a generically stable generic type for $G$, such that $Stab(p)=G$.  \begin{enumerate}
    \item Let $B \supseteq A$ and $a$ realizing $p|_{B}$. Then, the element $a^{-1}$ realizes $p|_{B}$. In other words, we have $p^{-1} = p$.
    \item Let $g \in G$ and $a$ realizing $p|_{Ag}$. Then, the element $a\cdot g$ realizes $p|_{Ag}$. In other words, the right stabilizer of $p$ is also equal to the whole group $G$.
    \item The type $p$ is the unique generic type of the group $G$.
    
    \item Any element of $G$ is the product of two realizations of $p$.
\end{enumerate}
\end{lemma}

\begin{proof}
Let $\alpha, \beta$ realize $(p\otimes p)|_B$. Then, since $\beta^{-1} \in Stab(p) = G$ and $\alpha \models p|_{B\beta}$, we know that $\beta^{-1} \cdot \alpha$ realizes $p|_B$. Then, by Fact \ref{fact_gen_stable_commute_with_lui_meme}, we know that $\alpha \beta \equiv_B \beta \alpha$. So $\alpha^{-1} \cdot \beta$ realizes $p|_B$. Then, since $p|_B$ is a complete type, we have shown that, for all elements $c$ realizing $p|_B$, the element $c^{-1}$ realizes $p|_B$, as desired.

Let us then prove  the second point. If $g$ is in $ G$ and $a$ realizes $p|_{Ag}$, then, by the first point, we know that $a^{-1}$ realizes $p|_{Ag}$. So, by hypothesis on the stabilizer of $p$, $g^{-1} \cdot a^{-1}$ realizes $p|_{Ag}$. Then, again by the first point, the element $a \cdot g = (g^{-1} \cdot a^{-1})^{-1}$ realizes  $p|_{Ag}$, as stated.

Let us now prove the third point. Let $q \in S(B)$ be a definable generic type for $G$. Let $M$ be a sufficiently saturated model containing $A$ and $B$. Then, by Proposition \ref{prop_criterion_connectedness}, we have $G = Stab(q)$. Now, let $(a,b)$ realize $(p \otimes q) |_M$. Then, by the second point, we know that $a \cdot b $ realizes $p|_{Mb}$. On the other hand, since $Stab(q)=G$, we know that $a\cdot b$ realizes $q|_{Ma}$. So $q|_M=p|_M$. So $p$ and $q$ are equal as definable types.

Finally, let us prove the fourth point. Let $g \in G$, and $a \models p|_{Ag}$. Then, we have $g= (g \cdot a) \cdot a^{-1}$,where $g \cdot a$ realizes $p$ because $Stab(p) = G$, and $a^{-1}$ realizes $p$ by the first point.
\end{proof}

\begin{coro}\label{coro_G0_exists}\emph{(\cite{HruRK-MetaGp}, Lemma 3.9)}
Let $G$ be a generically stable type-definable group. Then, the connected component $G^{00}$ of $G$, i.e. the smallest type-definable subgroup of bounded index, exists. The group $G^{00}$ is the stabilizer of any generic type of $G$, and has a unique generic type.

\end{coro}

\begin{proof}
Let us first prove the existence of the connected component $G^{00}$. Let $p \in S(M)$ be a generically stable generic type for $G$, where $M$ is sufficiently saturated. Then, as $H=Stab(p)$ is of bounded index, every left coset and every right coset of $Stab(p)$ is represented in $M$. So $p$ concentrates on some right coset of $H$. Now, let $q$ be the image of $p \otimes p$ under the map $(a, b) \mapsto a \cdot b^{-1}$. Then, $q$ concentrates on $H$. 

\begin{claim}
 We have $Stab(q) \geq Stab(p) = H$. 

\end{claim}
\begin{proof}
This is a consequence of Lemma \ref{lemma_pb_is_generic}, applied to $f : g \mapsto g \cdot b^{-1}$.
\end{proof}

Thus, $q$ concentrates on $H_1:= Stab(q)$. So, by Proposition \ref{prop_criterion_connectedness}, the type-definable group $H_1$ has no proper type-definable subgroup of bounded index. Moreover, as $H \leq H_1 \leq G$, we know that $H_1$ is of bounded index in $G$. Hence, $H_1$ is indeed the smallest type-definable subgroup of bounded index of $G$. This implies that $H_1$ is normal, and even invariant under definable automorphisms, in $G$. Also, by Lemma \ref{lemma_gen_stable_generics_in_a_group}, the unique generic type of $H_1$ is the type $q$.

Now, let $p_1$ be some generic of $G$. By definition, $Stab(p_1)$ is of bounded index, so it contains $H_1$. On the other hand, the complete type $p_1$ concentrates on some coset of $H_1$, so $Stab(p_1)$ is contained in some conjugate of $H_1$. Since the latter is normal, we have in fact $Stab(p_1) = H_1$, as desired.
\end{proof}

\begin{prop}\label{prop_gen_stable_generics_in_a_space}

Let $G$ be a generically stable type-definable group acting definably and transitively on a type-definable space $X$. 

\begin{enumerate}
\item If $G$ is connected, then $X$ has a unique generic type, whose stabilizer is $G$. 

\item In general, the space $X$ has generically stable generics, and they are left translates of each other.
\end{enumerate}

\end{prop}

\begin{proof}
Let us show that the second point follows from the first. We know that the connected component $G^{00}$ of $G$ exists: it is the stabilizer of any generic of $G$. Then, we consider the action of $G^{00}$ on $X$. By the first point, each $G^{00}$-orbit contains a unique generic type, whose stabilizer is $G^{00}$. Note that, since $G^{00}$ is normal and of bounded index in $G$ and the action is transitive, there are only boundedly many $G^{00}$-orbits. Now, let $M$ be a sufficiently saturated model over which everything is defined. So $M$ contains a point in each $G^{00}$-orbit. Let $q_1, q_2 \in S(M)$ be two generic types of $X$. Let $x_1, x_2 \in X(M)$ be in the $G^{00}$-orbits of (the realizations of) $q_1$ and $q_2$ respectively, and let $g \in G(M)$ be such that $g \cdot x_1 = x_2$. We shall prove that $g$ sends the type $q_1$ to $q_2$. Since $G^{00}$ is normal, we have $g \cdot G^{00} = G^{00} \cdot g$, so $g \cdot q_1$ concentrates on the same $G^{00}$-orbit as $q_2$, and is still generic. So, by the first point, we have $g \cdot q_1 = q_2$, as desired.

Let us now prove the first point. Assume that $G$ is connected. By Lemma \ref{lemma_gen_stable_generics_in_a_group}, let $p$ be the unique generic type of $G$. Let $M$ be a big enough model containing all the parameters involved, and let $x_0 \in X(M)$. Let $f: G \rightarrow X$ be the definable map $g \mapsto g \cdot x_0$. By transitivity of the action, it is onto. Let $q$ be the definable type $f_*(p)$. It is easy to check that $Stab(q)(M) = G(M)$, because $Stab(p) = G$. So, since $M$ is sufficiently saturated, we have $Stab(q) = G$, so $q$ is generic in $X$. 

Now, let $q_1$ be another generic type in $X$. Without loss of generality, we may assume that $q_1$ is $M$-definable. We want to show that $q_1 = q$. Let $x_1 \models q_1|_M$, and $g_1 \in G$ such that $f(g_1) = x_1$. Let  $g \models p|_{Mg_1}$. Then, by Lemma \ref{lemma_gen_stable_generics_in_a_group}, we have $g \cdot g_1 \models p|_{M g_1}$, so $f(g \cdot g_1) \models f_*(p)|_{M g_1} = q|_{M g_1}$. In particular, we have $f(g \cdot g_1) = g \cdot x_1 \models q|_M$. On the other hand, since $q_1$ is generic in $X$, and $G$ is connected, we have $Stab(q_1) = G$. Moreover, by Fact \ref{fact_gen_stable_commute_with_lui_meme}, we have $(x_1, g) \models (q_1 \otimes p)|_M$. So $g \cdot x_1 \models q_1|_{M g}$. Therefore, $g\cdot x_1$ realizes both $q|_M$ and $q_1|_M$. So $q=q_1$, as desired.
\end{proof}

\begin{defi}\label{defi_configuration_groupe}
Let $A$ be a set of parameters.   A \emph{regular group configuration} over $A$ is a tuple $(a_1,a_2, a_3, b_1, b_2, b_3)$ of elements satisfying the following properties:

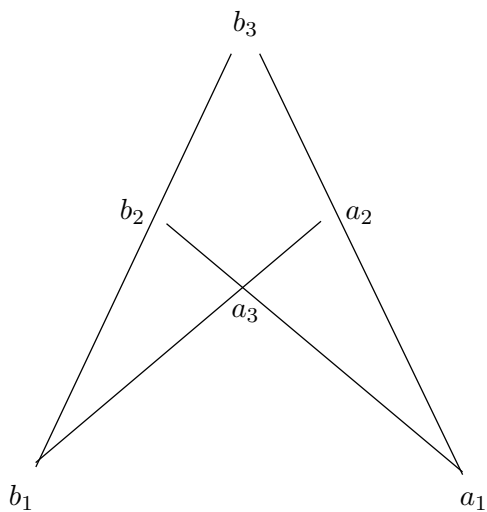
\begin{center}
    \begin{tikzcd}

&& b_3 \arrow[llddddd, dash, shorten= 1 mm] \arrow[rrddddd, dash, shorten= 1 mm]

\\
 
\\
 & b_2 & & a_2
\\
&& a_3

\\

\\

b_1 \arrow[rrruuu, dash, shift left = 1.4 ex]&&&& a_1\arrow[llluuu, dash, shift right = 1.0 ex]

    \end{tikzcd}
\end{center}

\begin{enumerate}    
    \item If $\alpha,\beta, \gamma$ are three non-colinear points in the diagram above, then the triplet $(\alpha,\beta, \gamma)$ is an independent family over $A$.
    \item If $\alpha,\beta, \gamma$  are three colinear points in the diagram above, then $\alpha \in \mathrm{acl}(A\beta \gamma)$.

\end{enumerate}

    A \emph{definable group configuration} over $A$ is a tuple  of elements $(a_1,a_2, a_3, b_1, b_2, b_3)$ satisfying the following properties:\begin{enumerate}    \item The type $tp(a_1 a_2 a_3 b_1 b_2 b_3 / A)$ is definable.    \item If $\alpha,\beta, \gamma$ are three non-colinear points in the diagram above, then the types $tp(\alpha \beta/ A\gamma)$ and $tp(\alpha/ A\beta \gamma)$ are $A$-definable. %the triplet $(x,y,z)$ is an independent family over $A$, and  %realizes the tensor product $tp(x/A) \otimes tp(y/A) \otimes tp(z/A)$.    
    \item The equalities $acl(A b_1 b_2) = acl(A b_1 b_3) = acl(A b_2 b_3)$ hold.   \item For all natural numbers $i,j,k$ such that $\lbrace i,j,k \rbrace = \lbrace 1,2,3 \rbrace$, the elements $a_j$ and $a_k$ are interalgebraic over $A b_i$, and the element $b_i$ is interalgebraic over $A$ with the canonical basis $Cb(a_j a_k / acl(A b_i))$.\end{enumerate}

A \emph{generically stable (resp. generically stable regular) group configuration} over $A$ is a definable (resp. regular) group configuration $(a_1,a_2, a_3, b_1, b_2, b_3)$ over $A$, such that the type $tp(a_1 a_2 a_3 b_1 b_2 b_3 / A)$ is generically stable, in the sense of Definition \ref{defi_gen_stable_abus}.

%We also use the name ``group configuration'', instead of ``group configuration''.

\end{defi}
 We might call ``quadrangle''  a $6$-tuple of elements which has not been proven to be a definable or regular group configuration (yet).

\begin{rem}\label{rem_defi_group_config}
\begin{enumerate}
\item Recall that, by Definition \ref{defi_definable_types}(5), a type $tp(a/ A b)$ is $A$-definable if and only if it admits an $A$-definable extension to a model. This implies that $a \downfree_A b$. In particular, this implies that, in definable group configurations, non-colinear triples are independent families.
\item In the case of a generically stable $6$-tuple, independence over $M$ of the non-colinear triplets can be checked more easily, using Lemma \ref{lemma_switch}. For instance, the set $\lbrace a_1, a_2, a_3 \rbrace$ being independent over $M$ is equivalent to having $a_1 \downfree_M a_2$ and $a_3 \downfree_M a_1 a_2$.  
\end{enumerate}

\end{rem}

\begin{prop}\label{prop_how_to_build_quadrangles}
Let $G$ be a  type-definable  connected group acting definably on a type-definable connected space $X$. Assume that the action is free (resp. faithful) and transitive. Let $p,q$ be the generics of $G$ and $X$ respectively. Assume that $p$ and $q$ are generically stable. Let $(g_1, g_2, x)$ be a triplet realizing $p^{\otimes 2}\otimes q|_M$, where $M$ is a sufficiently saturated model over which everything is defined. Then, the following family is a regular group configuration (resp. a definable group configuration) over $M$:

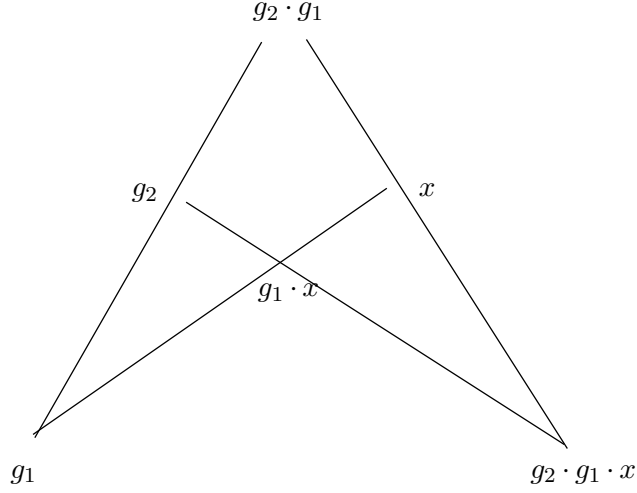
\begin{center}
    \begin{tikzcd}

&& g_2\cdot g_1 \arrow[llddddd, dash, shorten= 2 mm, shift right = 0.5 ex] \arrow[rrddddd, dash, shorten= 1 mm, shift right = 0.15 ex]

\\
 
\\
 & \,\,g_2\,\, & & \,\,x\,\, 
\\
&&    g_1 \cdot x

\\

\\

g_1 \arrow[rrruuu, dash, shift left = 1.9 ex, shorten= -0 mm]&&&& g_2\cdot g_1 \cdot x\arrow[llluuu, dash, shift right = 0.9 ex, shorten= -0.0 mm]

    \end{tikzcd}
\end{center}

Such a quadrangle is called a ``group configuration for $(G, X)$''. Note that this quadrangle is generically stable  over $A$ if and only if the tensor product $p^{\otimes 2} \otimes q$ is generically stable over $A$.

\end{prop}

\begin{proof}
The algebraicity relations are clear in the case of a free action. In fact, let us deal with the more subtle case of a faithful transitive action. Since any free action is in particular faithful, our proof will also include a proof of the independence relations for the case of a free action.

First, note that the type $tp(g_1, g_2, x / M)$ is definable, since it is the tensor product $p\otimes p \otimes q$. Then, by Lemma \ref{definissabilite_et_acl}, the type over $M$ of the sextuple is definable. This is the first point of the definition. The third point is easy to check, since $b_1 = g_1$, $b_2 = g_2$ and $b_3 = g_2 \cdot g_1$. 

Let us now prove the second point. First, note that $tp(g_1\cdot x / M) = tp(g_2 \cdot g_1 \cdot x / M) = tp(x / M) = q$, because $Stab(q) = G$, and $x \models q|_{M g_1 g_2}$. Similarly, $tp(g_2\cdot g_1 / M) = p$. So, to prove the second point, we may use stationarity and commutativity for generically stable types. Let us prove that $(g_2 \cdot g_1, \, g_1 \cdot x, \, g_2 \cdot g_1 \cdot x)$ realizes the tensor product $p \otimes q \otimes q |_M$, the other cases being similar. By saturation of $M$, let $x_0 \in X(M)$, and let $f: G \rightarrow X$ be the definable map $g \mapsto g \cdot x_0$. By transitivity of the action,  this map is onto. Since $p$ is generic in $G$, the type $f_*(p)$ is generic in $X$. Thus, uniqueness of the generic of $X$ (Proposition \ref{prop_gen_stable_generics_in_a_space} (1)) implies that $f_*(p) = q$. %Now, using this property, it is easy to check that $g_2 \cdot g_1 \cdot x \downfree_M  g_1 \cdot x$. It remains to show that $g_2 \cdot g_1 \downfree_M g_1 \cdot x \flex g_2 \cdot g_1 \cdot x$. To that end, it suffices to prove that $g_2 \cdot g_1 \downfree_M g_1 \cdot x \flex g_2$. We know that $x \models q|_{M g_1 g_2}$, so, by genericity, $g_1 \cdot x \models q|_{M g_1 g_2}$, i.e. $g_1 \cdot x \downfree_M g_1 \flex g_2$. So, by Lemma \ref{lemma_switch} and symmetry, we get $g_1 \downfree_M g_1 \cdot x \flex g_2$. Thus, it follows by genericity that $g_2 \cdot g_1 \downfree_M g_1 \cdot x \flex g_2$, as desired.

We will prove that $g_2 \cdot g_1 \cdot x \downfree_M g_2 \cdot g_1 \,\,(1)$, and then show $g_1 \cdot x \downfree_M g_2 \cdot g_1 \textbf{\textasciicircum} g_2 \cdot g_1 \cdot x \, \, (2)$. Using stationarity of generically stable types, this implies the result. We know that $x \models q|_{M g_1 g_2}$. By genericity of $q$, we then have $g_2 \cdot g_1 \cdot x \models q|_{M g_1 g_2}$, so $g_2 \cdot g_1 \cdot x \downfree_M g_2 \cdot g_1$, which is $(1)$. Let us now prove $(2)$. We know that $x \models q|_{M g_1 g_2}$, so there exists $g \models p|_{M g_1 g_2}$ such that $f(g) = x$.  By genericity, we have $g_1 \downfree_M g_2 \cdot g_1$. Then, using Lemma \ref{lemma_switch} and symmetry, we deduce from $g \downfree_M g_1 \textbf{\textasciicircum} g_2\cdot g_1$ that $g_1 \downfree_M g_2 \cdot g_1\textbf{\textasciicircum} g$. Then, since $Stab^r(p) = G$ (by Lemma \ref{lemma_gen_stable_generics_in_a_group}), we have  $g_1 \cdot g \models p|_{M g_2 \cdot g_1\textbf{\textasciicircum} g}$, i.e. $g_1 \cdot g \downfree_M g_2 \cdot g_1 \textbf{\textasciicircum}  g$, so $f(g_1 \cdot g) \downfree_M g_2 \cdot g_1 \textbf{\textasciicircum} f(g_2 \cdot g_1 \cdot g)$. In other words, we have $g_1 \cdot x \downfree_M  g_2 \cdot g_1 \textbf{\textasciicircum}  g_2 \cdot g_1 \cdot x$, which is (2).

Finally, let us prove the fourth point of the definition of a definable group configuration. The first part follows from definability over $M$ of the group action. Let us now prove the statements dealing with the canonical bases. We follow \cite{pillay} (Chapter 5, Remark 4.1). By Proposition \ref{prop_description_base_canonique_af_c(a)}, to prove that, say $g_1$ is interalgebraic over $M$ with $Cb(x, g_1 \cdot x / acl(M g_1))$, it suffices to prove that $g_1$ is interalgebraic over $M$ with its germ at $q$. Let us show that, in fact, for any $h \in G$, the element $h$ is definable over $M  \textbf{\textasciicircum} [h]$. 

Note that, since $Stab(q)=G$, composition of germs at $q$ of elements of $G$ is well-defined, and yields a type-definable group of germs $\Gamma$, such that $g \in G \mapsto [g] \in \Gamma$ is an $A$-definable group homomorphism. Moreover, by the strong germ property (Corollary \ref{germs_forts_cas_gen_stable}), the generic action of $G$ on $q$ induces a generic action of $\Gamma$ on $q$. So, it suffices to prove that, if $g \in G$ is such that $[g] = 1$, then $g=1$. 

Let $g$ be such an element. Let $a \in X$ be arbitrary. We want to show that $g \cdot a=a$, which will imply that $g=1$. Let $b \models q|_{Mg}$. By transitivity, let $h \in G$ be such that $h \cdot b = a$. Let $\gamma \models p|_{M g\textbf{\textasciicircum} h\textbf{\textasciicircum} a \textbf{\textasciicircum}b}$. Then, $\gamma \cdot h$ and $\gamma \cdot g \cdot h$ are generic over $M b$. So, by commutativity for generically stable types, we have $b \models q|_{M \gamma \cdot h}$ and $b \models q|_{M \gamma \cdot g \cdot h}$. Now, since germs can be composed, and $[g] = 1$, we have $[\gamma \cdot g \cdot h] = [\gamma \cdot h]$. Also, by the strong germ property, since $b \models q|_{M \gamma \cdot h}$ and $b \models q|_{M \gamma \cdot g \cdot h}$, we have $[\gamma \cdot h]\cdot b = (\gamma \cdot h) \cdot b$ and  $[\gamma \cdot g \cdot h]\cdot b = (\gamma \cdot g \cdot h) \cdot b$. As  $[\gamma \cdot g \cdot h] = [\gamma \cdot h]$, we deduce that $\gamma \cdot g \cdot h \cdot b = \gamma \cdot h \cdot b$, i.e. $\gamma \cdot g \cdot a = \gamma \cdot a$. So $g \cdot a = a$. As $a$ was arbitrary, faithfulness implies that $g=1$, as desired.
\end{proof}

In fact, such a configuration captures the structure of the group and its action, up to some notion of correspondence.

\begin{prop}\label{prop_recovering_group_from_config}
Let $(G, X)$ and $(H, Y)$ be type-definable transitive faithful actions, where $G, H$ and $X, Y$ are connected, type-definable with generically stable generics. Let $M$ be a sufficiently saturated model over which everything is defined. Let $(g_1, g_2, g_2 \cdot g_1, g_2 \cdot g_1 \cdot x, x, g_1 \cdot x)$ and $(h_1, h_2, h_2 \cdot h_1, h_2 \cdot h_1 \cdot y, y, h_1 \cdot y)$ be configurations built as in Proposition \ref{prop_how_to_build_quadrangles}, for $(G, X)$ and $(H, Y)$ respectively. Assume that these configurations are equivalent over $M$. Then, there exist type-definable sets $S \leq G \times H$ and $T \subseteq X \times Y$, and finite normal subgroups $N_1 \lhd G$, $N_2 \lhd H$ such that: 

\begin{enumerate}
\item The projection of the subgroup $S$ to $G/N_1 \times H/N_2$ is the graph of a group isomorphism $G/ N_1 \simeq H / N_2$.

\item The set $T$ is an $S$-invariant  finite-to-finite surjective correspondence between $X$ and $Y$.

\end{enumerate}

\end{prop}

\begin{proof}
Let $C \subset M$ be a small algebraically closed set of parameters over which everything is defined, and which captures the interalgebraicities. Thus, we have the following configurations, which are equivalent over $C$: 

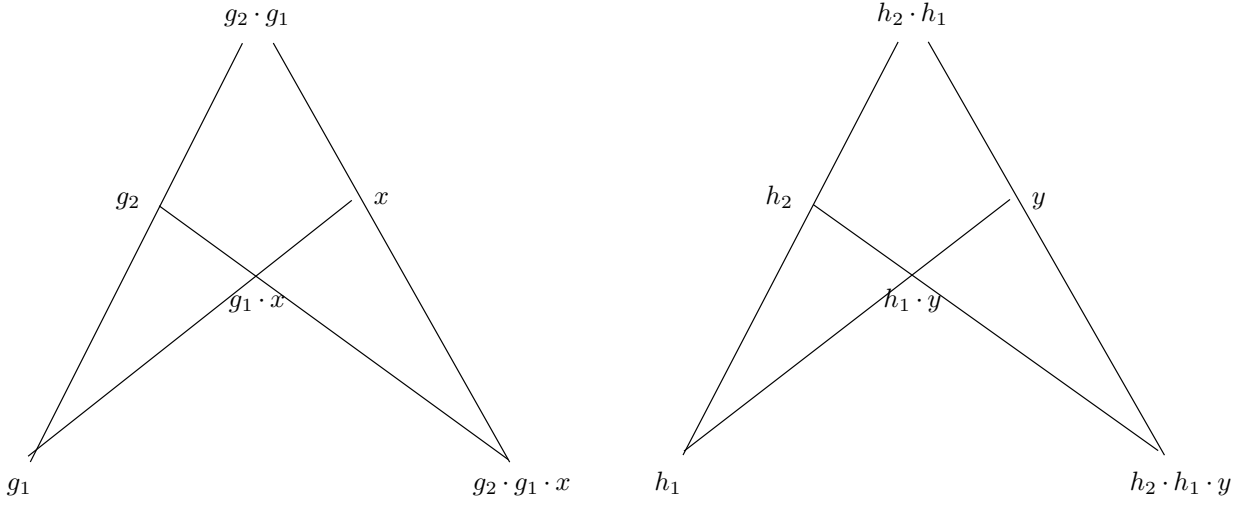
\small

\begin{center}
    \begin{tikzcd}

&& g_2 \cdot g_1 \arrow[llddddd, dash, shorten= 1 mm] \arrow[rrddddd, dash, shorten= 1 mm] &&

&

&&h_2 \cdot h_1 \arrow[llddddd, dash, shorten= 1 mm] \arrow[rrddddd, dash, shorten= 1 mm] &&   

\\
 
\\
 & g_2 & & x &    
 
 &
 
  & h_2 & & y &    

\\
&&g_1 \cdot x &&  

&

&& h_1 \cdot y&&

\\

\\

g_1 \arrow[rrruuu, dash, shift left = 1.4 ex]&&&& g_2 \cdot g_1 \cdot x\arrow[llluuu, dash, shift right = 1.0 ex]  

&

h_1 \arrow[rrruuu, dash, shift left = 1.4 ex]&&&& h_2 \cdot h_1 \cdot y\arrow[llluuu, dash, shift right = 1.2 ex, shorten= -0 mm]

    \end{tikzcd}
\end{center}

\normalsize

For $1 \leq i \leq 2$, let $c_i=(g_i, h_i) \in G \times H$. Also, let $c_3=c_2 \cdot c_1 = (g_2 \cdot g_1, h_2 \cdot h_1) \in G \times H$. Let $p_i=tp(c_i / M)$, for $1 \leq i \leq 3$. \begin{claim} The $p_i$ are generically stable over $acl(C) = C$. \end{claim}
\begin{proof}
Recall that, by construction, we have $tp(g_1 / M) = tp(g_2/M) = tp(g_2 \cdot g_1 / M)$, and this type is the unique generically stable generic of $G$. The last equality follows from the fact $g_2 \in G = Stab(tp(g_1 / M))$ and $tp(g_1/M g_2)$ is the generic of $G$.  So, this type is generically stable over $C = acl(C)$. Then, by interalgebraicity and Proposition \ref{prop_gen_stabilite_et_acl}, the types $p_i$ are generically stable over $acl(C) = C$.
\end{proof}

By assumption and interalgebraicity, we have $c_2 \downfree_M c_1$ and $ c_3 \downfree_M c_1$. Moreover, by definition, we have $c_2 \cdot c_1 = c_3$.  Thus, we have $c_1=(g_1, h_1)  \in Stab^r(p_2, p_3)$. Let us also define $S=Stab^r(p_2)$ and $Z = Stab^r(p_2, p_3)$. Let $c^{'}_1  \in M$ be such that $c^{'}_1 \equiv_C c_1$. So, there exist $g \in G(M)$ and $h \in H(M)$ such that $c^{'}_1 = (g,h) \in G \times H$.

\begin{claim}The following equality of type-definable sets holds:  $ S \cdot c^{'}_1 = Z$.

\end{claim}

\begin{proof} Let $s \in S$. Let $\gamma$ realize $p_2|_{C c^{'}_1 s}$. Then, by definition of $S$ as a stabilizer, we have $\gamma \cdot s \models p_2|_{C c^{'}_1 s}$. Moreover, we know that $c_1 \in Stab^r(p_2, p_3)$, so $c^{'}_1$ is in $Stab^r(p_2, p_3)$ as well. Thus, as $\gamma \cdot s$ realizes $ p_2|_{C c^{'}_1 s}$, we deduce that $  \gamma\cdot s\cdot c^{'}_1 $ realizes $ p_3|_{C c^{'}_1 s}$. Since we chose $\gamma$ realizing  $p_2|_{C c^{'}_1 s}$, we can conclude that $s \cdot c^{'}_1 \in Stab^r(p_2, p_3) = Z$. As $s \in S$ was arbitrary, we have just proved that $S \cdot c^{'}_1 \subseteq Z$. The other inclusion is proved similarly, by picking an arbitrary element $\alpha$ in $Z$, and letting the product $\alpha \cdot {c^{'}_1}^{-1}$  act by right-translation on (a realization of) $p_2$. \end{proof}

Let $\pi: G \times H \rightarrow G$ be the canonical projection. Since $\pi$ is a group morphism, the claim implies that $\pi(S) \cdot \pi(c^{'}_1) = \pi(Z)$. In particular, as $c_1=(g_1, h_1) \in Z$, we get ${g}^{-1} \cdot g_1 \in \pi(S)$. However, by construction, $g_1$ is, in the group $G$, generic over $M$. Since $g$ is in $M$, this implies that ${g}^{-1} \cdot g_1$ is also generic over $M$. Therefore, the $M$-type-definable subgroup $\pi(S) \leq G$ contains an element generic over $M$. So, the generic of $G$ concentrates on $\pi(S)$. Since, by Lemma \ref{lemma_gen_stable_generics_in_a_group} (4), any element of $G$ is the product of two realizations of this generic type, we deduce that $\pi(S)$ is equal to the whole group $G$. Symmetrically, the projection of $S$ to the second coordinate is equal to $H$.

Let $N_1 = \lbrace n_1 \in G \, | \, (n_1,1) \in S \rbrace$ and $N_2 = \lbrace n_2 \in H \, | \, (1, n_2) \in S \rbrace$. We then have $N_1 \times \lbrace 1 \rbrace \lhd S$, so $N_1 = \pi(N_1 \times \lbrace 1 \rbrace) \lhd \pi(S) = G$. Similarly, we have $N_2 \lhd H$. Thus, projecting the subgroup $S \leq G \times H$, we get a subgroup $\Sigma \leq G/N_1 \times H/N_2$, which is the graph of a definable group isomorphism $G/N_1 \rightarrow H/N_2$. Indeed, since $S$ projects surjectively onto $G$ and $H$, the subgroup $\Sigma$ projects surjectively onto $G/N_1$ and $H/N_2$.

\begin{claim}The groups $N_1$ and $N_2$ are finite. \end{claim}
\begin{proof}
 Let us show that $N_1$ is finite; the proof for $N_2$ will be symmetric. Let $n_1 \in N_1(M)$. Then, we know that $c_2=(g_2, h_2)$ realizes $p_2|_{C n_1}$. Then by definition of $S$, the element $(g_2\cdot n_1, h_2)$ realizes $p_2|_{Cn_1}$. In particular, we have $ g_2\cdot n_1\in acl(C h_2) $. Thus, since we are in a group, we deduce $n_1 \in acl(C c_2)$. However, we know that $c_2 \downfree_C n_1$. So, by Fact \ref{fact_forking} (5), we have $n_1 \in acl(C)$. As $n_1 \in N_1$ was arbitrary, we deduce by compactness that $N_1$ is finite. 
\end{proof}

Now, by saturation of $M$, let $(x_0, y_0) \in X(M) \times Y(M)$ be such that $(x_0, y_0) \equiv_C (g_1 \cdot x, h_1 \cdot y) $. Let us consider the following type-definable set: $T= S \cdot (x_0, y_0) \subseteq X \times Y$. Let $q_3 = tp(g_1 \cdot x, h_1 \cdot y / M)$ and $q_1 = tp(g_2 \cdot g_1 \cdot x, h_2 \cdot h_1 \cdot y / M)$. Then, since $tp(g_1 \cdot x / M)$ and $tp(g_2 \cdot g_1 \cdot x /M)$ are generically stable over $C$, Proposition \ref{prop_gen_stabilite_et_acl} implies that $q_1$ and $q_3$ are generically stable over $C=acl(C)$. Also, by hypothesis, we have $(g_2, h_2) \in Stab(q_3, q_1)$.

%From now on, for ease of notation, we will, somewhat abusively, identify $G$ with $G / N_1$ and with $H / N_2$. We also assume that, for $i=1, 2$, we have $g_i = h_i$. Let $q_2=tp(xy / M)$ and $q_3 = tp(g_1 \cdot x, h_1 \cdot y / M)$. By saturation, let $(x_0, y_0) \in M$ realize $q_2|_C$. Note that this implies that $acl(C x_0) = acl(C y_0)$. Letting $G$ act diagonally on $X \times Y$, we consider the type-definable set $ T:= G\cdot (x_0, y_0) \subseteq X \times Y$.

%\begin{claim}
%The stabilizer in $G$ of the generically stable type $q_2=tp(x_0, y_0 / C) = tp(x, y / C)$ is $G$ itself.
%\end{claim}

%\begin{proof}
%Just as above, consider the type-definable set $Stab(q_2, q_3)$. Then, by assumption, $g_1$ belongs to $Stab(q_2, q_3)$. As above, picking some $g \in G(M)$ such that $g \equiv_C g_1$, we have $g \cdot Stab(q_2) = Stab(q_2, q_3)$. Then, the element $g^{-1} \cdot g_1$ is in $Stab(q_2)$. Since $g_1$ is generic in $G$ over $M$, so is $g^{-1} \cdot g_1$. Thus, $Stab(q_2)$ is a  $C$-type-definable subgroup of the connected group $G$ containing a generic type. So $Stab(q_2) = G$, as stated.
%\end{proof}

\begin{claim}
The set $T$ is closed under the action of $S \leq G \times H$. Moreover, for all $x_1 \in X, \, y_1 \in Y$, the sets $T \cap (\lbrace x_1 \rbrace \times Y)$ and $T \cap  (X \times \lbrace y_1 \rbrace)$ are finite and nonempty.
\end{claim}

\begin{proof}

Since $T$ is the $S$-orbit of a point in the space $X \times Y$, it is closed under the action of $S$. Moreover, we have proved that $S$ projects onto $G$ and onto $H$. Therefore, by transitivity of the actions of $G$ on $X$ and $Y$, the projections $T \rightarrow X$ and $T \rightarrow Y$ are onto.

Thus, it remains to prove finiteness of the fibers, so to speak. For any $y \in Y$, let $T_y$ denote $T \cap (X \times \lbrace y \rbrace)$.  By symmetry, it suffices to prove that, for any $y_1 \in Y$, the type-definable set $T_{y_1}$ is finite. Note that, if $y_1, y_2 \in Y$, if $\gamma \in S$ is an element such that $\gamma \cdot y_1 =y_2$, then $\gamma \cdot T_{y_1} \subseteq T_{y_2}$. Then, since we also have $\gamma^{-1} \cdot y_2 = y_1$, the equality $\gamma \cdot T_{y_1} = T_{y_2}$ holds. So, by transitivity of the action of $G$ on $Y$, it suffices to prove that $T_{y_1}$ is finite, for \emph{some} $y_1 \in Y$. We shall prove that $T_{y_0}$ is finite.

By compactness and saturation of $M$, it suffices to prove that, if $(a, y_0) \in T_{y_0}(M)$, then $a \in acl(C y_0)$. So, let $a \in X(M)$ be such that $(a, y_0) \in T_{y_0}$. Let us prove that $a \in acl(C y_0)$. By saturation of $M$, and definition of $T$, there exists $(g, h) \in S(M)$ such that $(a, y_0) = (g \cdot x_0, h \cdot  y_0)$. In particular, we have $h\cdot y_0 = y_0$. Now, recall that $Stab^r(p_2) = S$, and $(g_2, h_2) \models p_2|_M$. Therefore, we have $(g_2 \cdot g, h_2 \cdot h) \models p_2|_M$. In particular, as $(g_2, h_2) \in Stab(q_3, q_1)$, the element $c:=(g_2 \cdot g, h_2 \cdot h)$ is also in $Stab(q_3, q_1)$. Also, by commutativity for the generically stable types $q_3|_C$ and $p_2|_C$, we have $(x_0, y_0) \models q_3|_{C c}$. Thus, $c \cdot (x_0, y_0) \models q_1|_{C c}$. In particular, we have $g_2 \cdot g \cdot x_0 \in acl(C h_2 \cdot h \cdot y_0)$. So $g\cdot x_0 \in acl(C g_2, h_2 \cdot h \cdot y_0) \subseteq acl(C g_2, h_2, h\cdot y_0)$. Now, recall that $acl(C g_2) = acl(C h_2)$, that $h\cdot y_0 = y_0$, and $a=g\cdot x_0$. Thus, we have $a \in acl(C g_2, y_0)$. Since $a, y_0 \in M$ and $g_2$ is generic over $M$, we have $g_2 \downfree_C a y_0$, so $g_2 \downfree_{C y_0} a$. Therefore, by Fact \ref{fact_forking} (5), we have  $a \in acl(C y_0)$, as desired.
\end{proof}

Thus, we can use the set $T$ to define an $S$-equivariant  finite correspondence between $X$ and $Y$.
\end{proof}

\begin{prop}

Any generically stable regular group configuration over $A$ is a generically stable group configuration over $A$.

\end{prop}

\begin{proof}
Let $(a_1, a_2, a_3, b_1, b_2, b_3)$ be a  generically stable regular group configuration over $A$. Let us show that it is a  generically stable group configuration over $A$.
Looking at the definitions, the only properties to check are those regarding the canonical bases $Cb(a_j a_k / acl(A b_i))$, for $i,j,k$ such that $\lbrace i,j,k \rbrace = \lbrace 1,2,3 \rbrace$. By symmetry of the context, it suffices to show that $b_1$ and $Cb(a_2 a_3 / acl(A b_1))$ are interalgebraic over $A$.

First, let $N$ be a sufficently saturated model containing $A b_1$, such that $a_2 \downfree_{A b_1} N$, which exists by extensibility of $tp(a_2 / A b_1)$. Then, by Fact \ref{fact_gen_stable_abus}, the type $tp(a_2 / N)$ is generically stable over $A b_1$: it is in fact the unique nonforking extension of $tp(a_2 / A)$. Then, by Proposition \ref{prop_gen_stabilite_et_acl}, the type $tp(a_2 a_3 / N)$ is generically stable over $acl(A b_1)$, and stationarity yields the equality: $Cb(a_2 a_3 / acl(A b_1)) = Cb(a_2 a_3 / N)$.

On the one hand, it is clear that $Cb(a_2 a_3 / acl(A b_1)) \subseteq acl(A b_1)$. Conversely, let us show that $b_1$ is algebraic over $A \cup  Cb(a_2 a_3 / acl(A b_1))$. Let $C=A \cup Cb(a_2 a_3 / acl(A b_1))$. By assumption, we know that $b_1 \in acl (A a_2 a_3) \subseteq acl( C a_2 a_3)$. Moreover, since  the type $tp(a_2 a_3 / N)$ is generically stable over $acl(Ab_1)$, and $C \subset N$ contains $Cb(a_2 a_3 / acl(A b_1))$, the type $tp(a_2 a_3 / N)$  is generically stable over $C$. In particular, $a_2 a_3 \downfree_C N$. Thus, since $b_1 \in N$, we have $a_2 a_3 \downfree_C b_1$. Finally, since $b_1 \in acl(C a_2 a_3)$, we can apply Fact \ref{fact_forking} (5), which implies that $b_1 \in acl(C)$, as desired. \end{proof}

\begin{defi}{Equivalent quadrangles }

Let $A$ be a set of parameters. Let $(a_1,a_2, a_3, b_1, b_2, b_3)$ and $(a^{'}_1,a^{'}_2, a^{'}_3, b^{'}_1, b^{'}_2, b^{'}_3)$ be quadrangles. We say that these quadrangles are \textit{equivalent over $A$}, or \textit{interalgebraic over $A$} if, for $i= 1, 2, 3$, we have  $acl(A a_i) = acl(A a^{'}_i)$ and $acl(A b_i)=acl(A b^{'}_i)$.

\end{defi}

\begin{prop}\label{prop_quadrangle_equivalent}

If a quadrangle is equivalent over $A$ to a generically stable (regular) group configuration over $A$, then it is itself a generically stable (regular) group configuration over $acl(A)$.

\end{prop}

\begin{proof}
Indeed, Proposition \ref{prop_gen_stabilite_et_acl} implies that  generic stability is preserved.  Moreover,  the algebraicity relations  are preserved, and, in the case of generically stable regular group configurations, so are the independence relations, thanks to Fact \ref{fact_forking} (4). In the case of generically stable group configurations, the second point of the definition (i.e. definability of the relative types) is preserved thanks to Lemma \ref{definissabilite_et_acl}: let $\alpha, \beta$ be elements such that $tp(\alpha / A \beta)$ is $A$-definable, and let $\alpha^{'}, \beta^{'}$ be such that $acl(A \alpha) = acl(A \alpha^{'})$ and $acl(A \beta) = acl(A \beta^{'})$. Then, by Definition \ref{defi_definable_types} (5), there exists a model $M \supset A \beta$ and a type $p \in S(M)$ extending $tp(\alpha / M \beta)$ such that $p$ is $A$-definable. We may assume that $\alpha$ realizes $p$. Then, $tp(\alpha / M)$ is $acl(A)$-definable, so $tp(\alpha^{'} / M)$ is also definable over $acl(A)$. Since $M \supset acl(A \beta) = acl(A \beta^{'})$, we have indeed shown that $tp(\alpha^{'} / acl(A) \beta^{'})$ is definable over $acl(A)$.

Lemma \ref{definissabilite_et_acl} also implies the required properties for the canonical bases: let $(b_1, a_2, a_3)$ and $(b^{'}_1, a^{'}_2, a^{'}_3)$ be such that $acl(A b_1) = acl(A b^{'}_1)$, $acl(A a_2) = acl(A a^{'}_2)$ and $acl(A a_3) = acl(A a^{'}_3)$. Assume that $tp(a_2 a_3 / acl(A b_1))$ is definable, and $b_1$ is interalgebraic over $A$ with the canonical basis of this type. Then, applying Lemma \ref{definissabilite_et_acl} once, we deduce that the type $tp(a^{'}_2 a^{'}_3 / acl(A b^{'}_1))$ is definable. Let $C \subseteq acl(A b^{'}_1) = acl(A b_1)$ denote its canonical basis. 
\begin{claim}
We have $acl(AC) = acl(A b^{'}_1)$.

\end{claim}
\begin{proof}
Applying Lemma \ref{definissabilite_et_acl} again, the type $tp(a_2 a_3 / acl(A b_1))$ is definable over $acl(A C)$. So, by the property of $Cb(a_2 a_3 / acl(A b_1))$, we have $b_1 \in acl(AC)$. So $b^{'}_1 \in acl(AC)$. Since $C \subseteq acl(A b^{'}_1)$, we are done.
\end{proof} This concludes the proof.
\end{proof}

\begin{notation}
If $p(x,y)$ is a complete type in several variables, where $x,y$ are tuples of variables, we let $p_{x}$ denote the restriction of $p$ to the tuple of variables $x$. 
\end{notation}

%Then, $p_x$ is a complete type in the tuple of variables $x$.

\begin{prop}\label{prop_descente_pour_quadrangles}

Let $(a_1,a_2,a_3,b_1,b_2,b_3)$ be a generically stable (regular) group configuration over a set $A$, and $A_0 \subseteq A$. Assume that the type  $tp(a_1a_2a_3b_1b_2b_3/A)$ is generically stable over $A_0$. Then $(a_1,a_2,a_3,b_1,b_2,b_3)$ is a  generically stable (regular) group configuration  over $A_0$.

\end{prop}

\begin{proof}

The interalgebraicity relations are proved using Proposition \ref{prop_b_in_acl(Aa)_cas_gen_stable}. Let us prove one of the independence relations, say $a_1 a_2 \downfree_{A_0} b_1$. We know that $a_1 a_2 \downfree_A b_1$, that $a_1 a_2 \downfree_{A_0} A$ and that $tp(a_1 a_2 / A_0)$ is generically stable. So, by Proposition \ref{transitivite_gen_stable}, we have $a_1 a_2 \downfree_{A_0} A b_1$, which implies that $a_1 a_2 \downfree_{A_0} b_1$, as desired. Note that this independence relation, along with stationarity of generically stable types, imply that $tp(a_1 a_2 / A_0 b_1)$ is $A_0$-definable.  
\end{proof}

\end{subsection}
\begin{subsection}{Statement of the theorem}

\begin{theo}\label{theo_config_groupe_faible}
Let $M$ be a $|\mathcal{L}|^+$-saturated model. Let $(a^0_1, a^0_2, a^0_3, b^0_1, b^0_2, b^0_3)$ be a  generically stable group configuration over $M$. Let $p_0(x_1, x_2, x_3, y_1, y_2, y_3)=tp(a^0_1, a^0_2, a^0_3, b^0_1, b^0_2, b^0_3/ M) $. Let $C_0 \subset M$ be the canonical basis of the type $p_0^{}$. Then, there exists a type-definable group $\Gamma$ acting transitively, faithfully and definably on a type-definable set $X$, elements $b^{''}_1, b^{''}_2, b^{''}_3, b^{'}_3 \in M$, and elements  $g_1, g_2 \in \Gamma$, $x \in X$, whose types over $M$ are generic, such that: 

\begin{enumerate}    \item The tuple $  b^{'}_1 b^{''}_2 b^{'''}_3$ realizes $((p_0)_{y_1, y_2 , y_3})|_{C_0}$.    

\item The group $\Gamma$ is connected and type-definable over $acl(C_0 b^{'}_1 b^{''}_2 b^{'''}_3)$. The space $X$ is connected and type-definable over $acl(C_0 b^{'}_1 b^{''}_2 b^{'''}_3)$.    

\item The following is a generically stable group configuration over $M$ which is equivalent, over $M$, to the quadrangle $(a^0_1, a^0_2, a^0_3, b^0_1, b^0_2, b^0_3)$: 

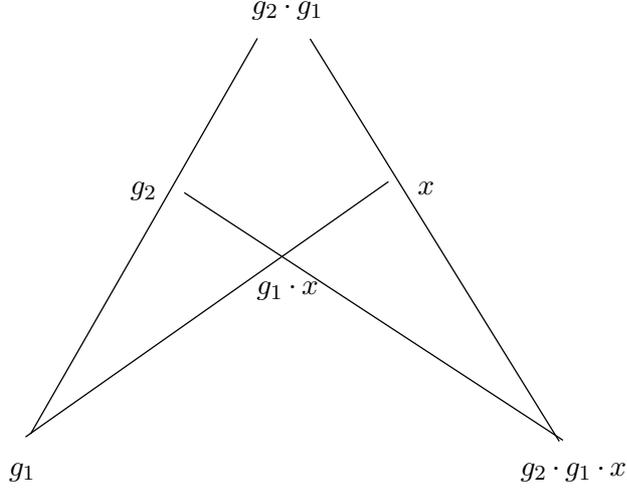
\begin{center}
    \begin{tikzcd}

&& g_2\cdot g_1 \arrow[llddddd, dash, shorten= 2 mm, shift right = 0.8 ex] \arrow[rrddddd, dash, shorten= 1.5 mm, shift right = -0.15 ex]

\\
 
\\
 & \,\,g_2\,\, & & \,\, x
\\
&& g_1\cdot x

\\

\\

g_1 \arrow[rrruuu, dash, shift left = 1.9 ex, shorten= -1 mm]&&&& g_2 \cdot g_1\cdot x\arrow[llluuu, dash, shift right = 1.4 ex, shorten= -0.7 mm]

    \end{tikzcd}
\end{center}

\end{enumerate}
Moreover, let $(R)$ denote the property ``$p_0$ is a generically stable \emph{regular} group configuration over $M$''. Then, if $(R)$ holds, the action of $\Gamma$ on $X$ can be assumed to be regular, i.e. transitive and free, instead of merely being faithful.

\end{theo}

Note that, in particular, $\Gamma$ and $X$ are type-definable over $M$, which is enough information if one does not need to control parameters.

\end{subsection}

The following proof is adapted from that of Elisabeth Bouscaren \cite{Bouscaren1989TheGC}, with ideas from \cite{pillay} (Chapter 5, Remark 1.10 and Theorem 4.5)  for the general case of a faithful transitive action.

The proof is divided into two steps: first, we find a group configuration which is equivalent to the original one, where some algebraicity relations have been replaced by definability. This is the content of Proposition \ref{prop_WMA_dcl} below. Then, using this stronger property, we consider some definable bijections permuting elements in (copies of) the new group configuration, and build a type-definable group from the germs of such maps. This is done in Section 3.

\begin{subsection}{Replacing algebraicity with definability}

The goal of this subsection is the following proposition, which enables us, at the cost of enlarging the basis, to replace some of  the algebraicity  in the configuration by definability, while keeping a generically stable type. This will then enable us to consider definable bijections that permute elements of the configuration.

\begin{prop}\label{prop_WMA_dcl}
Under the hypotheses of Theorem \ref{theo_config_groupe_faible}, there exist elements $b^{'}_1, b^{''}_2, b^{'''}_3 \in M$ and a configuration $(a_1, a_2, a_3, b_1, b_2, b_3)$ equivalent over $M$ to $(a^0_1, a^0_2, a^0_3, b^0_1, b^0_2, b^0_3)$, such that \begin{enumerate}
\item The tuple $  b^{'}_1 b^{''}_2 b^{'''}_3$ realizes $((p_0)_{y_1, y_2 , y_3})|_{C_0}$. 
\item The type $tp(a_1, a_2, a_3, b_1, b_2, b_3 / M)$ is generically stable over $C := acl(C_0  b^{'}_1 b^{''}_2 b^{'''}_3)$.

\item We have $a_1 \in dcl(C a_2 b_3)$, $a_2 \in dcl(C a_3 b_1)$, and $a_3 \in dcl(C a_1 b_2) \cap dcl(C a_2 b_1)$.
\end{enumerate}

\end{prop}

\begin{rem}
Note that, if the property (R) in Theorem \ref{theo_config_groupe_faible} held, the proposition would yield a configuration equivalent over $M$ to a  generically stable \emph{regular} group configuration, so still a \emph{regular} group configuration over $M$, generically stable over $acl(C_0  b^{'}_1 b^{''}_2 b^{'''}_3 )$. Then, by Proposition \ref{prop_descente_pour_quadrangles}, the new configuration would be a generically stable \emph{regular} group configuration over $acl(C_0  b^{'}_1 b^{''}_2 b^{'''}_3)$.
\end{rem}

The following result will be useful in this subsection.

\begin{lemma} \label{lemma_interalg_for_codes}
Let $A$, $B$, $B^{'}$ be parameter sets, and $a$ be an element such that $a \in acl(B) \cap acl(B^{'})$. Let $\alpha$ be the code of the set of conjugates of $a$ over $A BB^{'}$. Assume that $B \downfree_{Aa} B^{'}$. Then, $a$ and $\alpha$ are interalgebraic over $A$.

\end{lemma}

\begin{proof}
First, note that $a$ belongs to the finite set coded by $\alpha$, so $a \in acl(A \alpha)$. Let us prove that $\alpha$ is algebraic over $Aa$. We know that $\alpha$ codes a subset of the set of conjugates of $a$ over $B$, so $\alpha \in acl(B a)$. Similarly, $\alpha \in acl(B^{'}a)$. Thus, by Fact \ref{fact_forking} (5) and the hypothesis $B \downfree_{Aa} B^{'}$, we have $\alpha \in acl(Aa)$, as required.
\end{proof}

    Let us now prove the proposition.
    
    \begin{proof}[Proof of Proposition \ref{prop_WMA_dcl}]
    
    By saturation of $M$, let $b^{'}_1, b^{''}_2 \in M$ be such that $b^{'}_1 b^{''}_2 \equiv_{C_0} b^0_1 b^0_2$.

    \begin{claim}
    We have $a^0_1 b^0_2 a^0_3 b^{'}_1 \equiv_{C_0} a^0_1 b^0_2 a^0_3 b^0_1$ and $b^0_1 a^0_2 a^0_3 b^{''}_2 \equiv_{C_0} b^0_1 a^0_2 a^0_3 b^0_2$.
    
    \end{claim}
    \begin{proof}
    The type $tp(a^0_1 b^0_2 a^0_3 / M b^0_1)$ is the unique nonforking extension of the generically stable type $tp(a^0_1 b^0_2 a^0_3 / C_0)$. Indeed, we have $a^0_1 b^0_2 \downfree_{M}   b^0_1$, and $a^0_1 b^0_2 \downfree_{C_0} M $, so by transitivity for generically stable types (Proposition \ref{transitivite_gen_stable}), we have $a^0_1 b^0_2 \downfree_{C_0} M b^0_1$. Since $a^0_3 \in acl(C_0 a^0_1 b^0_2)$, we have $a^0_1 b^0_2 a^0_3 \downfree_{C_0} M b^0_1$. Thus, the type $tp(a^0_1 b^0_2 a^0_3 / M b^0_1)$ is generically stable over $C_0$, in particular it is $C_0$-invariant. As $b^{'}_1 \in M$ has the same type over $C_0$ as $b^0_1$, we have indeed $a^0_1 b^0_2 a^0_3 b^{'}_1 \equiv_{C_0} a^0_1 b^0_2 a^0_3 b^0_1$.

    The other result is proved similarly, using the fact that $tp(b^0_1 a^0_2 a^0_3 / M b^0_2)$ is generically stable over $C_0$, so $C_0$-invariant.
    \end{proof}
    
    Then, let $b^{'}_3, a^{'}_2$ and $b^{''}_3, a^{''}_1$ be such that  $a^0_1 b^0_2 a^0_3 b^{'}_1  a^{'}_2 b^{'}_3  \equiv_{C_0} a^0_1 b^0_2 a^0_3 b^0_1 a^0_2 a^0_3$ and $b^0_1 a^0_2 a^0_3 a^{''}_1 b^{''}_2 b^{''}_3 \equiv_{C_0} b^0_1 a^0_2 av_3 a^0_1 b^0_2 b^0_3$. In other words, the following are copies of the original configuration:

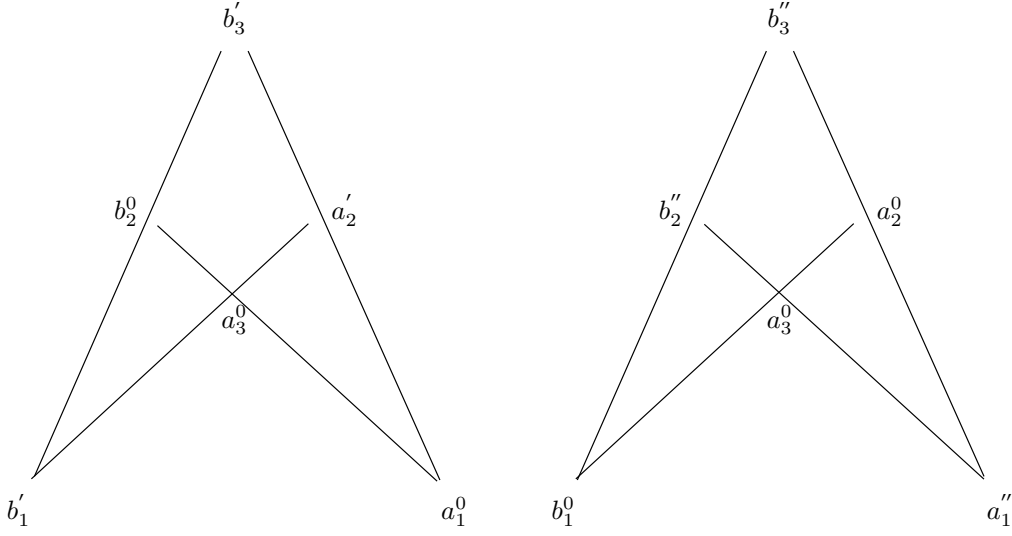
\small

\begin{center}
    \begin{tikzcd}

&& b^{'}_3 \arrow[llddddd, dash, shorten= 1 mm] \arrow[rrddddd, dash, shorten= 1 mm] &&

&

&& b^{''}_3  \arrow[llddddd, dash, shorten= 1 mm] \arrow[rrddddd, dash, shorten= 1 mm] &&   

\\
 
\\
 & b^0_2 & &  a^{'}_2  &    
 
 &
 
  &  b^{''}_2 & &  a^{0}_2  &    

\\
&& a^0_3 &&  

&

&& a^0_3 &&

\\

\\

 b^{'}_1  \arrow[rrruuu, dash, shift left = 1.4 ex]&&&& a^0_1 \arrow[llluuu, dash, shift right = 1.0 ex]  

&

b^0_1 \arrow[rrruuu, dash, shift left = 1.4 ex]&&&&  a^{''}_1 \arrow[llluuu, dash, shift right = 1.2 ex, shorten= -0 mm]

    \end{tikzcd}
\end{center}

\normalsize
    
    Then, let $\widetilde{a_1}$ be the code of the set of conjugates of $a_1^0$ over $C_0 b^0_2 a^0_3  a^{'}_2 b^{'}_3$. Similarly, let  $\widetilde{a_2}$ be the code of the set of conjugates of $a_2^0$ over $C_0 b^0_1 a^0_3  a^{''}_1 b^{''}_3$.
    
    \begin{claim}
    We have $acl(C_0 \widetilde{a_1}) = acl(C_0 a_1^0)$ and $acl(C_0 \widetilde{a_2}) = acl(C_0 a_2^0)$.
    
    \end{claim}
    
    \begin{proof}
    For the first interalgebraicity, we wish to apply Lemma \ref{lemma_interalg_for_codes} to $a = a^0_1$, $A = C_0$, $B = C_0 b^0_2 a^0_3$, and $B^{'} = C_0 b^{'}_3 a^{'}_2$. For the second one, we apply the same Lemma to $a = a^0_2$, $A = C_0$, $B = C_0 b^0_1 a^0_3$, and $B^{'} = C_0 b^{''}_3 a^{''}_1$. The only hypotheses that do not follow immediately from the constructions are the independence properties. To prove $b^0_2 a^0_3 \downfree_{C_0 a^0_1} b^{'}_3 a^{'}_2$ and $b^0_1 a^0_3 \downfree_{C_0 a^0_2} b^{''}_3 a^{''}_1$, we use the fact that the copies have the same type over $C_0$ as the original configuration. So, it suffices to prove $b^0_2 a^0_3 \downfree_{C_0 a^0_1} b^{0}_3 a^{0}_2$ and $b^0_1 a^0_3 \downfree_{C_0 a^0_2} b^{0}_3 a^{0}_1$. Using the algebraicity properties of the configuration, along with Fact \ref{fact_forking} (4)(b), it suffices to prove $b^0_2 \downfree_{C_0 a^0_1} b^{0}_3$ and $b^0_1 \downfree_{C_0 a^0_2} b^{0}_3$. By Fact \ref{fact_forking} (1), these follow from $b^0_2 \downfree_{C_0 }a^0_1 b^{0}_3$ and $b^0_1 \downfree_{C_0 }a^0_2 b^{0}_3$, which hold by Remark \ref{rem_defi_group_config} (1).
    \end{proof}
    
    The motivation for building these copies is that $\widetilde{a_1} \in dcl(C_0 b^0_2 a^0_3 b^{'}_3 a^{'}_2)$ and $\widetilde{a_2} \in dcl(C_0 b^0_1 a^0_3 b^{''}_3 a^{''}_1)$.
    
    Then, let $(\alpha_1, \alpha_2, \alpha_3, \beta_1, \beta_2, \beta_3)$ denote the following configuration:

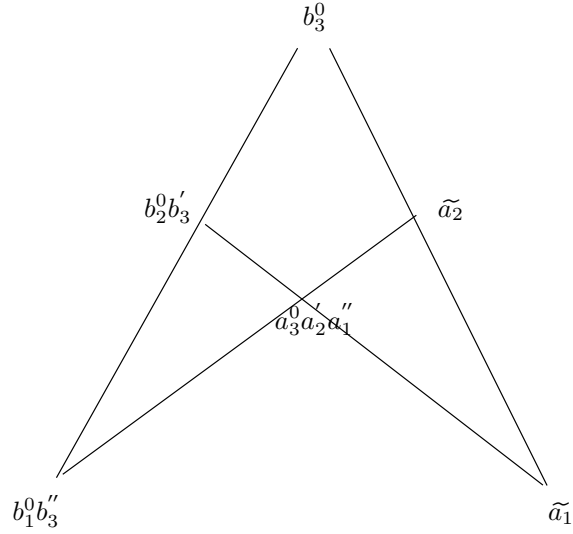
\small
\begin{center}
    \begin{tikzcd}

&& b^0_3 \arrow[llddddd, dash, shorten= 1 mm] \arrow[rrddddd, dash, shorten= 1 mm]

\\
 
\\
 & b^0_2 b^{'}_3 & & \widetilde{a_2}
\\
&& a^0_3 a^{'}_2 a^{''}_1  

\\

\\

b^0_1 b^{''}_3  \arrow[rrruuu, dash, shift left = 1.4 ex]&&&& \widetilde{a_1}\arrow[llluuu, dash, shift right = 1.0 ex]

    \end{tikzcd}
\end{center}
\normalsize

    Let $C_1 = acl(C_0 b^{'}_1 b^{''}_2) \subset M$.

    \begin{claim}
    The quadrangle $(\alpha_1, \alpha_2, \alpha_3, \beta_1, \beta_2, \beta_3)$ is equivalent over $C_1$ to $(a^0_1, a^0_2, a^0_3, b^0_1, b^0_2, b^0_3)$.
    
    \end{claim}
    \begin{proof}
    The only things left to check are $b^{'}_3 \in acl(C_1 b^0_2)$, $ b^{''}_3 \in acl(C_1 b^0_1)$,  and that $a^{'}_2 a^{''}_1$ is in $acl(C_1 a^0_3)$. These can be checked by inspecting the above copies of the configuration, keeping in mind that $C_1$ contains $C_0 b^{'}_1 b^{''}_2$.
    \end{proof}
    
    Thus, by Proposition \ref{prop_quadrangle_equivalent}, the tuple $(\alpha_1, \alpha_2, \alpha_3, \beta_1, \beta_2, \beta_3)$ is a generically stable group configuration over $C_1$. What we gained is that $\alpha_1 \in dcl(C_1 \alpha_3 \beta_2)$ and $\alpha_2 \in dcl(\alpha_3 \beta_1)$. Also note that $\beta_3 = b^0_3$.

    Now, let $\widetilde{\alpha_3}$ be the code of the set of conjugates of $\alpha_3$ over $C_1 \beta_1 \alpha_2 \beta_2 \alpha_1$. As above, using Lemma \ref{lemma_interalg_for_codes}, we have $acl(C_1 \alpha_3) = acl(C_1 \widetilde{\alpha_3})$. So, by Proposition \ref{prop_quadrangle_equivalent}, the quadrangle $(\alpha_1, \alpha_2, \widetilde{\alpha_3}, \beta_1, \beta_2, \beta_3)$ is a generically stable group configuration over $C_1$. By construction, we get $\widetilde{\alpha_3} \in dcl(C_1 \beta_1 \alpha_2 \beta_2 \alpha_1)$.

    \begin{claim}
    We still have $\alpha_1 \in dcl(C_1 \widetilde{\alpha_3} \beta_2)$ and $\alpha_2 \in dcl(\widetilde{\alpha_3} \beta_1)$.
    
    \end{claim}
    \begin{proof}
    Recall that  $\alpha_1 \in dcl(C_1 \alpha_3 \beta_2)$ and $\alpha_2 \in dcl(\alpha_3 \beta_1)$. So, let $j$ and $l$ be some $C_1$-definable maps such that $j(\alpha_3, \beta_2) = \alpha_1$ and $l(\alpha_3, \beta_1) = \alpha_2$. Then, for any element $\alpha$ which is a conjugate of $\alpha_3$ over  $C_1 \beta_1 \alpha_2 \beta_2 \alpha_1$, we have $j(\alpha, \beta_2) = \alpha_1$ and $l(\alpha, \beta_1) = \alpha_2$. As $\widetilde{\alpha_3}$ is the code of that set of conjugates, we are done.
    \end{proof}

    Now, let $\beta^{'''}_3 \in M$ realize $tp(\beta_3 / C_1) = tp(b^0_3 / C_1)$. Then, by $C_1$-invariance of the types $tp(\widetilde{\alpha_3} \beta_2 \alpha_1 / M \beta_3)$ and $tp(\alpha_3 \alpha_2 \beta_1 / M \beta_3)$, we have $\widetilde{\alpha_3} \beta_2 \alpha_1 \beta_3 \equiv_{C_1} \widetilde{\alpha_3} \beta_2 \alpha_1 \beta^{'''}_3$ and $\widetilde{\alpha_3} \alpha_2 \beta_1 \beta_3 \equiv_{C_1} \widetilde{\alpha_3} \alpha_2 \beta_1 \beta^{'''}_3$. So, let $\beta^{'''}_1, \alpha^{'''}_2$ and $\beta^{'''}_2, \alpha^{'''}_1$ be such that $\widetilde{\alpha_3}\beta_2 \alpha_1 \beta_3 \beta_1 \alpha_2  \equiv_{C_1} \widetilde{\alpha_3} \alpha_2 \beta_1 \beta^{'''}_3 \alpha^{'''}_2 \beta^{'''}_1 $ and $\widetilde{\alpha_3} \alpha_2 \beta_1 \beta_3 \beta_2 \alpha_1 \equiv_{C_1} \widetilde{\alpha_3} \alpha_2 \beta_1 \beta^{'''}_3 \beta^{'''}_2 \alpha^{'''}_1$.

    In other words, we have two copies of the configuration $(\alpha_1, \alpha_2, \widetilde{\alpha_3}, \beta_1, \beta_2, \beta_3)$:

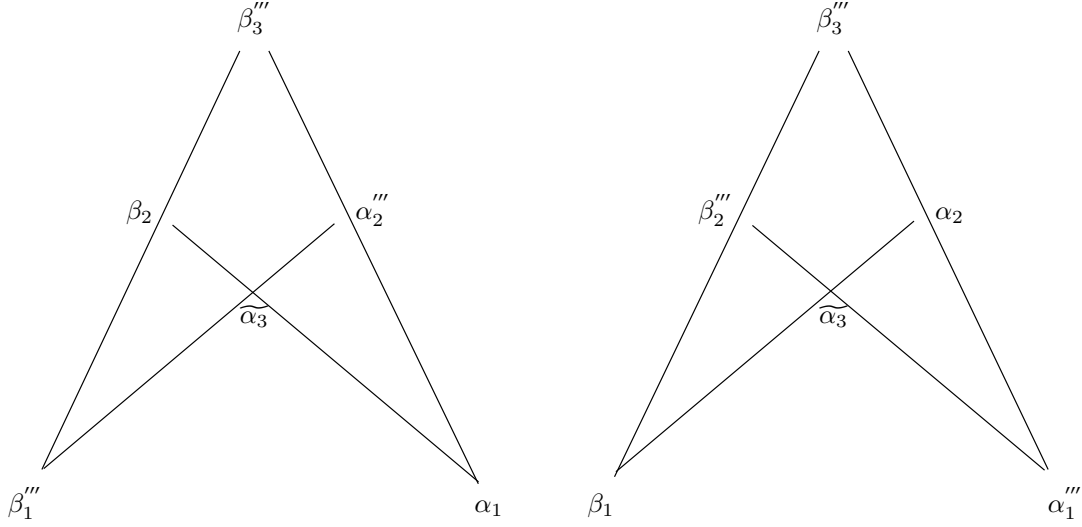
\small

\begin{center}
    \begin{tikzcd}

&& \beta^{'''}_3 \arrow[llddddd, dash, shorten= 1 mm] \arrow[rrddddd, dash, shorten= 1 mm] &&

&

&& \beta^{'''}_3  \arrow[llddddd, dash, shorten= 1 mm] \arrow[rrddddd, dash, shorten= 1 mm] &&   

\\
 
\\
 & \beta_2 & &  \alpha^{'''}_2  &    
 
 &
 
  &  \beta^{'''}_2 & &  \alpha_2  &    

\\
&& \widetilde{\alpha_3}&&  

&

&& \widetilde{\alpha_3} &&

\\

\\

 \beta^{'''}_1  \arrow[rrruuu, dash, shift left = 1.4 ex]&&&& \alpha_1 \arrow[llluuu, dash, shift right = 1.0 ex]  

&

\beta_1 \arrow[rrruuu, dash, shift left = 1.4 ex]&&&&  \alpha^{'''}_1 \arrow[llluuu, dash, shift right = 1.2 ex, shorten= -0 mm]

    \end{tikzcd}
\end{center}

\normalsize

Finally, the configuration $(a_1, a_2, a_3, b_1, b_2, b_3)$ we consider is $(\alpha_1 \alpha^{'''}_2, \alpha_2 \alpha^{'''}_1, \widetilde{\alpha_3}, \beta_1 \beta^{'''}_2, \beta_2 \beta^{'''}_1, \beta_3)$. We let $C$ denote $acl(C_1 \beta^{'''}_3) = acl(C_0 b^{'}_1 b ^{''}_2 \beta^{'''}_3)$.

    \begin{claim}
    We have $b^{'}_1 b ^{''}_2 \beta^{'''}_3 \models ((p_0)_{y_1, y_2, y_3})|_{C_0} = tp(b^0_1 b^0_2 b^0_3 / C_0)$.
    \end{claim}
    \begin{proof}
    We already know that $b^{'}_1 b ^{''}_2 \models ((p_0)_{y_1, y_2})|_{C_0}$. Then, by stationarity and by Remark \ref{rem_defi_group_config}, the type $tp(b^0_3 / C_0 b^{0}_1 b ^{0}_2 )$ is equal to $(p_0)_{y_3} |_{C_0 b^{0}_1 b ^{0}_2}$. Since $(p_0)_{y_3}$ is generically stable over $C_0$, and $\beta^{'''}_3$ realizes $(p_0)_{y_3}|_{C_0 b^{'}_1 b ^{''}_2}$, we conclude by $C_0$-invariance.
    \end{proof}

    \begin{claim}
    We have $a_3 \in dcl(C a_2 b_1) \cap dcl(C a_1 b_2)$.
    \end{claim}
    \begin{proof}
    Recall that $\widetilde{\alpha_3} \in dcl(C_1 \beta_1 \alpha_2 \beta_2 \alpha_1)$. So, in the copies above, we have $\widetilde{\alpha_3} \in dcl(C_1 \beta^{'''}_1 \alpha^{'''}_2 \beta_2 \alpha_1)$ and $\widetilde{\alpha_3} \in dcl(C_1 \beta_1 \alpha_2 \beta^{'''}_2 \alpha^{'''}_1)$. Thus, since $a_1 = \alpha_1 \alpha^{'''}_2$, $a_2 = \alpha_2 \alpha^{'''}_1$, $b_1 = \beta_1 \beta^{'''}_2$, and $b_2 = \beta_2 \beta^{'''}_1$, we get the result.
    \end{proof}
    
    \begin{claim}
    We have $a_1 \in dcl(C a_3 b_2)$ and $a_2 \in dcl(C a_3 b_1)$.
    \end{claim}
    
    \begin{proof}
    Recall that $a_1 = \alpha_1 \alpha^{'''}_2$. Also, we proved above that $\alpha_1 \in dcl(C_1 \widetilde{\alpha_3} \beta_2)$ and $\alpha_2 \in dcl(C_1 \widetilde{\alpha_3} \beta_1)$. So, in the first copy (the one on the left in the diagrams above), we also have $\alpha^{'''}_2 \in dcl(C_1 \widetilde{\alpha_3} \beta^{'''}_1)$. Also, we defined $b_2$ as $\beta_2 \beta^{'''}_1$. Thus, we get $\alpha_1 \alpha^{'''}_2 \in dcl(C_1 \widetilde{\alpha_3} \beta_2 \beta^{'''}_1)) = dcl(C_1 a_3 b_2) \subseteq dcl(C a_3 b_2)$.

    The other result is proved similarly: we have $\alpha^{'''}_1 \in dcl(C_1 \widetilde{\alpha_3} \beta^{'''}_2)$ and $\alpha_2 \in dcl(C_1 \widetilde{\alpha_3} \beta_1)$, so $a_2 = \alpha_2 \alpha^{'''}_1 \in dcl(C_1  \widetilde{\alpha_3} \beta^{'''}_2 \beta_1) = dcl(C_1 a_3 b_1) \subseteq dcl(C a_3 b_1).$
    \end{proof}

    \begin{claim}
    The quadrangle $(a_1, a_2, a_3, b_1, b_2, b_3)$ is equivalent over $C$ to $(a^0_1, a^0_2, a^0_3, b^0_1, b^0_2, b^0_3)$.
    
    \end{claim}
    \begin{proof}
    This amounts to checking that  $(a_1, a_2, a_3, b_1, b_2, b_3)$ is equivalent over $C$ to the configuration $(\alpha_1, \alpha_2, \widetilde{\alpha_3}, \beta_1, \beta_2, \beta_3)$. It relies on the following algebraicity properties:  $\alpha^{'''}_2 \in acl(C_1  \beta^{'''}_3 \alpha_1) = acl(C \alpha_1)$, $\alpha^{'''}_1 \in acl(C \alpha_2)$, $\beta^{'''}_2 \in acl(C \beta_1)$ and $\beta^{'''}_1 \in acl(C \beta_2)$.
    \end{proof}

    This finishes the proof, by setting $b^{'''}_3 = \beta^{'''}_3$.
    \end{proof}

\end{subsection}

\end{section}

\begin{section}{Constructing a group using germs of definable bijections}

In this section, construct an appropriate typ-definable group, and finish the proof of the theorem.

\begin{subsection}{Composition of germs}

The aim of this subsection is to build an appropriate group from germs of definable bijections. In this subsection and the next, the context is that of the conclusion of Proposition \ref{prop_WMA_dcl}.

\begin{defi}Let $f_{b_1}$ and $g_{b_2}$ be definable bijections sending respectively $a_2$ to $a_3$ and $a_3$ to $a_1$, where $f$ and $g$ are $C$-definable families of definable bijections. Let $h_{b_1 b_2}$ be the composite $g_{b_2} \circ f_{b_1}$.

\end{defi}

Then, the independence hypotheses on the configuration imply that, in the sense of Definition \ref{defi_f_defined_at_p}, the functions $f_{b_1}$ and $h_{b_1 b_2}$ are well-defined at $tp(a_2/C)$, and that the function $g_{b_2}$ is well-defined at $tp(a_3 / C)$. Moreover, we will show that the germs of these functions can be composed.

\begin{prop}\label{prop_f_b(a)_independant_de_Ab_cas_gen_stable}
Let $A, D \supseteq C$ be sets of parameters. Let $\beta_1, \beta_2, \alpha_1, \alpha_2, \alpha_3$ be realizations of $tp(b_1 / M)|_D$, $tp(b_2 / M)|_D$, $tp(a_1 / M)|_D$, $tp(a_2 / M)|_D$ and $tp(a_3 / M)|_D$ respectively. 
\begin{enumerate}
    \item If $\alpha_2 \downfree_D A\beta_1 $, then $f_{\beta_1}(\alpha_2) \downfree_D A\beta_1$.
    \item If $\alpha_3 \downfree_D A\beta_2 $, then $g_{\beta_2}(\alpha_3) \downfree_D A\beta_2$.
    \item If $\alpha_3 \downfree_D A \beta_1$, then $(f_{\beta_1})^{-1}(\alpha_3) \downfree_D A\beta_1$.
    \item If $\alpha_1 \downfree_D A\beta_2 $, then $(g_{\beta_2})^{-1}(\alpha_1) \downfree_D A\beta_2$.
\end{enumerate}
\end{prop}

\begin{proof}
Let us prove the first implication, the other ones being similar. Assume that $\alpha_2 \downfree_D A\beta_1 $. We know, by hypothesis, that $\alpha_2 \downfree_C D$. Moreover, the type $tp(\alpha_2/C)$ is generically stable. So, by transitivity, $\alpha_2 \downfree_C AD\beta_1$. Then, $\alpha_2 \downfree_{C\beta_1} AD $, so $f_{\beta_1}(\alpha_2) \downfree_{C\beta_1} AD$ $(*)$.

On the other hand, since $\alpha_2 \downfree_C AD\beta_1$, we have $\alpha_2 \downfree_C \beta_1$.
So, by stationarity, $\alpha_2 \beta_1 \equiv_C a_2 b_1$, so $f_{\beta_1}(\alpha_2)\alpha_2 \beta_1 \equiv_C f_{b_1}(a_2) a_2 b_1 = a_3 a_2 b_1$.
Since $a_3 \downfree_C b_1$, we deduce that $f_{\beta_1}(\alpha_2) \downfree_C \beta_1$. 

Recall that the type $tp(f_{\beta_1}(\alpha_2) / C) = tp(a_3 / C)$ is generically stable. So, by $(*)$ and transitivity, we have $f_{\beta_1}(\alpha_2) \downfree_{C} \beta_1 AD$. So, by monotonicity, we have indeed $f_{\beta_1}(\alpha_2) \downfree_{D} \beta_1 A$.
\end{proof}

%\begin{rem}These results are stated in terms of independence and stationarity. One may ask if they can be rephrased with definable types, writing the proofs in a different manner.\end{rem}

\begin{coro}\label{coro_composee_germs_bien_def_cas_gen_stable}
Let $\beta_1, \beta^{'}_1, \beta^{''}_1, \beta^{'''}_1$ be realizations of $tp(b_1 / C)$, $\beta_2, \beta^{'}_2$ be realizations of $tp(b_2 / C)$.

Then, the following germs are well-defined, i.e. only depend on the germs of the functions involved: $[g_{\beta_2}]^{-1} \circ [g_{\beta^{'}_2}]$, $[g_{\beta_2}] \circ [f_{\beta_1}]$, $[f_{\beta_1}]^{-1} \circ [f_{\beta^{'}_1}]$, $[f_{\beta_1}] \circ [f_{\beta^{'}_1}]^{-1}$, and 

\noindent $[f_{\beta_1}]^{-1} \circ [f_{\beta^{'}_1}] \circ [f_{\beta^{''}_1}]^{-1} \circ [f_{\beta^{'''}_1}]$.
\end{coro}

\begin{proof}
We shall use Proposition \ref{prop_f_b(a)_independant_de_Ab_cas_gen_stable}, taking realizations of the generically stable types involved, independent from all the parameters that appear.

Let us prove for instance that the germ of $(g_{\beta_2})^{-1} \circ g_{\beta^{'}_2}$ only depends on the germs $[g_{\beta_2} ]$ and $[g_{\beta^{'}_2}]$. Let $b_2^1, b_2^2$ be realizations of $tp(b_2 / C)$ such that $[g_{\beta_2} ] = [g_{b^1_2} ]$ and $[g_{\beta^{'}_2}]= [g_{b^2_2} ]$. Let us show that the germs of the functions $(g_{\beta_2})^{-1} \circ g_{\beta^{'}_2}$ and $(g_{b^1_2})^{-1} \circ g_{b^2_2}$ are equal.

Let $\alpha_3$ realize the type of $a_3$ over $C$, such that $\alpha_3 \downfree_C b_2^1 b_2^2 \beta_2 \beta^{'}_2$. Then, since $[g_{\beta^{'}_2} ] = [g_{b^2_2} ]$, we have $g_{\beta^{'}_2}(\alpha_3)= g_{b^2_2}(\alpha_3) $.  Moreover, by Proposition \ref{prop_f_b(a)_independant_de_Ab_cas_gen_stable} applied to $D=C$ and $A=b_2^1 \beta_2$, we know that $g_{b^2_2}(\alpha_3) \downfree_C b_2^1 \beta_2$. Since we have assumed that $[g_{\beta_2}]= [g_{b^1_2} ]$, we can deduce that $(g_{b^1_2}^{-1} \circ g_{b^2_2})(\alpha_3)=(g_{\beta_2}^{-1} \circ g_{b^2_2})(\alpha_3)= (g_{\beta_2}^{-1} \circ g_{\beta^{'}_2})(\alpha_3)$. As $\alpha_3 \downfree_C b_2^1 b_2^2 \beta_2 \beta^{'}_2 $, we have shown that the germs of the composites $(g_{\beta_2})^{-1} \circ g_{\beta^{'}_2}$ and $(g_{b^1_2})^{-1} \circ g_{b^2_2}$  are equal.
\end{proof}

\begin{defi}
Let $F$, resp. $G$, be the type-definable set of the germs of functions of the form $f_{\beta_1}$, resp. $g_{\beta_2}$, where $\beta_1$ realizes $tp(b_1 / C)$, resp. $\beta_2$ realizes $tp(b_2 / C)$.
%Then, by Fact \ref{fact_definissabilite_des_germs}, these sets are type-definable.
\end{defi}

\begin{rem}\label{rem_F_et_G_et_p_commutent}
By completeness of $tp(b_1/C)$ and $tp(b_2 / C)$, the partial types over $C$ defining $F$ and $G$ are in fact complete. By Proposition \ref{prop_gen_stabilite_et_acl} and Fact \ref{fact_gen_stable_abus}, these types are generically stable. Moreover, since the type $r=tp(a_1a_2a_3b_1b_2b_3/M)$ is generically stable over $C$, it commutes with itself. So, by Corollary \ref{coro_limage_dun_type_symetrique_est_symetrique}, the definable types $F$ and $G$ commute with $r$. Then, applying this corollary again, we deduce, by associativity, that any tensor product whose factors are among $F$, $G$, or $r$, is commutative. In other words, the family $\lbrace F, G, r \rbrace$ is commutative.
\end{rem}

Recall that $h_{b_1 b_2}$denotes the composite $g_{b_2} \circ f_{b_1}$.

\begin{lemma}\label{h_b1b_2_algebrique_sur_b_3_cas_gen_stable}
The germ of the definable map $h_{b_1 b_2}$ is interalgebraic over $C$ with $b_3$.
\end{lemma}

\begin{proof}
We know that we have a partial group configuration over $C$. Thus, the element $b_3$ is interalgebraic over $C$ with the canonical basis of the generically stable type $tp(a_1 a_2 / acl(C b_3))$. Moreover, we have $a_1 a_2 \downfree_{acl(C b_3)} b_1 b_2$, because $a_1 \downfree_C b_1 b_2 b_3$ and $a_2 \in acl(C a_1 b_3)$. Thus, by stationarity of $tp(a_1 a_2 / acl(C b_3))$, we have the following equalities: $$Cb(a_1 a_2 / acl(C b_3)) = Cb(a_1 a_2 / acl(C b_1 b_2 b_3)) = Cb(a_1 a_2 / acl(C b_1 b_2)).$$ Since $a_1 \downfree_C b_1 b_2$, we have by stationarity $Cb(a_1 / acl(C b_1 b_2)) = Cb(a_1 / C) \subseteq C$. Then, Proposition \ref{prop_description_base_canonique_af_c(a)} implies that $Cb(a_1 a_2 / acl(C b_1 b_2))$ is interdefinable over $C$ with $[h_{b_1 b_2}]$. This concludes the proof.
\end{proof}

\begin{defi}\label{defi_Gamma}
Let $K$ be the set of germs of the form $f^{-1} \circ f'$ where $(f, f') \models F\otimes F |_{C}$.

Similarly, let $L$ be the set of germs of the form $f' \circ f^{-1}$ where $(f, f') \models F\otimes F |_{C}$.  Finally, let $\Gamma$  be the $C$-type-definable set of germs of the form $k \circ k'$ where $k,k' \in K$.
\end{defi}

\begin{rem}\label{rem_they_act_generically}
The set $K$ is then defined by a complete $C$-definable type, also denoted as $K$. Indeed, the type $K$ is the image of the tensor product $F \otimes F$ under the definable map which composes a germ with the inverse of another germ. This map is well-defined, by Corollary \ref{coro_composee_germs_bien_def_cas_gen_stable}. 
Note that, thanks to the strong germs property (Corollary \ref{germs_forts_cas_gen_stable}) and Corollary \ref{coro_composee_germs_bien_def_cas_gen_stable}, the realizations of $F$, $K$ and $\Gamma$ act generically, in the sense of Definition \ref{defi_acting_generically}, on the definable type $tp(a_2 / C)$, and those of $G$ and $L$ on the type $tp(a_3/C)$.
\end{rem}

The type-definable set $\Gamma$ is the underlying set of the group we are going to construct.

\begin{rem}\label{K=s(p_otimes_p)_cas_gen_stable}
For ease of notation, let us write $p=tp(b_1 / C)$, $ q=tp(b_2 / C)$.

\begin{enumerate}

\item Since $F \otimes F$ is a complete type, the image of $p \otimes p$ under the $C$-definable map  $(b, b^{'}) \mapsto ([f_{b}^{}], [f_{b^{'}}])$ is equal to the type $F \otimes F$. More generally, any finite tensor product whose factors are $F$ and $G$ is the image under the appropriate function  of the tensor product of corresponding factors $p$ and $q$. This follows from the definitions.

\item Therefore, the type $K$ is the image of the type $p \otimes p$ under the $C$-definable map  $(b, b^{'}) \mapsto [f_{b}^{-1} \circ f_{b^{'}}]$. Similarly, the type $L$ is the image of $p\otimes p$ by the function $(b, b^{'}) \mapsto [f_{b^{'}} \circ f_{b}^{-1}]$.

\item By Remark \ref{rem_F_et_G_et_p_commutent}, we know that $F\otimes F$ commutes with $F$, $G$ and $r$, where $r$ is the type $tp(a_1 a_2 a_3 b_1 b_2 b_3 / M)$. Then, by Corollary \ref{coro_limage_dun_type_symetrique_est_symetrique}, the type $K$ commutes with $F$, $G$ and $r$. Again, we deduce that $F$, $G$, $K$, and $r$ are in a commutative family, then that the family $\lbrace F, G, K, L, r \rbrace$ is commutative.

\item Since $F$ commutes with itself, the inverse of a germ $k$ realizing $K|_D$ is still a realization of $K|_D$, for all $D \supseteq C$. Similarly for $L$.
\end{enumerate} 

\end{rem}

The following lemma shows that the collection of germs is, in some sense, homogeneous. It will be used for several key results.

\begin{lemma}\label{gf=g2f2_cas_gen_stable}

Let $D$ be a small set containing $C$.

\begin{enumerate}
    \item Let $g,f_1, f_2$ be such that $g f_1 f_2$ realizes $G \otimes F \otimes F|_D$. Then, there exists $b^2_2$ realizing $tp(b_2 / C)|_D$ such that $g \circ f_1 = [g_{b^2_2}]\circ f_{2} $, $b^2_2 \downfree_C D g f_1$ and $b^2_2 \downfree_C D g f_2$.
    \item Let $g, g', f_1 $ be  such that $g g' f_1$ realizes $G \otimes G \otimes F|_D$. Then, there exists $b_1^2$ realizing $tp(b_1 / C)|_D$ such that $g\circ f_1 = g'\circ [f_{b_1^2}] $, $b_1^2 \downfree_C D g f_1$ and $b_1^2 \downfree_C D g' f_1$.
    
    %\item Let $b_1^1, b_1^2$ realizations of $tp(b_1 / C)$, and $ b_2^1$ a realization of $tp(b_2 / C)$, such that $b_2^1 \downfree_C b_1^1$ and $b_2^1 b_1^1 \downfree_C b_1^2$. Then, il existe $b^2_2$ realizing $tp(b_2 / C)$ such that $[g_{b_2^1}] \circ [f_{b_1^1}] = [g_{b^2_2}]\circ [f_{b_1^2}] $ and $b^2_2 \downfree_C b_2^1 b_1^1$.
    %\item Let $b_1^1$ a realization of $tp(b_1 / C)$, and $ b_2^1, b_2^2$ realizations of $tp(b_2 / C)$, such that $b_2^1 \downfree_C b_1^1$ and $b_2^1 b_1^1 \downfree_C b_2^2$. Then, il existe $b^2_1$ realizing $tp(b_1 / C)$ such that $[g_{b_2^1}] \circ [f_{b_1^1}] = [g_{b^2_2}]\circ [f_{b_1^2}] $ and $b^2_1 \downfree_C b_2^1 b_1^1$.

\end{enumerate}

\end{lemma}

%
%\begin{claim}
%It suffices to prove the result in the case $D=C$.
%
%\end{claim}
%
%\begin{proof}
%Our proof in the case $D=C$ will only use the fact that $(a_1, a_2, a_3, b_1, b_2, b_3)$ is a generically stable group configuration over $C$
%
%\end{proof}

\begin{proof} Let us prove the first result, and then explain how to prove the second one.

\begin{claim}
Assuming $b_2^2 \models tp(b_2 / C)|_D$, to prove that $b^2_2 \downfree_C D g f_1$ and $b^2_2 \downfree_C D g f_2$, it suffices to show $b^2_2 \downfree_D g f_1$ and $b^2_2 \downfree_D g f_2$.
\end{claim}

\begin{proof}
This is a direct application of transitivity (Proposition \ref{transitivite_gen_stable}).
\end{proof}

To simplify notations, let us assume that $tp(a_1 a_2 a_3 b_1 b_2 b_3 / D) = tp(a_1 a_2 a_3 b_1 b_2 b_3 / C)|_D$, i.e. $a_1 a_2 a_3 b_1 b_2 b_3 \downfree_C D$. Recall that $a_1 a_2 a_3 b_1 b_2 b_3$ is the configuration we built in Proposition \ref{prop_WMA_dcl}, and $D$ is not necessarily contained in $M$, so this assumption is not vacuous. We wish to be able to simply write, say $\alpha \equiv_D a_1$, instead of $\alpha \models tp(a_1/C)|_D$.

Let $q(x,y,z)$ denote the tensor product $tp(b_2 / D)(x) \otimes tp(b_1/D)(y) \otimes tp(b_1 / D)(z)$. By the first point of Remark \ref{K=s(p_otimes_p)_cas_gen_stable},  there exist $b^1_2, b_1^1, b_1^2$  such that $[g_{b^1_2}] = g$, $[f_{b_1^1}]=f_1$,  $[f_{b_1^2}]=f_2$ and $b^1_2 b_1^1 b_1^2 \models q(x,y,z) $. We look for a suitable element $b^2_2$. Let $\beta_3 \in acl(C b^1_2 b_1^1)$ be such that $\beta_3 b^1_2 b_1^1 \equiv_D b_3 b_2 b_1$.
Let $\alpha_2 $ be a realization of $tp(a_2/M)|_{ D b_1^1 b_1^2 b^1_2 \beta_3} $. Then, by stationarity, we have $\alpha_2 b_1^1 b^1_2 \beta_3 \equiv_D a_2 b_1 b_2 b_3$. Let $a_3^1 = f_{b_1^1}(\alpha_2)$ and $\alpha_1 = g_{b^1_2}(a_3^1)$. Thus, we have $\alpha_1 \alpha_2 a_3^1 b_1^1 b^1_2 \beta_3 \equiv_D a_1 a_2 a_3 b_1 b_2 b_3$.

By choice of $\alpha_2$, we have $\alpha_2 \downfree_D b_1^1 b^1_2 b_1^2$. By construction and commutativity, we also have $b_1^2 \downfree_D b_1^1 b^1_2$. Thus, we may apply Lemma \ref{lemma_switch}, to deduce that $\alpha_2 b_1^1 b^1_2 \downfree_D b_1^2$. So, by Fact \ref{fact_forking}(4)(b), we have

\begin{equation}\label{eq_1} \alpha_1 \alpha_2 a_3^1 b_1^1 b^1_2 \beta_3 \downfree_D b_1^2  \end{equation} 

Then, by symmetry and stationarity, $\alpha_1 \alpha_2 \beta_3  b_1^2 \equiv_D a_1 a_2 b_3 b_1$. Let $b^2_2, a_3^2$ be such that 
$b^2_2 a_3^2 \alpha_1 \alpha_2 \beta_3  b_1^2 \equiv_D b_2 a_3 a_1 a_2 b_3 b_1$. We will show that $b_2^2$ has the required properties.

We end up with the following generically stable group configurations, which have $(\beta_3, \alpha_2, \alpha_1)$-line in common, and whose type over $D$ is $r|_D = tp(a_1, a_2, a_3, b_1, b_2, b_3 / D)$.

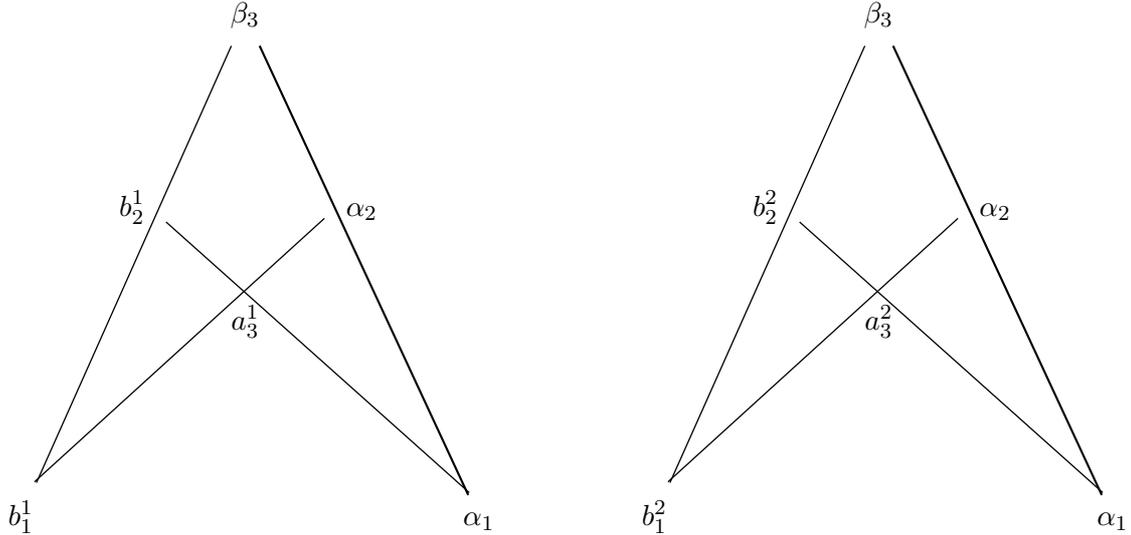
\begin{center}
    \begin{tikzcd}

&& \beta_3 \arrow[llddddd, dash, shorten= 1 mm] \arrow[rrddddd, dash, thick, shorten= 1 mm] &&    

&&

&& \beta_3 \arrow[llddddd, dash, shorten= 1 mm] \arrow[rrddddd, dash, thick, shorten= 1 mm] &&   

\\
 
\\
 & b^1_2 & & \alpha_2 &    
 
 &&
 
  & b^2_2 & & \alpha_2 &    

\\
&& a^1_3 &&  

&&

&& a^2_3 &&

\\

\\

b^1_1 \arrow[rrruuu, dash, shift left = 1.4 ex]&&&& \alpha_1\arrow[llluuu, dash, shift right = 1.0 ex]  

&&

b^2_1 \arrow[rrruuu, dash, shift left = 1.4 ex]&&&& \alpha_1\arrow[llluuu, dash, shift right = 1.0 ex]

    \end{tikzcd}
\end{center}

\normalsize

We know that $b_2^2 \downfree_D b_1^2$, because $b_2 \downfree_D b_1$. Let us show that $b_2^2 \downfree_D b^1_2 b^1_1$. By \ref{eq_1}, we know that $\alpha_1 \alpha_2 a_3^1 b_1^1 b^1_2 \beta_3 \downfree_D b_1^2$. Then, as $tp(\alpha_1 \alpha_2 a_3^1 b_1^1 b^1_2 \beta_3/D)$ is generically stable, we can apply symmetry, to get $b_1^2 \downfree_D \alpha_1 \alpha_2 a_3^1 b_1^1 b^1_2 \beta_3$. So, by Fact \ref{fact_forking}(1), we have  $b_2^2 \downfree_{D\beta_3} \alpha_1 \alpha_2 a_3^1 b_1^1 b^1_2$. Since $b_2^2 \downfree_D \beta_3$, by transitivity, we have $b_2^2 \downfree_{D}\beta_3 \alpha_1 \alpha_2 a_3^1 b_1^1 b^1_2$. By monotonicity, we deduce $b_2^2 \downfree_{D} b_1^1 b^1_2$. Hence, we have indeed $b_2^2 \downfree_D g f_1$.

In order to prove that $b_2^2 \downfree_D g f_2$, we will show that $b_2^2 \downfree_{D}b^1_2 b_1^2 $. By construction, we know that $b_1^2 \downfree_D b_1^1 b_2^1$, so $b_1^2 \downfree_D \beta_3 b_2^1 $. Since $\beta_3 \downfree_D b_2^1$, by Lemma \ref{lemma_switch}, we have $\beta_3 b_1^2 \downfree_D  b_2^1 $, so $b_2^2 b_1^2 \downfree_D  b_2^1 $. Applying Lemma \ref{lemma_switch} again, we have $b_2^2  \downfree_D  b_2^1 b_1^2$. Therefore, we have $b_2^2 \downfree_D g f_2$.

Besides, we know that $g_{b^1_2}\circ f_{b_1^1} (\alpha_2) = \alpha_1 = g_{b^2_2}\circ f_{b_1^2}(\alpha_2) $. It then remains to show that $\alpha_2 \downfree_D b^1_2 b_1^1 b^2_2 b_1^2$, so that we can conclude equality of the germs of  $g_{b^1_2}\circ f_{b_1^1}$ and $ g_{b^2_2}\circ f_{b_1^2} $. It suffices to prove that $\alpha_2 \downfree_D b_1^2 b^1_2 \beta_3$, which is true by choice of $\alpha_2$.

In order to prove the second statement, it suffices to swap the roles of $b^2_2$ and $b_1^2$. The nine points that are obtained have the same properties in both cases, only the order in which they are defined is different.
\end{proof}

\begin{lemma}\label{composee_germs_indep_cas_gen_stable}
Let $D$ be a small set containing $C$.
Let $f, f' \in F$ be such that $(f , f') \models F \otimes F |_{D}$. Then $f^{-1} \circ f' \models K|_{ D f}$ and $f^{-1} \circ f' \models K|_{ D f^{'}}$. On the other hand, $f' \circ f^{-1} \models L|_{ D f}$ and $f' \circ f^{-1} \models L|_{ D f^{'}}$.
\end{lemma}

\begin{proof}
Let us start by proving the statements about $f^{-1} \circ f'$. By the first point of Remark \ref{K=s(p_otimes_p)_cas_gen_stable}, we can apply Remark \ref{rem_h(p)} to the case where $q=tp(b_1/M)^{\otimes 2}$ and $h: ({b,b^{'}}) \mapsto [f_b^{-1} \circ f_{b^{'}}]$. We then find $\beta_1, \beta^{'}_1$ such that $([f_{\beta_1}], [f_{\beta^{'}_1}])=(f , f')$ and $(\beta_1, \beta^{'}_1) \models q|_{D}$.

Besides, the definition of $K$ implies that  $f^{-1} \circ f' \models K|_{ D}$, since we assumed that  $(f , f') \models F \otimes F |_{D}$. Recall that we want to show that $(f^{-1} \circ f', f) \models K\otimes F|_{ D}$ and $(f^{-1} \circ f', f') \models K\otimes F|_{ D}$. By the third point of Remark \ref{K=s(p_otimes_p)_cas_gen_stable}, $K$ and $F$ commute, so it is equivalent to prove that $(f, f^{-1} \circ f') \models F\otimes K|_{ D}$ and $(f', f^{-1} \circ f') \models F\otimes K|_{ D}$. By stationarity of $F$, as we have seen above that $f^{-1} \circ f' \models K|_{ D}$, it suffices to show that  $f \downfree_C D f^{-1} \circ f'$ and $f' \downfree_C D f^{-1} \circ f'$. By symmetry of the hypotheses on $f$ and $f'$, which comes from commutativity (see for instance Remark \ref{K=s(p_otimes_p)_cas_gen_stable} (3), or Fact \ref{fact_gen_stable_commute_with_lui_meme}), it is enough to prove that $f \downfree_C D f^{-1} \circ f'$.

Let $g = [g_{\beta_2}] \in G$, where $\beta_2 \models tp(b_2 / M) |_{D \beta_1 \beta^{'}_1}$. By stationarity, we know that $tp(b_1 b_2 /M)=tp(b_1 /M) \otimes tp(b_2 /M)$. Since $M$ is a model, we know by Fact \ref{fact_unique_invariant_extension} that $tp(b_1 b_2 /M)|_{E}=(tp(b_1 /M) \otimes tp(b_2 /M))|_{E}$ for all $E \supseteq C$. Then, $\beta_2 \beta_1$ and $\beta_2 \beta^{'}_1$ realize the type $tp(b_2 b_1/M)|_D$. Let $\beta_3, \beta^{'}_3$ be such that $\beta_3 \beta_2 \beta_1 $ and $ \beta^{'}_3 \beta_2 \beta^{'}_1$ realize $tp(b_3 b_2 b_1 / M)|_D$.
Then, by Lemma \ref{h_b1b_2_algebrique_sur_b_3_cas_gen_stable}, we have $[g_{\beta_2}]\circ [f_{\beta_1}] \in acl(C \beta_3)$ and $[g_{\beta_2}]\circ [f_{\beta^{'}_1}] \in acl(C \beta^{'}_3)$. Thus $f^{-1} \circ f' = ([g_{\beta_2}]\circ [f_{\beta_1}])^{-1} \circ [g_{\beta_2}]\circ [f_{\beta^{'}_1}] \in acl(C \beta_3 \beta^{'}_3)$. 

If we manage to prove that $\beta_1 \downfree_C D\beta_3\beta^{'}_3$, we can then deduce by Fact \ref{fact_forking} (4)(b) that $f \downfree_C D f^{-1} \circ f'$, which will finish the proof. By transitivity, it suffices to show that $\beta_1 \downfree_D \beta_3\beta^{'}_3$.

 By construction, we have  $\beta_2\beta_1 \beta^{'}_1 \models (tp(b_2 /M) \otimes q)|_D $, so $\beta_1 \downfree_D \beta_2 \beta^{'}_1$. Since $\beta^{'}_3 \in acl(D \beta_2 \beta^{'}_1)$, this implies by Fact \ref{fact_forking}(4)(b) that $\beta_1 \downfree_D \beta_2 \beta^{'}_3$. By choice of $\beta^{'}_3$, we have $\beta_2 \downfree_D \beta^{'}_3$. So, by Lemma \ref{lemma_switch}, we have 
$\beta_1 \beta_2 \downfree_D \beta^{'}_3$. As $\beta_3 \in acl(D \beta_1 \beta_2)$, this implies $\beta_1 \beta_3 \downfree_D \beta^{'}_3$. We also know that $\beta_1 \downfree_D \beta_3$. So, again by Lemma \ref{lemma_switch}, we can conclude that   $\beta_1 \downfree_D  \beta_3\beta^{'}_3$, as desired.

Now, to prove the result for $f^{'} \circ f^{-1}$, we use Lemma \ref{gf=g2f2_cas_gen_stable}. As above, by commutativity and stationarity, it suffices to show $f \downfree_C D f^{'} \circ f^{-1}$ and  $ f' \downfree_C Df^{'} \circ f^{-1}$. By transitivity, it suffices to show $f \downfree_D f^{'} \circ f^{-1}$ and  $ f' \downfree_Df^{'} \circ f^{-1}$. Let $g \models G|_{D f f^{'}}$. Then, by Lemma \ref{gf=g2f2_cas_gen_stable}, there exists $g' \models G|_{D}$ such that $g\circ f = g' \circ f'$, $g' \downfree_C D g f$ and $g' \downfree_C D g f'$. In particular, we have $g' \downfree_D g f$ and $g' \downfree_D g f'$. So, by Lemma \ref{lemma_switch} and symmetry, we have $f \downfree_D g \flex g' $ and $f'\downfree_D g \flex g'$. From the equality $g\circ f = g' \circ f'$, we deduce ${g'}^{-1} \circ g =  f^{'} \circ f^{-1}$. This implies $f \downfree_D f^{'} \circ f^{-1}$ and  $ f' \downfree_Df^{'} \circ f^{-1}$, as required.
\end{proof}

%
%Let us now sketch the proof for $f^{'} \circ f^{-1}$. It is similar in spirit, one needs to ``swap the roles of $b_2$ and $b_3$''. We shall only explain how to prove $f \downfree_C D f^{'} \circ f^{-1}$.  Let $\alpha_3 \models tp(a_3/C)|_{D \beta_1 \beta^{'}_1}$. Let $\alpha_2 = f^{-1}(\alpha_3) = f_{\beta_1}^{-1}(\alpha_3)$.
%
%\begin{claim}
%We have $\beta^{'}_1\downfree_C D \alpha_2$.
%\end{claim}
%
%\begin{proof}
%By commutativity, we can check that $\beta^{'}_1\downfree_C D \beta_1 \alpha_3$. The claim follows by Fact \ref{fact_forking}(4)(b).
%\end{proof}
%
%
%Then, let $\alpha^{'}_3 = f^{'}(\alpha_2) = f^{'} \circ f^{-1}(\alpha_3)$.

%Similar to the first point, let $\beta^{''}_3 \models tp(b_3 / C)|_{D \beta_1 \beta^{'}_1}$. As above, we know that $\beta_1 \beta^{''}_3$ and $\beta^{'}_1 \beta^{''}_3$ realize $tp(b_1 b_3 / C)|_{D}$, so, let $\beta^{'}_2$ and $\beta^{''}_2$ be such that $\beta_1 \beta^{''}_2 \beta^{''}_3 $and $\beta^{'}_1 \beta^{'}_2 \beta^{''}_3$ realize $tp(b_1 b_2 b_3 / C)|_D$.

%We have $f^{'} \circ f^{-1} = [g_{\beta_2}^{-1}] \circ ([g_{\beta_2}]\circ [f_{\beta^{'}_1}]) \circ ([g_{\beta_2}]\circ [f_{\beta_1}])^{-1} \circ [g_{\beta_2}]$. Thus, $f^{'} \circ f^{-1} \in acl(C \beta_2 \beta_3 \beta^{'}_3)$. Moreover, we have $\beta_2 \downfree_C D \beta_1 \beta^{'}_1$.

% thanks to Lemma \ref{h_b1b_2_algebrique_sur_b_3_cas_gen_stable}, we know that $[h_{\beta_1 \beta_2}]$ is interalgebraic over $C$ with $\beta_3$, and 

% Remark: this lemma does not use the point a_i in the configuration, except through the use of Lemma h_b1b_2_algebrique_sur_b_3_cas_gen_stable ! 

\begin{coro}\label{coro_description_of_Gamma}
\begin{enumerate}
\item The $C$-definable type $K$ is generically stable.

\item The type-definable set $\Gamma$ is the set of germs of the form $f^{-1} \circ f^{'}$, where $f, f^{'}$ realize $F|_C$.

\end{enumerate}

\end{coro}

\begin{proof}

1. Let $f$ realize $F|_C$. Then, by Lemma \ref{composee_germs_indep_cas_gen_stable}, there exists a $Cf$-definable bijection $F|_{Cf} \simeq K|_{Cf}$, which maps any $f'$ to $f^{-1} \circ f'$. Thus, by Proposition \ref{prop_gen_stabilite_et_acl}, the type $K|_{Cf}$ is generically stable, in the sense of \ref{defi_gen_stable_abus}. However, this type is definable over $C$, so $K|_C$ is generically stable, as stated.

2. Let $\gamma$ be an element of $\Gamma$. By definition, there exist $k_1, k_2$ realizing $K|_C$ such that $\gamma = k_1 \circ k_2$. Let $f \models F|_{C k_1 k_2}$. Then, by Lemma \ref{composee_germs_indep_cas_gen_stable}, and completeness of the type $K|_{Cf} = K^{-1}|_{Cf}$ (see Remark \ref{K=s(p_otimes_p)_cas_gen_stable}(4)), there exist $f_1, f_2 \in F$ such that $k_1 = f_1^{-1} \circ f$ and $k_2 = f^{-1} \circ f_2$. Then, we compute that $\gamma =f_1^{-1} \circ f_2$, as desired. \end{proof}

The following lemma will be used to prove transitivity of the action of the group $\Gamma$, and regularity in the case of a regular group configuration.

\begin{lemma}\label{lemma_action_reguliere_de_K_cas_gen_stable}
Let $\alpha_2, \alpha^{'}_2$ be realizations of $tp(a_2 / C)$ such that $\alpha_2 \downfree_C \alpha^{'}_2$. Then, there exists a germ $k \in K$ such that $k \models K|_{C \alpha_2} \cup K|_{C \alpha^{'}_2}$ and $k(\alpha_2) = \alpha^{'}_2$.

Moreover, under the hypothesis $(R)$ of Theorem \ref{theo_config_groupe_faible}, i.e. if we started from a generically stable \emph{regular} group configuration, there exist only finitely many germs $k$ in $K$ such that $k \models K|_{C \alpha_2}$ and $k(\alpha_2) = \alpha^{'}_2$. 
\end{lemma}

\begin{proof}
First, recall that if $k \in K$ realizes $K|_{C \alpha_2}$, then the couple $(k, \alpha_2)$ realizes $K \otimes tp(a_2 / M)|_{C}$. So, by commutativity (see the third point of Remark \ref{K=s(p_otimes_p)_cas_gen_stable}), $\alpha_2$ realizes $tp(a_2/M)|_{Ck}$. Then, since $K$ acts generically on $tp(a_2 / C)$ (see Remark \ref{rem_they_act_generically}), the element $k(\alpha_2)$ realizes $tp(a_2/M)|_{Ck}$.  

Let $p=tp(b_1 / M)$ and $t=tp(a_2 / M)$. By stationarity of $t$, the hypothesis is equivalent to $\alpha_2 \alpha^{'}_2$ realizing $t \otimes t|_C$.

Let us prove existence: Let $\alpha_3 \equiv_C a_3$ be such that $\alpha_3 \downfree_C \alpha_2 \alpha^{'}_2$, and let $\beta_1 \beta^{'}_1$ realize $p\otimes p |_{C \alpha_3}$.

\begin{claim}\label{claim_f_beta_1_-1_alpha_3}
The pair $(f_{\beta_1}^{-1}(\alpha_3), f_{\beta^{'}_1}^{-1}(\alpha_3))$ realizes the type $t\otimes t |_{C\alpha_3}$.
\end{claim}

\begin{proof}
By definition, we have $\beta_1 \models p|_{C \beta^{'}_1 \alpha_3}$. So $\beta_1 \downfree_{C\alpha_3} \beta^{'}_1$, so $f_{\beta_1}^{-1}(\alpha_3)\downfree_{C \alpha_3}f_{\beta^{'}_1}^{-1}(\alpha_3)$. Also, by stationarity, since $\beta_1 \downfree_C \alpha_3$ (and $b_1 \downfree_C a_3$), we have $\beta_1 \alpha_3 \equiv_C b_1 a_3$. Moreover, in the configuration $a_1a_2a_3b_1b_2b_3$, we have $f_{b_1}^{-1}(a_3)=a_2$, and $a_2 \downfree_C a_3$. So $f_{\beta_1}^{-1}(\alpha_3) \downfree_C \alpha_3$, and similarly $f_{\beta^{'}_1}^{-1}(\alpha_3) \downfree_C \alpha_3$. So, by transitivity, $f_{\beta_1}^{-1}(\alpha_3)\downfree_{C}f_{\beta^{'}_1}^{-1}(\alpha_3) \alpha_3$. In other words, $f_{\beta_1}^{-1}(\alpha_3)f_{\beta^{'}_1}^{-1}(\alpha_3)$ realizes the type $t\otimes t |_{C\alpha_3}$, as required.
\end{proof}

 By choice of $\alpha_3$, the tuple $\alpha_2 \alpha^{'}_2$ also realizes the type $t\otimes t |_{C\alpha_3}$. So, up to changing $\beta_1$ and $\beta^{'}_1$, we may assume that $f_{\beta^{}_1}^{-1}(\alpha_3)=\alpha_2$ and $f_{\beta^{'}_1}^{-1}(\alpha_3)=\alpha^{'}_2$, without changing the fact that $\beta_1 \beta^{'}_1$ realizes $p\otimes p |_{C \alpha_3}$, and without changing $\alpha_3$.

Then $(f_{\beta^{'}_1}^{-1}\circ f_{\beta_1})(\alpha_2) = \alpha^{'}_2$, considering the definable maps, and not their germs. We want to show that $k=[f_{\beta^{'}_1}^{-1}\circ f_{\beta_1}]$ has the required properties.

\begin{claim}
We have $\beta_1 \beta^{'}_1 \models p^{\otimes 2}|_{C\alpha_2}\cup p^{\otimes 2}|_{C\alpha^{'}_2} $.

\end{claim}
\begin{proof}
By construction of $\beta_1, \beta^{'}_1$, we know that $\beta_1 \models p|_{C \alpha_3}$ and that $\beta^{'}_1 \models p|_{C \alpha^{}_3}$. Then, by the third point of Proposition \ref{prop_f_b(a)_independant_de_Ab_cas_gen_stable}, we deduce that  $\beta_1 \models p|_{C \alpha_2}$ and that $\beta^{'}_1 \models p|_{C \alpha^{'}_2}$. Since $p$ commutes with itself, it remains to show that $\beta^{'}_1 \models p|_{C \alpha^{}_2\beta_1}$ and $\beta_1 \models p|_{C \alpha^{'}_2\beta^{'}_1}$. By symmetry of the construction, it suffices to prove the second point. Recall that, by definition of $\beta_1 \beta^{'}_1$ (right before Claim \ref{claim_f_beta_1_-1_alpha_3}), we have $\beta_1 \models p|_{C \beta^{'}_1 \alpha_3}$. Since $\alpha^{'}_2 \in dcl(C  \beta^{'}_1\alpha_3)$, we have indeed that $\beta_1 \models p|_{C \alpha^{'}_2\beta^{'}_1}$.
\end{proof}

We know by Corollary \ref{coro_limage_dun_type_symetrique_est_symetrique} that $p(x) \otimes p(y) \otimes t(z) = t(z) \otimes p(x) \otimes p(y)$. So, the claim implies that $\alpha_2 \models t|_{C \beta_1 \beta^{'}_1}$. Thus, $k$ has the required properties.

\medskip

Now, let us assume that the hypothesis $(R)$ of Theorem \ref{theo_config_groupe_faible} holds, and prove finiteness. 
Let $k$ in $K$ be such that $k \models K|_{C \alpha_2} $ and $k(\alpha_2) = \alpha^{'}_2$. Let us show that $k \in acl(C \alpha_2 \alpha^{'}_2)$. As $K$ is the image of $p\otimes p$ (see the second point of Remark \ref{K=s(p_otimes_p)_cas_gen_stable}), we can apply Remark \ref{rem_h(p)}. We then find $\beta_1, \beta^{'}_1$ realizations of $p=tp(b_1 / C)$ such that  $\beta_1 \beta^{'}_1 \models p\otimes p|_{C \alpha_2}$ and $k = [f_{\beta^{'}_1}^{-1}\circ f_{\beta_1}]$. Then, by commutativity, $\alpha_2$ realizes $t|_{C\beta_1 \beta^{'}_1}$. Thus, we have $$f_{\beta^{'}_1}^{-1}\circ f_{\beta_1}(\alpha_2)= k(\alpha_2)=\alpha^{'}_2.$$

In order to symmetrize the information on $\alpha_2$ and $\alpha^{'}_2$, let us prove the following \begin{claim} We have $\beta_1 \beta^{'}_1 \models p\otimes p|_{C \alpha^{'}_2}$.\end{claim}

\begin{proof} Let $\alpha_3 = f_{\beta_1}(\alpha_2) = f_{\beta^{'}_1}(\alpha^{'}_2)$. Then, as $\beta_1 \beta^{'}_1  \models p\otimes p|_{C \alpha_2}$, we have, by commutativity, $ \beta^{'}_1 \models p|_{C \beta_1 \alpha_2}$, so $\beta^{'}_1 \models p|_{C \beta_1 \alpha_2\alpha_3} $. Moreover, $\beta_1 \models p|_{\alpha_2}$, so, by Proposition \ref{prop_f_b(a)_independant_de_Ab_cas_gen_stable} (1), commutativity and stationarity, we have $\beta_1 \models p|_{\alpha_3}$. So $\beta^{'}_1 \beta_1 \models p\otimes p|_{C \alpha_3}$. So, by commutativity,  $ \beta_1 \beta^{'}_1\models p\otimes p|_{C \alpha_3}$. So $\beta_1 \models p|_{C \alpha_3  \beta^{'}_1}$, so $\beta_1 \models p|_{C  \beta^{'}_1 \alpha^{'}_2}$. Also, applying Proposition \ref{prop_f_b(a)_independant_de_Ab_cas_gen_stable}(3) (and commutativity and stationarity) to the hypothesis ``$\beta^{'}_1 \models p|_{C \alpha_3}$'' , we have $\beta^{'}_1 \models p|_{C \alpha^{'}_2}$. Then, by definition of a tensor product, we have $\beta_1 \beta^{'}_1 \models p\otimes p|_{C \alpha^{'}_2}$, as stated. \end{proof}

Let $\beta_2$ realize $tp(b_2 / M)|_{C \alpha_3 \beta_1 \beta^{'}_1 \alpha_2 \alpha^{'}_2}$. By stationarity, we have $\beta_1 \beta_2 \equiv_C \beta^{'}_1 \beta_2 \equiv_C b_1 b_2$. Then, let $\beta_3 \beta^{'}_3$ be a couple such that $\beta_1 \beta_2 \beta_3 \equiv_C \beta^{'}_1 \beta_2 \beta^{'}_3 \equiv_C b_1 b_2 b_3$. Then, since $\beta_2 \downfree_C \beta_1 \alpha_2$, and $\beta_1 \downfree_C \alpha_2$, we have, by Lemma \ref{lemma_switch}, $\beta_2 \beta_1 \downfree_C  \alpha_2$, so $\alpha_2 \downfree_C  \beta_2 \beta_1$, so $\alpha_2 \downfree_C  \beta_3\beta_2 \beta_1$. Then, by stationarity, we have $\alpha_2 \beta_1 \beta_2 \beta_3 \equiv_C a_2 b_1 b_2 b_3$. So, we have $$\alpha_2 \alpha_3 \beta_1 \beta_2 \beta_3 \equiv_C a_2 a_3 b_1 b_2 b_3,$$ because $\alpha_3 =  f_{\beta_1}(\alpha_2) $. Let $\alpha_1:= g_{\beta_2}(\alpha_3)$. Then $\alpha_1 \alpha_2 \alpha_3 \beta_1 \beta_2 \beta_3 \equiv_C a_1 a_2 a_3 b_1 b_2 b_3$. By symmetric arguments, we also have $\alpha_1 \alpha^{'}_2 \alpha_3 \beta^{'}_1 \beta_2 \beta^{'}_3 \equiv_C a_1 a_2 a_3 b_1 b_2 b_3$.

Thus, we get the following configurations, which have the $(\beta_2, \alpha_3, \alpha_1)$-line in common: 

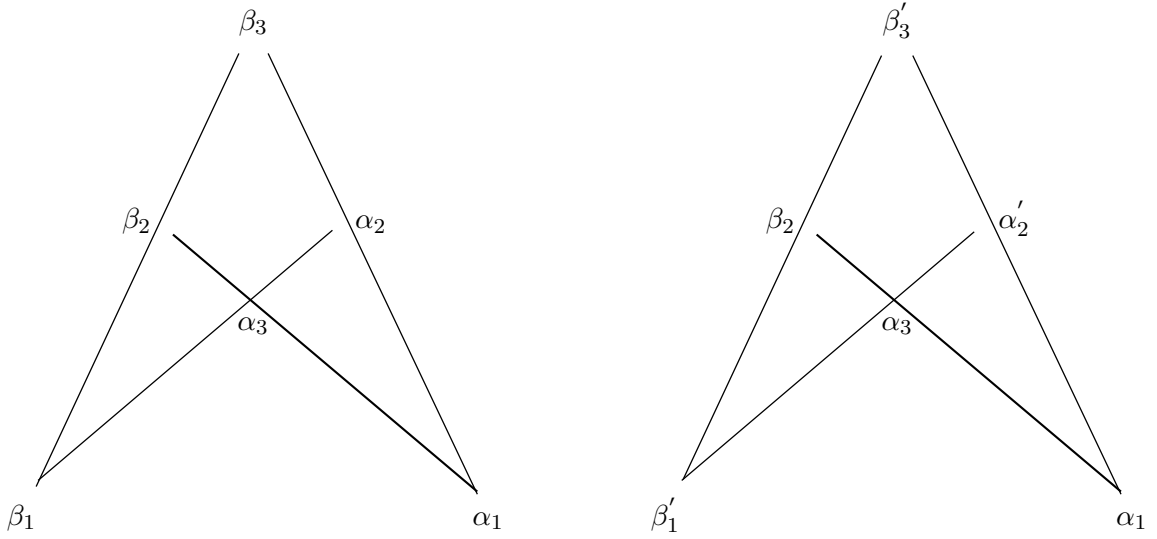
\begin{center}
    \begin{tikzcd}

&& \beta_3 \arrow[llddddd, dash, shorten= 1 mm] \arrow[rrddddd, dash, shorten= 1 mm] &&    

&&

&& \beta^{'}_3 \arrow[llddddd, dash, shorten= 1 mm] \arrow[rrddddd, dash, shorten= 1 mm] &&   

\\
 
\\
 & \beta_2 & & \alpha_2 &    
 
 &&
 
  & \beta_2 & & \alpha^{'}_2 &    

\\
&&\alpha_3 &&  

&&

&& \alpha_3 &&

\\

\\

\beta_1 \arrow[rrruuu, dash, shift left = 1.4 ex]&&&& \alpha_1\arrow[llluuu, dash,thick, shift right = 1.0 ex]  

&&

\beta^{'}_1 \arrow[rrruuu, dash, shift left = 1.4 ex]&&&& \alpha_1\arrow[llluuu, dash,thick, shift right = 1.0 ex]

    \end{tikzcd}
\end{center}

By Lemma \ref{h_b1b_2_algebrique_sur_b_3_cas_gen_stable}, the germ  $k=[f_{\beta^{'}_1}^{-1}\circ f_{\beta_1}] = [f_{\beta^{'}_1}^{-1} \circ g_{\beta_2}^{-1} \circ g_{\beta_2}\circ f_{\beta_1}] $ is algebraic over $C \beta_3 \beta^{'}_3$. Besides, using the hypothesis $(R)$ of Theorem \ref{theo_config_groupe_faible}, we know that $\beta_1 \beta^{'}_1 \in acl(C\alpha_3 \alpha_2 \alpha^{'}_2)$, so $k \in acl(C\alpha_3 \alpha_2 \alpha^{'}_2)$. If we show that $\alpha_3 \downfree_{C \alpha_2 \alpha^{'}_2}  \beta_3 \beta^{'}_3$, we can then apply Fact \ref{fact_forking} (5), to deduce that $k \in acl(C \alpha_2 \alpha^{'}_2)$. To that end, using Fact \ref{fact_forking} (1) and (4), and recalling that $\alpha^{'}_2 \in acl(C \beta^{'}_3 \alpha_1) \subseteq acl(C \beta^{'}_3 \beta_3 \alpha_2)$, it suffices to prove the
\begin{claim} We have $\alpha_3 \downfree_C \alpha_2 \beta_3 \beta^{'}_3$.\end{claim}

\begin{proof}
We have seen above that $ \beta^{'}_1 \models p|_{C \beta_1 \alpha_2}$. So $\beta^{'}_1 \downfree_C \beta_1 \alpha_2$. Moreover, by choice of $\beta_2$, we have $\beta_2 \downfree_C \beta^{'}_1\beta^{}_1\alpha_2$. The type $tp(\beta^{}_1\alpha_2/C)$ being generically stable, because $\beta^{}_1\alpha_2 \equiv_C b_1 a_2$, we may apply Lemma \ref{lemma_switch}, which yields $\beta_2\beta^{}_1\alpha_2 \downfree_C \beta^{'}_1$. By symmetry, we deduce that $ \beta^{'}_1 \downfree_C \beta_2\beta^{}_1\alpha_2$, so $ \beta^{'}_1 \downfree_{C \beta_2} \beta^{}_1\beta_2\beta_3 \alpha_1\alpha_2\alpha_3$. Since $\beta^{'}_3 \in acl(C \beta_2 \beta^{'}_1)$, this implies by Fact \ref{fact_forking}(4)(b) that $$ \beta^{'}_3 \downfree_{C \beta_2} \beta^{}_1\beta_2\beta_3 \alpha_1\alpha_2\alpha_3.$$ As $\beta^{'}_3 \downfree_C \beta_2$, we have, by transitivity for generically stable types, $ \beta^{'}_3 \downfree_{C} \beta^{}_1\beta_2\beta_3 \alpha_1\alpha_2\alpha_3$. So $ \beta^{'}_3 \downfree_{C} \beta_3 \alpha_2\alpha_3$. Since $\alpha_3 \downfree_C \beta_3 \alpha_2$, and $tp(\beta_3 \alpha_2/C)$ is generically stable, we can apply Lemma \ref{lemma_switch} again, to get $ \beta^{'}_3  \beta_3 \alpha_2 \downfree_{C}\alpha_3$. Now, the type $tp( \beta^{'}_3  \beta_3 \alpha_2/C)$ is extensible, for it is the tensor product $tp( \beta^{'}_3 /C) \otimes tp(\beta_3 \alpha_2/C)$. Thus, we may apply Proposition \ref{prop_symetrie}, to deduce that $\alpha_3 \downfree_C \alpha_2 \beta_3 \beta^{'}_3$, as desired.
\end{proof}

Thus, we have proved that $k \in acl(C \alpha_2 \alpha^{'}_2)$. This holds for all realizations of the partial type defined by ``$k \models  K|_{C \alpha_2}$ and $k(\alpha_2)=\alpha^{'}_2 $''  so, by compactness, there are only finitely many $k \in K$ satisfying $ k \models  K|_{C \alpha_2}$ and $k(\alpha_2)=\alpha^{'}_2$.
\end{proof}

%On vérifie then qu'elle reste vraie for the tensor product $F \otimes F$, puis for $K$, which is the image of $F \otimes F$ by a map $C$-definable.

We can now prove that $K$ behaves like the generic of a group: 

\begin{coro}\label{coro_produit_generique_dans_K_cas_gen_stable}
Let $(f_1, f_2, f_3, f_4)$ be a family of elements of $F$ which realizes the tensor product $F^{\otimes 4}|_C$. Let $E \supseteq C$. Assume that $f_1 f_2$ realizes $F \otimes F|_{E f_3 f_4}$. Then, there exist $f, f' \in F$ such that \begin{enumerate}     
   \item We have $(f_1^{-1}\circ f_2)\circ  (f_3^{-1} \circ f_4)=f^{-1}\circ f' $
    \item The pair $(f,  f')$ realizes $F \otimes F |_{C f_1^{-1}\circ f_2}$ and $F \otimes F|_{E f_3  f_4}$.
    \item We have $(f_1^{-1}\circ f_2)\circ  (f_3^{-1} \circ f_4) \models K|_{C f_1^{-1}\circ f_2}$ and $(f_1^{-1}\circ f_2)\circ  (f_3^{-1} \circ f_4) \models K|_{E f_3  f_4}$.

\end{enumerate}
\end{coro}

Note that neither the hypotheses nor the conclusion are symmetric: $f_3, f_4$ do not necessarily realize $F|_E$, whereas  $f_1$, $f_2$, $f$, and $f^{'}$ do.

\begin{proof} Let $g_2$ realize $G|_{E f_1 f_2 f_3 f_4}$. Then, by Lemma \ref{gf=g2f2_cas_gen_stable}, there exists $g_3 \in G$ such that $g_2 \circ f_3 = g_3 \circ f_2$ and $g_3 \downfree_C g_2 f_3$. We know that $g_2 \models G|_{C f_1 f_2 f_3 f_4}$, and $f_2 f_3 f_4\models F^{\otimes 3}|_{C f_1 }$. By choice of $g_2$, this implies $g_2 f_2 f_3 f_4 \models G \otimes F^{\otimes 3}|_{C f_1}$. Then, by commutativity, we have $f_2 f_3 \models F\otimes F |_{C g_2 f_1 f_4}$.

Then, by Lemma \ref{composee_germs_indep_cas_gen_stable} applied to $D= C g_2 f_1 f_4$, we have $f_3\circ f_2^{-1} \models L|_{C g_2 f_1 f_3 f_4}$. So $g_2\circ f_3 \circ f_2^{-1} \downfree_{C g_2} f_1 f_3 f_4$. In other words, $g_3 \downfree_{C g_2}f_1 f_3 f_4$. Since $g_3 \downfree_C g_2$, we then have $g_3 \downfree_C g_2 f_1 f_3 f_4$. Thus, by stationarity, $g_3 \models G|_{C g_2 f_1 f_3 f_4}$. Moreover, we have chosen $g_2$ so that $g_2  \models G|_{E f_1 f_2 f_3 f_4}$, hence $$g_3 g_2 f_4 \models G\otimes G \otimes F|_{C f_1 f_3}.$$ Then, we can again apply Lemma \ref{gf=g2f2_cas_gen_stable}, for the germs $g_2, g_3$ and $f_4$. We thus obtain a germ $f_5 \in F$ such that $g_3 \circ f_5 = g_2 \circ f_4$ and $f_5 \downfree_C g_2 f_4$. We will show that $f = f_1$ and $f'=f_5$ have the required properties.

Compute: $f_1^{-1} \circ f_2 \circ f_3^{-1} \circ f_4 = f_1^{-1} \circ g_3^{-1} \circ g_2 \circ f_4 = f_1^{-1} \circ f_5 $, these being equalities of germs.
In other words, $f_5 = f_2 \circ f_3^{-1} \circ f_4$. In particular, $f_5 \downfree_C f_1$, so, by stationarity, $f_1 f_5 \models F \otimes F|_C$.
%
%Then, let us prove that $f_1^{-1} \circ f_5 \models K|_{C f_1^{-1} \circ f_2}$. By the fourth point of Remark  \ref{K=s(p_otimes_p)_cas_gen_stable}, it is the same as showing that  $f_5^{-1} \circ f_1 \models K|_{C f_1^{-1} \circ f_2}$. By definition of $K$, it then suffices to show that $f_5 f_1 \models F\otimes F|_{C f_1^{-1} \circ f_2}$.

\begin{claim}
The pair $(f_1, f_5)$ realizes $F \otimes F |_{C f_1^{-1}\circ f_2}$.
\end{claim}

\begin{proof}
By Lemma \ref{composee_germs_indep_cas_gen_stable}, we know that $ f_1^{-1} \circ f_2 \models K|_{E f_1}$, so $( f_1^{-1} \circ f_2, f_1) \models K\otimes F|_E$. So, by commutativity, $f_1$ realizes $F|_{E  f_1^{-1} \circ f_2}$, so a fortiori $f_1$ realizes $F|_{C  f_1^{-1} \circ f_2}$. It then remains to show that $f_5 \models F|_{C f_1 \cup (f_1^{-1} \circ f_2)}$, i.e. $f_5 \models F|_{C f_1 f_2}$. By stationarity, it is enough to prove that $f_5 \downfree_C f_1 f_2$. Using the hypotheses on $(f_1, f_2, f_3, f_4)$, we have, by commutativity, $f_2 f_3 \models F \otimes F|_{C  f_1 f_4}$. So, by Lemma \ref{composee_germs_indep_cas_gen_stable} applied to $(f_3, f_2)$, with $D = C f_1 f_4$, we have $f_2 \circ f_3^{-1} \models L|_{C f_1 f_2 f_4}$. So $ f_5 = f_2 \circ f_3^{-1} \circ f_4  \downfree_{C f_4 } f_1 f_2$. By construction, we have $f_5  \downfree_{C } g_2 f_4$. Thus, by transitivity, $f_5 \downfree_{C } f_1 f_2 f_4$. In particular, $f_5 \downfree_{C} f_1 f_2$, as desired. 
\end{proof}

\begin{claim}
The pair $(f_1, f_5)$ realizes $F \otimes F |_{E f_3  f_4}$.
\end{claim}

\begin{proof}By commutativity, it suffices to show that $f_5 f_1$ realizes $F\otimes F|_{E f_3 f_4}$. We know by hypothesis that $f_1$ realizes $F|_{E f_3 f_4}$. By stationarity of $tp(f_5 / C)$, it remains to show that $f_5 \downfree_C E f_1 f_3 f_4$. On the one hand, by hypothesis (and symmetry), we have $f_2 \downfree_C E f_1 f_3 f_4$. So $f_2\circ f_3^{-1} \circ f_4 \downfree_{C f_3 f_4} E f_1$, i.e. 

\begin{equation}\label{eq_2}
f_5 \downfree_{C f_3 f_4} E f_1
\end{equation}

On the other hand, Lemma \ref{composee_germs_indep_cas_gen_stable} applied with $D=Cf_1 f_4$ implies that $f_2\circ f_3^{-1}$ realizes $L|_{C f_1 f_3 f_4}$, so $f_2\circ f_3^{-1} \downfree_C f_1 f_3 f_4$, so $f_5 \downfree_{C f_4} f_1 f_3$. We also know that $f_5 \downfree_C g_2 f_4$. So, by transitivity, $f_5 \downfree_C f_1 f_3 f_4$. So, by transitivity in \ref{eq_2}, we have $f_5 \downfree_C E f_1 f_3 f_4$, as desired. 
\end{proof}

Finally, that the third point follows from the first two points and the definition of $K$.
\end{proof}

\begin{coro}\label{coro_produit_generique_dans_K_simplifie} Let $k_1, k_2$ be realizations of $K|_C$, and $D \supseteq C$, such that $k_1 \models K|_{D k_2}$. Then $k_1 \circ k_2 \models K|_{D k_2}$ and $k_2 \circ k_1 \models K|_{D k_2}$.
\end{coro}

\begin{proof}
Let us show that $k_2 \circ k_1 \models K|_{D k_2}$, the other result being more straightforward. By definition of $K$, there are $f_3, f_4$ such that $(f_3, f_4) \models F \otimes F|_C$ and $ f_3^{-1} \circ f_4 = k_2$. Let $f_1, f_2$ in $F$ be such that $(f_1, f_2) \models F\otimes F|_{D f_3 f_4}$. Then, if $k'= f_1^{-1} \circ f_2$, we know that $k' \models K|_{D f_3 f_4}$, so in particular $k' \models K|_{D k_2}$, so $k' \equiv_{D k_2} k_1$. So, it suffices to prove that $k_2 \circ k' \models K|_{D k_2}$. 

Since $F$ commutes with itself, we can apply Corollary \ref{coro_produit_generique_dans_K_cas_gen_stable} to the family $(f_2,f_1 , f_4, f_3)$. It yields that $f_2^{-1} \circ f_1 \circ f_4^{-1} \circ f_3$ realizes $K|_{D k_2}$. Then, by Remark \ref{K=s(p_otimes_p)_cas_gen_stable}(4), the inverse $f_3^{-1} \circ f_4 \circ f_1^{-1} \circ f_2$ still realizes $K|_{D k_2}$. In other words, $k_2 \circ k' \models K|_{D k_2}$, as stated.

To show that $k_1 \circ k_2 \models K|_{D k_2}$, we also apply Corollary \ref{coro_produit_generique_dans_K_cas_gen_stable}, without permuting functions, nor considering inverses.
\end{proof}

Recall that, by Definition \ref{defi_Gamma}, the set $\Gamma$ is the set of composites $k_1 \circ k_2$, where $k_1, k_2$ realize $K|_C$. 

\begin{prop}
The $C$-type-definable-set $\Gamma$ is closed under composition of germs.
\end{prop}

\begin{proof}

Let $k_1, k_2, k_3, k_4$ realize $K|_C$. Let $D = C k_1k_2 k_3 k_4$. Let $k \models K|_D$. Then, using Corollary \ref{coro_produit_generique_dans_K_simplifie} four times, we can easily show that $k_1 \circ k_2 \circ k_3 \circ k_4 \circ k \models K|_D$.

Finally, we notice that $k_1 \circ k_2 \circ k_3 \circ k_4 = (k_1 \circ k_2 \circ k_3 \circ k_4\circ k) \circ (k^{-1})$. Since $k^{-1}$ and $(k_1 \circ k_2 \circ k_3 \circ k_4\circ k)$ are in $K$, the germ $k_1 \circ k_2 \circ k_3 \circ k_4$ is indeed in $K \circ K = \Gamma$.
\end{proof}

\begin{coro}
Composition of germs induces a  definable group structure on the type-definable set $\Gamma$.
\end{coro}

\begin{proof}
Corollary \ref{coro_composee_germs_bien_def_cas_gen_stable} implies that the composition of germs is associative, for it is induced by composition of functions. More precisely, let $\gamma_1, \gamma_2, \gamma_3 \in \Gamma$. By Corollary \ref{coro_description_of_Gamma}, let $\beta_1, \beta^{'}_1, \beta_2, \beta^{'}_2, \beta_3, \beta^{'}_3$ be realizations of $tp(b_1 / C)$ such that, for $i = 1, 2, 3$, we have $\gamma_i = [f_{\beta_i}^{-1} \circ f_{\beta^{'}_i}]$. Then, by definition, we have $(\gamma_1 \cdot \gamma_2) \cdot \gamma_3 = ([f_{\beta_1}^{-1} \circ f_{\beta^{'}_1}] \cdot [f_{\beta_2}^{-1} \circ f_{\beta^{'}_2}]) \cdot [f_{\beta_3}^{-1} \circ f_{\beta^{'}_3}]$. This is equal to the germ of the definable map $( f_{\beta_1}^{-1} \circ f_{\beta^{'}_1} \circ f_{\beta_2}^{-1} \circ f_{\beta^{'}_2}) \circ f_{\beta_3}^{-1} \circ f_{\beta^{'}_3}$. By associativity of composition for functions, that definable map is equal to  $f_{\beta_1}^{-1} \circ f_{\beta^{'}_1} \circ (f_{\beta_2}^{-1} \circ f_{\beta^{'}_2} \circ f_{\beta_3}^{-1} \circ f_{\beta^{'}_3})$.
Thus, by computing the germs of these maps, we get $(\gamma_1 \cdot \gamma_2) \cdot \gamma_3  = \gamma_1 \cdot (\gamma_2 \cdot \gamma_3)$, as desired.

Moreover, by the fourth point of Remark \ref{K=s(p_otimes_p)_cas_gen_stable}, $K$ is closed under taking inverses, and so is $\Gamma$. Besides, we have proved that $\Gamma$ is closed under composition.

Finally, the germ of the identity is indeed in $\Gamma$, for, if $k \in K$, then  $id = k \circ k^{-1} \in \Gamma$.
\end{proof}

\end{subsection}

\begin{subsection}{Properties of the group}

\begin{prop}
The type-definable group $\Gamma$ is connected, with generic $K$.
\end{prop}

\begin{proof}
First, recall that $K$ is a $C$-definable type. We will prove the following:  \begin{equation}\label{Stab(K)_cas_gen_stable}   Stab_{\Gamma}(K)=\Gamma
\end{equation}

Let $g \in \Gamma(M)$. Let $k$ be a realization of $K|_M$. By definition of $\Gamma$, and since $M$ is sufficiently saturated, there exist $k_1, k_2 \in K(M)$ such that $g = k_1 \circ k_2$. Then, applying Corollary \ref{coro_produit_generique_dans_K_simplifie} twice, we show that $g \circ k = k_1 \circ k_2 \circ k$ still realizes $K|_M$. Thus (\ref{Stab(K)_cas_gen_stable}) holds. We then apply Lemma \ref{lemma_gen_stable_generics_in_a_group}, to conclude that $K$ is the unique generic of $\Gamma$.\end{proof}

%Then, by definition of genericity, $K|_M$ is bien a generic type in $\Gamma$. Il reste à prove that $\Gamma$ is connected. Pour simplifier les notations, on écrira here $K$ instead of $K|_M$. Let $H$ be a subgroup of $\Gamma$, of bounded index, which is relativement type-$C$-definable. By finitude and saturation, toutes les classes in $\Gamma$ modulo $H$ are représentées in $M$. So $K$ is contained in a class modulo $H$. So $Stab_\Gamma(K) \subseteq H$. By (\ref{Stab(K)_cas_gen_stable}), we deduce that $H = \Gamma$. So $\Gamma$ is connected.

\begin{defi}

Let $Y$ be the set of pairs $(\gamma, \alpha)$ where $\gamma \in \Gamma$ and $\alpha \models tp(a_2 / C)$. Let $E$ be the equivalence relation on $Y$ defined by $(\gamma, \alpha) \, E \, (\gamma^{'}, \alpha^{'})$ if and only if there exists $\sigma \models K|_{C \gamma \gamma^{'} \alpha \alpha^{'}}$ such that $(\sigma\cdot \gamma)(\alpha) = (\sigma\cdot \gamma^{'})(\alpha^{'})$.

\end{defi}

\begin{prop}
The set $Y$ is type-definable over $C$, and $E$ is relatively $C$-definable. Let $X$ be the type-definable set $Y / E$.
\end{prop}
\begin{proof}
Type-definability over $C$ of $Y$ is immediate. We claim that, since $K$ is a complete definable type, the formula $\phi(\gamma_1, \alpha_1, \gamma_2, \alpha_2)$ defines $E$ inside $Y \times Y$, where $\phi(\gamma_1, \alpha_1, \gamma_2, \alpha_2) = d_K z \,  [(z \cdot \gamma_1)(\alpha_1) = (z\cdot \gamma_2)(\alpha_2)]$.
\end{proof}

\begin{lemma}\label{lemma_k_a_E_1_ka}

Let $(k, a) \models K|_C \otimes tp(a_2 / C)$. Then we have $(k, a) E (1, k(a))$.

\end{lemma}

\begin{proof}
Let $\sigma \models K|_{C a k}$. Then, by definition of the product in the group $\Gamma$, we have $\sigma(k(a)) = (\sigma \cdot k)(a)$, which proves the result.
\end{proof}

\begin{prop}

For each $\sigma \in \Gamma$, the map $(\gamma, a) \mapsto (\sigma\cdot \gamma, a)$ factorizes through the equivalence relation $E$, and this induces a definable action of $\Gamma$ on $X$. 

\end{prop}

\begin{proof}

Let $\sigma \in \Gamma$. Pick $(\gamma, a)$ and $(\gamma^{'}, a^{'})$ that are in the same $E$-class. Let us show that $(\sigma\cdot \gamma, a)$ and $(\sigma \cdot \gamma^{'}, a^{'})$ are in the same $E$-class. By assumption, there exists $\tau \models K|_{C  \gamma \gamma^{'}  a a^{'}}$ such that $(\tau \cdot \gamma) \cdot a = (\tau \cdot \gamma^{'}) \cdot a^{'}$. In fact, by completeness of the type $K|_{C  \gamma \gamma^{'}  a a^{'}}$, the equality holds for all such $\tau$. Let $\tau_1$ realize $K|_{C \sigma \gamma \gamma^{'}  a a^{'}}$. Then, by genericity of $K$ and Lemma \ref{lemma_gen_stable_generics_in_a_group}, the element $\tau_2=\tau_1 \cdot \sigma$ also realizes $K|_{C \sigma \gamma \gamma^{'}  a a^{'}}$. Thus, we have $(\tau_2 \cdot \gamma) \cdot a = (\tau_2 \cdot \gamma^{'}) \cdot a^{'}$, which implies  $(\tau_1 \cdot \sigma \cdot \gamma) \cdot a = (\tau_1 \cdot \sigma \cdot \gamma^{'}) \cdot a^{'}$, which proves that the map does factor through $E$.

The fact that this induces a definable action follows from the universal property of the quotient map $\pi : Y \rightarrow Y/E$, and the fact that $\Gamma$ acts on itself by left translation. More explicitly, let $\sigma, \tau \in \Gamma$, and $c = (\gamma, a) \in Y$. By construction, we have $\tau\cdot (\sigma \cdot \pi(c)) = \tau \cdot (\pi(\sigma \cdot \gamma, a)) = \pi(\tau \cdot ( \sigma \cdot \gamma), a) = \pi((\tau\cdot \sigma) \cdot \gamma, a) = (\tau \cdot \sigma) \cdot \pi(\gamma, a)$, as desired.
\end{proof}

\begin{rem}

Let $P_2$ denote the type-definable set of realizations of $tp(a_2 / C)$. Then $P_2$ embeds definably into $X$, via the map $a \mapsto (1, a) / E$. Moreover, the action of $\Gamma$ on $X$ extends the generic action of $K$ on $P_2$.
\end{rem}

\begin{prop}\label{prop_properties_of_Gamma_and_X}

\begin{enumerate}
\item The action of $\Gamma$ on $X$ is transitive.

\item The (image of the) type $tp(a_2 / C)$ is generic in the space $X$, which is connected.

\item The action of $\Gamma$ on $X$ is faithful.

\item Under the hypothesis $(R)$ of Theorem \ref{theo_config_groupe_faible}, the action is almost free: the stabilizers are finite.

\end{enumerate}

\end{prop}

\begin{proof}
1. We start with the following \begin{claim} Let $a, a^{'}$ realize $tp(a_2 / C)$. Then, there exists $\sigma \in \Gamma$ such that $(\sigma, a) \, E \, (1, a^{'})$. \end{claim} \begin{proof} Given such $a, a^{'}$, let $a^{''}$ realize $tp(a_2 / C) |_{C a a^{'}}$. Then, by Lemma \ref{lemma_action_reguliere_de_K_cas_gen_stable}, there exist $\tau_1, \tau_2$ such that $\tau_1 \models K|_{Ca} \cup K|_{Ca^{''}}$, $\tau_2 \models K|_{Ca^{'}} \cup K|_{Ca^{''}}$, $\tau_1(a) = a^{''}$ and $\tau_2(a^{''}) = a^{'}$. Let $\sigma$ be $\tau_2 \cdot \tau_1 \in \Gamma$. Since we already know that $\Gamma$ acts on $X$, it suffices to note that $\tau_1$ sends the class $(1,a) / E$ to $(1, a^{''}) / E$, which is then sent by $\tau_2$ to $(1, a^{'}) / E$, as desired. \end{proof}

Then, let $(\gamma, a) / E$ be an arbitray element of $X$. By the claim, let $\sigma \in \Gamma$ such that $(\sigma, a) E (1, a_2)$. Then, $\gamma \cdot \sigma^{-1}(1, a_2) = \gamma \cdot \sigma^{-1}((\sigma, a) / E) = (\gamma, a) / E$. So, we have proved that any element is in the orbit of $(1, a_2) / E$, which shows transitivity.

%Let us check that $(\sigma, a) \, E \, (1, a^{'})$ holds. Let $\tau$ realize $K|_{C a a^{'} a^{''} \tau_1 \tau_2}$. We compute that $\tau \cdot \sigma (a)= \tau \cdot \tau_2 ( \tau_1 (a) )= \tau (\tau_2 (a^{''})) = \tau(a^{'})$. The first equality holds because $\tau_1 \models K|_{Ca}$ and $\tau \cdot \tau_2 \models K|_{C \tau_1 (a)}$. The second and third ones are proved similarly. Thus, we have found an element $\tau$ realizing $K|_{C a a^{'} \sigma }$ such that $\tau \cdot \sigma (a) = \tau (a^{'})$. By definition, this shows that $(\sigma, a)$ and $(1, a^{'})$ are in the same $E$-class, as desired.

2. Let us show that the class of $(1, a_2)$ is generic in $X$. 

\begin{claim}
The stabilizer of the type $q=tp([1, a_2]_E / C)$ contains $K|_C$.
\end{claim}

\begin{proof}
Let $k \models K|_C$, let $\alpha_2 \models tp(a_2/C)|_{Ck}$, and $c=(1, \alpha_2)/E$, so that $c$ realizes $q|_{Ck}$. Then, by Lemma \ref{lemma_action_reguliere_de_K_cas_gen_stable}, we have $k(\alpha_2) \models  tp(a_2/C)|_{Ck}$. On the other hand, by Lemma \ref{lemma_k_a_E_1_ka}, we have $(1, k(a_2)) E (k, a_2)$. Thus, we have $k(c) = (1, k(a_2))/ E \models q|_{Ck}$, so $k \in Stab(q)$, as required.
\end{proof}

Moreover, the stabilizer of $q$ is a $C$-type-definable subgroup of $\Gamma$. Since $K$ generates $\Gamma$, the stabilizer of $q$ is $\Gamma$ itself, which proves genericity of $q$ and connectedness of $X$.  

3. Let $g \in \Gamma$ be an element that acts trivially on $X$. Let us show that $g=1$. We know that there exist $k_1, k_2 \in K$ such that $g= k_1^{-1} \circ k_2$. Then, by the hypothesis on $g$, we deduce that, for all $a$ realizing $tp(a_2 / C)|_{C k_1 k_2}$, we have $k_1(a) = k_2(a)$. Thus, by definition of a germ, we get $k_1 = k_2$, which implies $g= 1$.

4. Now, let us work under the hypothesis $(R)$ of Theorem \ref{theo_config_groupe_faible}. By transitivity of the action, it suffices to show that the stabilizer of the $E$-class of $(1, a_2)$ is finite. Let $\sigma \in \Gamma$ such that $\sigma \cdot (1, a_2) \, E \, (1, a_2)$. By definition, there exists $\tau$ realizing $K|_{C \sigma a_2}$ such that $\tau \cdot \sigma (a_2) = \tau(a_2)$. 

\begin{claim}The element $\tau(a_2)$ realizes $tp(a_2 / C)|_{C a_2}$.\end{claim}

\begin{proof}
By stationarity, it suffices to show that $\tau(a_2) \equiv_C a_2$ and $\tau(a_2) \downfree_C a_2$. Since the type $K|_C \otimes tp(a_2 / C) = tp(\tau, a_2 / C)$ is complete, Lemma \ref{lemma_action_reguliere_de_K_cas_gen_stable} applied backwards implies that $a_2 \downfree_C \tau(a_2)$. Then, by symmetry, we have  $\tau(a_2) \downfree_C a_2$, as required.
\end{proof}

 Moreover, since $\tau \models K|_{C \sigma a_2}$, and $\sigma \in Stab^r(K) = \Gamma$ (by Lemma \ref{lemma_gen_stable_generics_in_a_group}(2)), we know that $\tau \cdot \sigma$ realizes $K|_{C \sigma a_2}$. Thus, by the second point of Lemma \ref{lemma_action_reguliere_de_K_cas_gen_stable}, we deduce that $\tau \cdot \sigma \in acl(C a_2 \tau(a_2))$. Then, $\sigma \in acl( C a_2 \tau)$. However, the element $\tau$ satisfies $\tau \downfree_C a_2 \sigma$, so we have $\tau \downfree_{C a_2} \sigma$. Then, by Fact \ref{fact_forking} (5), we have $\sigma \in acl(C a_2)$. Finally, by compactness, the stabilizer of the $E$-class of $(1, a_2)$ is finite, as desired. \end{proof}

\end{subsection}

\begin{subsection}{End of the proof}

Here, we return to the context of Theorem \ref{theo_config_groupe_faible}.

\begin{proof}[Proof of Theorem \ref{theo_config_groupe_faible}]

By Proposition \ref{prop_WMA_dcl}, there are elements $b^{'}_1, b^{''}_2, b^{'''}_3 \in M$ and a configuration $(a_1, a_2, a_3, b_1, b_2, b_3)$ equivalent over $M$ to $(a^0_1, a^0_2, a^0_3, b^0_1, b^0_2, b^0_3)$, such that \begin{enumerate}
\item The tuple $  b^{'}_1 b^{''}_2 b^{'''}_3$ realizes $((p_0)_{y_1, y_2 , y_3})|_{C_0}$. 
\item The type $tp(a_1, a_2, a_3, b_1, b_2, b_3 / M)$ is generically stable over $C := acl(C_0  b^{'}_1 b^{''}_2 b^{'''}_3)$.

\item We have $a_1 \in dcl(C a_2 b_3)$, $a_2 \in dcl(C a_3 b_1)$, and $a_3 \in dcl(C a_1 b_2) \cap dcl(C a_2 b_1)$.
\end{enumerate}

Then, by the results proved above in Section 3, namely Proposition \ref{prop_properties_of_Gamma_and_X}, there is a connected $C$-type-definable group $\Gamma$ with a (unique) generically stable generic $K$, and a $C$-type-definable set $X$ equipped with a transitive and faithful $C$-definable action of $\Gamma$, such that the $C$-type-definable set of realizations of $tp(a_2 / C)$ embeds $C$-definably into $X$. Moreover, if $(R)$ holds, then the stabilizers for this action are finite.

We shall now construct, in the general case, a definable group configuration equivalent over $M$ to the initial one. At the end of the proof, we will explain how to deal with the case where $(R)$ holds. To build a definable group configuration equivalent over $M$ to $(a^0_1, a^0_2, a^0_3, b^0_1, b^0_2, b^0_3)$, it suffices to build one equivalent over $M$ to $(a_1, a_2, a_3, b_1, b_2, b_3)$. Let $b^{1}_1, b^{1}_2, b^{1}_3 \in M$ be such that $b^{1}_1 b^{1}_2 b^{1}_3 a_3 \equiv_C b_1 b_2 b_3 a_3$. 
Consider the following quadrangle:

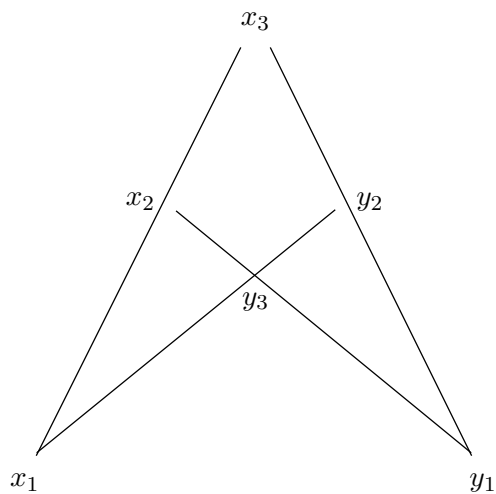
\begin{center}
    \begin{tikzcd}

&& x_3 \arrow[llddddd, dash, shorten= 1 mm] \arrow[rrddddd, dash, shorten= 1 mm]

\\
 
\\
 & x_2 & & y_2
\\
&& y_3

\\

\\

x_1 \arrow[rrruuu, dash, shift left = 1.0 ex]&&&& y_1\arrow[llluuu, dash, shift right = 1.0 ex]

    \end{tikzcd}
\end{center}

with the following definitions:

\begin{itemize}
    \item $x_1 = [f_{b^{1}_1}^{-1}] \circ [f_{b_1}]$
    \item $x_2 = [f_{b^{1}_1}^{-1}] \circ [g_{b^{1}_2}^{-1}] \circ [g_{b_2}^{}]\circ [f_{b^{1}_1}^{}]$
    \item $x_3 = x_2 \circ x_1^{} = [f_{b^{1}_1}^{-1}] \circ [g_{b^{1}_2}^{-1}] \circ [g_{b_2}^{}]\circ [f_{b_1}^{}]$
     \item $y_1 = x_3(a_2)$ 
     \item $y_2 = a_2$
  
    \item $y_3 = x_1(y_2)  $.
\end{itemize}

There are several facts to check, in order to make sure this is well-defined and equivalent to the original quadrangle. Note that, by Proposition \ref{prop_quadrangle_equivalent}, proving the equivalence with the original quadrangle yields that the quadrangle is a generically stable group configuration over $M$.

\begin{claim}

For $i=1,2,3$, we have $acl(M x_i) = acl(M b_i)$.

\end{claim}

\begin{proof}
For $i=1,2$, this is a consequence of the definition of definable group configurations, and of Proposition \ref{prop_description_base_canonique_af_c(a)}, applied to $a= a_2$ and $c=x_i$. For instance, the element $x_1$ is interalgebraic over $M$ with the germ $[f_{b_1}]$, which is, by Proposition \ref{prop_description_base_canonique_af_c(a)},  interdefinable over $C \subseteq M$ with the canonical basis $Cb(a_2 a_3 / acl(C b_1))$, which is, by point 4 of Definition \ref{defi_configuration_groupe}, interalgebraic over $C$ with $b_1$. For $i=3$, we use Lemma \ref{h_b1b_2_algebrique_sur_b_3_cas_gen_stable} and Proposition \ref{prop_description_base_canonique_af_c(a)}.
\end{proof}

\begin{claim}

The elements $x_1, x_2, x_3$ realize $K|_M$, and we have $x_1 \downfree_M x_2$.

\end{claim}

\begin{proof}
Since $K$ is generically stable, it suffices by stationarity (Proposition \ref{proprietes_type_gen_stable} (6)) to check that these elements realize $K|_C$, and that each is independent from $M$ over $C$. In fact, since $K$ is the unique generic of $\Gamma$, it suffices to check it for $x_1$ and $x_2$, and to prove the independence $x_1 \downfree_M x_2$. First note that, by definition of $K$, we have $tp(x_1 / C) = K|_C$.

For $x_1 = [f_{b^{1}_1}^{-1}] \circ [f_{b_1}]$, we know that $b_1 \downfree_C M$ and $b^{1}_1 \in M$, so that $x_1 \downfree_{C [f_{b^{1}_1}]} M$. Then, applying Lemma \ref{composee_germs_indep_cas_gen_stable}, we have $x_1 \downfree_C [f_{b^{1}_1}]$. Since $tp(x_1 / C) = K|_C$ is generically stable, we can apply transitivity, to get $x_1 \downfree_C M$, as desired.

For $x_2=  [f_{b^{1}_1}^{-1}] \circ [g_{b^{1}_2}^{-1}] \circ [g_{b_2}^{}]\circ [f_{b^{1}_1}^{}]$, the ideas are similar: one easily proves that $x_2 \downfree_{C [f_{b^{1}_1}]  [g_{b^{1}_2}]} M$. Also, we have by construction $[f_{b^{1}_1}] \flex  [g_{b^{1}_2}] \models (F \otimes G)|_C$ and $b_1 \downfree_C M$, so $b_1 \downfree_C [f_{b^{1}_1}] \flex  [g_{b^{1}_2}]$, thus $[f_{b^{}_1}] \downfree_C [f_{b^{1}_1}] \flex  [g_{b^{1}_2}]$. Therefore, by stationarity of $F$, we have $[f_{b^{}_1}] \models F|_{C [f_{b^{1}_1}] \flex  [g_{b^{1}_2}]}$. Hence, by Lemma \ref{gf=g2f2_cas_gen_stable}, there exists an $f \in F$ such that $ [g_{b_2}^{}]\circ [f_{b^{1}_1}^{}] =  [g_{b^{1}_2}^{}]\circ f$ and $f \downfree_C  [f_{b^{1}_1}]  [g_{b^{1}_2}]$, which implies that $x_2 \models K|_C$. Then, by Lemma \ref{composee_germs_indep_cas_gen_stable}, we also have $ [f_{b^{1}_1}^{-1}] \circ f \downfree_C [f_{b^{1}_1}]  [g_{b^{1}_2}]$, i.e. $x_2 \downfree_C  [f_{b^{1}_1}]  [g_{b^{1}_2}]$. We can then apply transitivity, just as before.
\end{proof}

\begin{claim}

The elements $y_1=x_3(y_2)$ and $y_3= x_1(y_2)$ are well-defined, and satisfy the following: $y_1 =    f_{b^{1}_1}^{-1} \circ g_{b^{1}_2}^{-1} (a_1)$ and $y_3 = f_{b^{1}_1}^{-1}(a_3)$.

\end{claim}

\begin{proof}

Let us first show that the elements are well-defined. Since $tp(a_2 / C)$ is generically stable, it suffices to check that $y_2 \downfree_C x_3$ and $y_2 \downfree_C x_1$. These verifications rely on the facts that $y_2 = a_2 \downfree_C M b_1 b_2$ and $b^1_1 b^1_2 \in M$, and  are left to the reader.

For the equalities, we have by definition that $b_1 b_2 \downfree_C a_2$ and $g_{b_2} \circ f_{b_1}(y_2) = a_1 \downfree_C  b^{1}_1 b^{1}_2$, so that $y_1 = x_3(y_2) =  [f_{b^{1}_1}^{-1}] \circ [g_{b^{1}_2}^{-1}] \circ [g_{b_2}^{}]\circ [f_{b_1}^{}] (a_2) = [f_{b^{1}_1}^{-1}] \circ [g_{b^{1}_2}^{-1}] (a_1) =  f_{b^{1}_1}^{-1} \circ g_{b^{1}_2}^{-1} (a_1)$.

Similarly, we have $a_2 \downfree_C b_1 b^{1}_1$, and $f_{b_1}(a_2) = a_3 \downfree_C b^{1}_1$, which implies the following: $y_3 = x_1(y_2) =  [f_{b^{1}_1}^{-1}] \circ [f_{b_1}](a_2) =  [f_{b^{1}_1}^{-1}] (a_3) =  f_{b^{1}_1}^{-1} (a_3).$ 
\end{proof}

\begin{claim}
For $i=1,2,3$, we have $acl(M y_i) = acl(M a_i)$.
\end{claim}

\begin{proof}
This is easily deduced from the identities proved above, and the fact that $y_2 = a_2$.\end{proof}

So, all that remains is the case where $(R)$ holds. Let $a \in X(M)$, and $N \leq \Gamma(M)$ be the finite stabilizer of $a$ for the action of $\Gamma$. Since the action is transitive, we have a $\Gamma$-equivariant (relatively) $M$-definable bijection $\rho : \Gamma / N \simeq X$. In particular, there is a $\Gamma$-equivariant (relatively) $M$-definable finite-to-one surjection $\pi : \Gamma \rightarrow X$.

\begin{claim}
We have $\pi_*(K|_C) = tp(a_2 / C)$. 
\end{claim}

\begin{proof}
Since $\pi$ is $\Gamma$-equivariant, this follows from Lemma \ref{lemma_pb_is_generic}, and uniqueness of the generic of $X$ (see Proposition \ref{prop_gen_stable_generics_in_a_space} (1)).
\end{proof}

So, let $g \models K|_C$ be such that $\pi(g) = a_2 = y_2$. In the previous configuration, replace $y_1$ with $x_3 \cdot g$, $y_2$ with $g$, and $y_3$ with $x_1 \cdot g$. Since $\pi$ is equivariant, it is straightforward to compute that $\pi(x_3 \cdot g) = x_3(a_2) = y_1$, $\pi(g) = y_2$, and $ \pi(x_1 \cdot g) = x_1(a_2) = y_3$. Since $\pi$ is $M$-definable and has finite fibers, this shows that the quadrangle $(x_3 \cdot g, g, x_1 \cdot g, x_1, x_2, x_3)$ is equivalent over $M$ to $(y_1, y_2, y_3, x_1, x_2, x_3)$, which concludes the proof. 
\end{proof}

\end{subsection}

\end{section}

\printbibliography[
heading=bibintoc,
title={References}
]

%%% Seuls les elements cités quelque part apparaissent in the biblio 

Paul Z. WANG, 

Ecole Normale Supérieure de Paris - PSL,

 45 rue d'Ulm, 75005 Paris,
 
 paul.wang@ens.psl.eu

\end{document}